\documentclass[12pt] {article}
\usepackage{lscape}
\usepackage{color}
\usepackage{amsmath,amsthm}
\usepackage{amssymb,latexsym}
\usepackage{comment}

\usepackage{graphicx}

\usepackage{amssymb}
\usepackage{amsfonts}
\usepackage{epic,eepic}

\usepackage{empheq}
\usepackage{amscd}
\usepackage{chngcntr}
\usepackage{enumitem}
\newlist{paragraphlist}{enumerate}{1}

\setlist[paragraphlist,1]{leftmargin=*,label={\bfseries \arabic*}}

\counterwithin{paragraphlisti}{subsubsection}

\title{Octonionic Calabi-Yau theorem.}
\date{}

\author{Semyon Alesker\footnote{Partially supported by ISF grants 865/16 and 743/22 and the  US - Israel BSF grant 2018115.}\\
{ \normalsize Department of Mathematics}\\
{ \normalsize Tel Aviv University, Ramat Aviv}\\
{ \normalsize 69978 Tel Aviv, Israel}\\
{ \normalsize e-mail: semyon@post.tau.ac.il}
\and
 Peter V. Gordon\footnote{Supported in a part by  the  US - Israel BSF grant 2020005 and Simons Foundation Collaboration Grant for Mathematicians 317882. }\\
{ \normalsize  Department of Mathematical Sciences}\\
{ \normalsize  Kent State University}\\
{ \normalsize  Kent, Ohio 44242, USA}\\
{ \normalsize  e-mail: gordon@math.kent.edu}}

\def\RR{\mathbb{R}}
\def\CC{\mathbb{C}}
\def\QQ{\mathbb{Q}}

\def\HH{\mathbb{H}}

\def\PP{\mathbb{P}}

\def\OO{\mathbb{O}}
\def\SS{\mathbb{S}}

\def\eps{\varepsilon}
\def\alp{\alpha}

\def\Ome{\Omega}
\def\lam{\lambda}

\def\to{\longrightarrow}

\def\qed { Q.E.D. }

\def\inj{\hookrightarrow}

\swapnumbers
\newtheorem{theorem}{Theorem}[section]
\newtheorem{corollary}[theorem]{Corollary}
\newtheorem{lemma}[theorem]{Lemma}
\newtheorem{proposition}[theorem]{Proposition}
\newtheorem{claim}[theorem]{Claim}
\theoremstyle{definition}

\newtheorem{example}[theorem]{Example}
\newtheorem{definition}[theorem]{Definition}
\newtheorem{remark}[theorem]{Remark}
\theoremstyle{conjecture}
\newtheorem{conjecture}[theorem]{Conjecture}
\theoremstyle{principle}

  \def\cc{{\cal C}}
\def\cd{{\cal D}}  
 \def\ch{{\cal H}} 
 \def\ck{{\cal K}} 
\def\cm{{\cal M}}  \def\co{{\cal O}}

\def\cs{{\cal S}} \def\ct{{\cal T}} \def\cu{{\cal U}}
\def\cv{{\cal V}}  \def\cx{{\cal X}}
 
\def\cm{{\cal M}}
\def\pt{\partial}

\numberwithin{equation}{section}

\def\dfq{\frac{\partial ^2 f}{\partial\bar q_i \partial q_j}}

\def\ddbq{\frac{\partial ^2 b}{\partial \bar q_i \partial q_j}}

\def\uch{\underline{\ch_2}(M)}

\begin{document}
\maketitle

\begin{abstract}
On a certain class of 16-dimensional manifolds a new class of Riemannian metrics, called octonionic K\"ahler, is introduced and studied. It is  an octonionic analogue of K\"ahler metrics on complex manifolds and
of HKT-metrics of hypercomplex manifolds. Then for this class of metrics an octonionic version of the Monge-Amp\`ere equation is introduced and solved under appropriate assumptions. The latter result
is an octonionic version of the Calabi-Yau theorem from K\"ahler geometry.
\end{abstract}

\def\bfm{{\cal B} (\phi, f)}
\def\bo{\partial \Omega}
\def\vz{v_{\zeta}}
\def\fss{F^{**}}
\def\ggo{\Gamma_0}
\def\oe{\omega(\eps)}
\def\xom{\chi_{\Omega}}

\tableofcontents

\section{Introduction.}\label{S:introduction}
\begin{paragraphlist}
\item Real and complex Monge-Amp\`ere (MA) equations play a central role in several areas of analysis and geometry.
In particular, the celebrated Calabi-Yau theorem \cite{yau} establishes existence of Ricci flat K\"ahler metrics on compact K\"ahler manifolds reducing the problem to solving
a complex MA equation. This result is a powerful analytic tool to study geometry and topology of
compact K\"ahler manifolds.

\item There were several attempts
to generalize the notion of MA equation and of plurisubharmonic function to a few other contexts.
Namely, the Calabi-Yau theorem was successfully extended to the more general class of compact Hermitian manifolds by Tosatti-Weinkove \cite{tosatti-weinkove} in full generality and previously by Cherrier \cite{cherrier} in special cases.

Furthermore, the first author \cite{alesker-bsm-03}
 introduced and studied the quaternionic MA operator and the class of quaternionic plurisubharmonic (qpsh) functions on the flat quaternionic space $\HH^n$.
 Later he introduced and solved, under appropriate assumptions, the Dirichlet
 problem for quaternionic MA equation in \cite{alesker-jga-03}. The latter problem was solved under weaker assumptions by Zhu \cite{zhu} and Kolodziej-Sroka \cite{sroka-kolodziej}.

 M. Verbitsky and the first author \cite{alesker-verbitsky-1} extended the definition of the quaternionic MA operator and the class of qpsh functions to the broader class of
 hypercomplex manifold; then the first author \cite{alesker-quatern-mflds} extended them to even broader class of quaternionic manifolds. M. Verbitsky and the first author
 \cite{alesker-verbitsky-2} introduced a quaternionic version of the Calabi problem. So far this problem is open in full generality, but it was proven in a few special cases, see \cite{alesker-calabi-yau}, \cite{getilli-vezzoni-calabi-yau}, \cite{dinew-sroka}.
 Some relevant a priori estimates were obtained in \cite{alesker-verbitsky-2}, \cite{alesker-shelukhin-1}, \cite{alesker-shelukhin-2}, \cite{sroka-C0-estimate}.

 A number of results in quaternionic pluripotential theory were obtained in \cite{wan-zhang-2015}, \cite{wan-kang-2017}, \cite{wan-wang-2017}.

\item The MA operator and a class of plurisubharmonic functions in two octonionic variables were introduced by the first author in \cite{alesker-octon}, where also a few properties of them were established and
an application to valuations on convex sets was given.

 \item At the same time there is a parallel development by Harvey and Lawson, e.g. \cite{harvey-lawson-2009a}, \cite{harvey-lawson-2009b}, \cite{harvey-lawson-2011}, \cite{harvey-lawson-2019}.
 The notion of plurisubharmonic function is developed in  great generality, as well as the theory of homogeneous MA equation. Being specialized to the quaternionic situation, this theory has an overlap
 with some of the results mentioned above.

 \item Let us describe the main results of this paper.

 As mentioned above, the MA operator in two octonionic variables was introduced by the first author in \cite{alesker-octon}.
 Based on this operator, the present paper introduces and studies the octonionic MA equation analogous to the Calabi problem in K\"ahler geometry.

The MA equation is stated for a special class of 16-dimensional manifolds with a kind of octonionic structure. We call such manifolds $GL_2(\OO)$-affine manifolds, they are introduced in Section \ref{S:G-mflds}. Basic examples of such manifolds
are torii $\OO^2/\Lambda$, where $\OO^2\simeq \RR^{16}$ is the octonionic plane, and $\Lambda\subset \OO^2$ is a lattice.

Furthermore in Section \ref{S:metrics-octon} we introduce on any $GL_2(\OO)$-affine manifold a class of Riemannian metrics which are octonionic analogues of K\"ahler metrics. We call them octonionic K\"ahler metrics.
In local coordinates, i.e. on $\OO^2$, such metrics are characterized by two conditions.
First, the restriction of each such metric on any (right) octonionic line is proportional to the standard Euclidean metric on $\OO^2\simeq \RR^{16}$.
Second, such metrics satisfy a system of linear first order PDEs (see (\ref{E:system-equations})) analogous to the closeness of K\"ahler form on a complex manifold.
We also prove an analogue of the local $dd^c$-lemma from complex analysis (see Theorem \ref{T:octon-kahler}): locally any such metric can be given by a smooth potential $u$, i.e. the metric can be identified with the octonionic Hessian of $u$.
The converse is also true provided the Hessian of $u$ is pointwise positive definite.

Next we introduce an octonionic MA equation on a $GL_2(\OO)$-affine manifold $M$. Let $G_0$ be an octonionic K\"ahler ($C^\infty$-smooth)  metric; in local coordinates it can be identified with a positive definite octonionic Hermitian matrix.
Let $f\colon M\to \RR$ be a $C^\infty$-smooth  function. Consider the equation which in local coordinates reads
\begin{eqnarray}\label{E:MA-main-equation}
\det(G_0+Hess_\OO(\phi))=e^f\det(G_0),
\end{eqnarray}
where $Hess_\OO$ is the octonionic Hessian (see Section \ref{S:oct-hessian} for the precise definition), $\det$ is the determinant of octonionic Hermitian $2\times 2$ matrices (see (\ref{D:det})).
The main result of this paper, Theorem \ref{T:main-result}, says that on the subclass of compact connected so called $Spin(9)$-affine manifolds
(example: tori $\OO^2/\Lambda$ where $\Lambda \subset\OO^2$ is a lattice) there exists a $C^\infty$-smooth solution $\phi$ of (\ref{E:MA-main-equation})
if and only if the function $f$ satisfies the normalization condition
$$\int_M (e^f-1)\det(G_0)\cdot \Theta=0,$$
where $\Theta$ is a parallel non-vanishing 3/4-density on $M$ (it always exists on $Spin(9)$-affine manifolds and is unique up to multiplication by a constant).

\item The class of metrics on $GL_2(\OO)$-affine manifolds we mentioned is analogous not only to K\"ahler metrics on complex manifolds, but also to Hessian metrics on affine manifolds \cite{cheng-yau}, and to HKT-metrics on hypercomplex manifolds \cite{howe-papadopoulos}, \cite{grantcharov-poon}.
Our main result on the solvability of the MA equation is analogous to and is motivated by the Calabi-Yau theorem \cite{yau} for K\"ahler manifolds. Also it is analogous to the Cheng-Yau theorem \cite{cheng-yau} for Hessian manifolds,
and to the similar conjecture  on quaternionic MA equation on HKT-manifolds
due to M. Verbitsky and the author \cite{alesker-verbitsky-2} which has been proven so far in a few special cases mentioned above.

\item In literature there are quite a few investigations of geometric structures related to octonions.
Thus Gibbons,  Papadopoulos, and  Stelle \cite{gibbons-papadopoulos-stelle} introduced a class of Octonionic K\"ahler with Torsion metrics on $8k$-dimensional manifolds in connection with black hole moduli spaces.

 Friedrich \cite{friedrich} introduced the notion of nearly parallel $Spin(9)$-structure on 16-manifolds whose structure group is reduced to $Spin(9)$.
The role of $Spin(9)$ in octonionic geometry is reviewed by Parton and Piccinni \cite{parton-piccinni} (see also \cite{kotrbaty}).

Very recently Kotrbat\'y and Wannerer \cite{kotrbaty-wannerer} described the multiplicative structure on $Spin(9)$-invariant valuations on convex sets on $\OO^2$ and gave applications to integral geometry.

\item {\bf Acknowledgements.} We thank F. Nazarov and M. Sodin for useful discussions and G. Grantcharov, G. Papadopoulos, I. Soprunov, and A. Swann for supplying us with some references.
We are grateful to the anonymous referee for very careful reading the paper and numerous remarks. 

\end{paragraphlist}

\section{Algebra of octonions.}\label{S:properties-of-octonions}

\begin{paragraphlist}
\item The octonions $\OO$ form an 8-dimensional algebra over the reals
$\RR$ which is neither associative nor commutative. The product on $\OO$ can
be described as follows. $\OO$ has a basis $e_0,e_1,\dots,e_7$ over
$\RR$ where $e_0=1$ is the unit, and the product of basis elements is given by the
following multiplication table:

\begin{center}
$$\begin{array}{c|c|c|c|c|c|c|c|}
   & e_1& e_2& e_3& e_4& e_5& e_6& e_7\\\hline
e_1&-1  &e_4 &e_7 &-e_2&e_6 &-e_5&-e_3\\\hline e_2&-e_4&-1
&e_5&e_1 &-e_3 &e_7&-e_6\\\hline
e_3&-e_7&-e_5&-1&e_6&e_2&-e_4&e_1\\\hline
e_4&e_2&-e_1&-e_6&-1&e_7&e_3&-e_5\\\hline
e_5&-e_6&e_3&-e_2&-e_7&-1&e_1&e_4\\\hline
e_6&e_5&-e_7&e_4&-e_3&-e_1&-1&e_2\\\hline
e_7&e_3&e_6&-e_1&e_5&-e_4&-e_2&-1\\\hline
\end{array}$$
\end{center}

\item There is also an easier alternative way to reproduce the product rule using the so
called {\itshape Fano plane} (see Figure 1 below). In Figure 1 each
pair of distinct points lies in a unique line (the circle is also
considered to be a line). Each line contains exactly three points,
and these points are cyclically oriented. If $e_i,e_j,e_k$ are
cyclically oriented in this way then
$$e_ie_j=-e_je_i=e_k.$$
We have to add two more rules:

$\bullet$ $e_0=1$ is the identity element;

$\bullet$ $e_i^2=-1$ for $i>0$.

All these rules define uniquely the algebra structure of $\OO$. The
center of $\OO$ is equal to $\RR$.



\begin{figure}
\centering \includegraphics[width=4in]{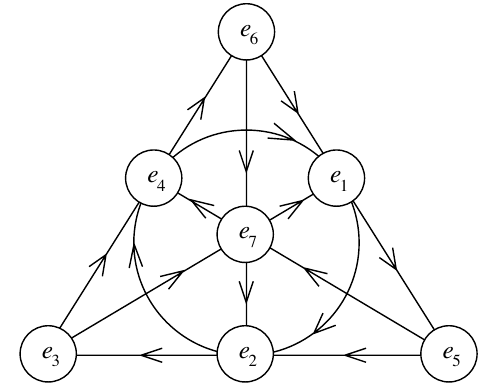}
 \caption{
 Fano plane, the figure is taken from \cite{baez}}
\label{figure:1}
\end{figure}

\item Every octonion $q\in \OO$ can be written uniquely in the form
$$q=\sum_{i=0}^7x_ie_i$$
where $x_i\in \RR$. The summand $x_0e_0=x_0$ is called the real
part of $q$ and is denoted by $Re(q)$.

One defines the octonionic conjugate of $q$ by
$$\bar q:=x_0-\sum_{i=1}^7x_ie_i.$$
It is well known that the conjugation is an anti-involution of
$\OO$:
$$\overline{\overline{q}}=q,\, \overline{a+b}=\bar a+\bar b,\, \overline{ab}=\bar b\bar a.$$
Let us define a norm on $\OO$ by $|q|:=\sqrt{q\bar q}$. Then
$|\cdot|$ is a multiplicative norm on $\OO$: $|ab|=|a||b|$. The
square of the norm $|\cdot|^2$ is a positive definite quadratic
form. Its polarization is a  positive definite scalar product
$<\cdot,\cdot>$ on $\OO$ which is given explicitly by
$$<x,y>=Re(x\bar y).$$

 Furthermore $\OO$ is a division algebra: any $q\ne 0$ has a unique
inverse $q^{-1}$ such that $qq^{-1}=q^{-1}q=1$. In fact
$$q^{-1}=|q|^{-2}\bar q.$$

\item We denote by $\HH$ the usual quaternions. It is associative division
algebra. We will fix once and for all an imbedding of algebras
$\HH\subset \OO$. Let us denote by $i,j\in \HH$ the usual
quaternionic units, and $k=ij$. Then $i,j,k$ are pairwise orthogonal
with respect to the scalar product $<\cdot,\cdot>$. Let us fix once
and for all an octonionic unit $l\in \OO, l^2=-1$ which is
orthogonal to $i,j,k$. Then $l$ anti-commutes with $i,j,k$. Every
element $q\in \OO$ can be written uniquely in the form
$$q=x+yl \mbox{ with } x,y\in \HH.$$
The multiplication of two such octonions is given by the formula
\begin{eqnarray}\label{E:product-oct}
(x+yl)(w+zl)=(xw-\bar z y)+(zx+y\bar w)l
\end{eqnarray}
where $x,y,w,z\in\HH$.

\item We have the following weak forms of the associativity in
octonions. A proof can be found in \cite{alesker-octon}, Lemma 1.1.1.
\begin{lemma}\label{L:oalg1}
Let $a,b,c\in \OO$. Then

(i) $Re((ab)c)=Re(a(bc))$ (this real number will be denoted by
$Re(abc)$).

(ii) $a(bc)+\bar b(\bar a c)=(ab+\bar b\bar a)c$.

(iii) $(ca)b+(c\bar b)\bar a=c(ab+\bar b \bar a)$.

(iv) Any subalgebra of $\OO$ generated by any two elements and
their conjugates is associative. It is always isomorphic either to
$\RR$, $\CC$, or $\HH$.

(v) $Re((\bar ab)( ca))=|a|^2Re(bc).$
\end{lemma}

The following identities are called Moufang identities (see e.g. \cite{schafer}).
\begin{lemma}\label{L:moufang}
Let $x,y,z\in \OO$. Then
\begin{eqnarray*}
(zx)(yz)=z((xy)z)=(z(xy))z.
\end{eqnarray*}
\end{lemma}
Polarizing the Moufang identities we immediately get a corollary.
\begin{corollary}\label{Cor-moufang}
Let $x,y,z_1,z_2\in \OO$. Then
\begin{eqnarray*}
( z_1x)(yz_2)+( z_2x)(yz_1)= z_1((xy)z_2)+ z_2((xy)z_1)=( z_1(xy))z_2+( z_2(xy))z_1.
\end{eqnarray*}
\end{corollary}

\end{paragraphlist}

\section{Octonionic projective line $\OO\PP^1$.}\label{S:oct-Hessian}
\begin{paragraphlist}
\item The results of this section are well known and elementary, we provide them for completeness of presentation.

\begin{definition}\label{D:op1}
An octonionic (right) line in $\OO^2$ is a real 8-dimensional subspace $L\subset\OO^2$ either of the form
$\{(q,aq)|q\in\OO\}$ with arbitrary $a\in \OO$,
or of the form $\{(0,q)|q\in\OO\}$.
\end{definition}
The latter option is considered as the line at infinity. The set of all octonionic lines is denoted by $\OO\PP^1$. It is easy to see that it is a compact subset of the Grassmannian of 8-dimensional linear subspaces of $\RR^{16}$.
Moreover it is a smooth submanifold diffeomorphic to the 8-sphere $\SS^8$: the line at infinity corresponds to the pole of this sphere, and the rest is parameterized by octonions $a\in\OO$ is diffeomorphic to $\RR^{8}$.

\item Let us summarize a few other basic and well known properties of $\OO\PP^1$.
\begin{lemma}\label{L:properties-op1}
(i) For any $a\in \OO$ the set $\{(aq,q)|q\in \OO\}$ is also an octonionic line.
\newline
(ii) Any two octonionic lines either intersect trivially or coincide.
\newline
(iii) For any non-zero vector $\xi\in \OO^2$ there is a unique octonionic line containing $\xi$. We call it the line spanned by $\xi$.
\newline
(iv) Let $\xi,\eta\in \OO^2$ be two vectors of norm one. The lines spanned by these vectors are equal if and only if one has for the $2\times 2$-matrices
$$\xi\cdot\xi^*=\eta\cdot\eta^*.$$
In this case we write $\xi\sim \eta$.
\end{lemma}
The proof of this lemma is elementary and straightforward.

\end{paragraphlist}

\section{Octonionic Hessian.}\label{S:oct-hessian} In this section we remind the  octonionic Hessian introduced in \cite{alesker-octon}. For $i=1,2$ write
$$q_i=\sum_{a=0}^7 x_i^ae_a \mbox{ with } x_i^a\in\RR.$$
 Let $F\colon \OO^2\to \OO$ be a sufficiently smooth function. Define for $i,j=1,2$
\begin{eqnarray*}
\frac{\pt F}{\pt \bar q_i}:=\sum_{a=0}^7 e_a\frac{\pt F}{\pt x_i^a},\\
\frac{\pt F}{\pt q_j}:=\sum_{a=0}^7 \frac{\pt F}{\pt x_j^a}\bar e_a.
\end{eqnarray*}

In general the operators $\frac{\pt }{\pt \bar q_i}$ and $\frac{\pt }{\pt q_j}$ do not commute. However, they do commute when applied to real valued functions.

Let $f\colon \OO^2\to \RR$ be a $C^2$-smooth real valued function. Define its octonionic Hessian by
\begin{eqnarray}\label{D:oct-hessian}
Hess_\OO(f):=\left(\frac{\pt^2 f}{\pt \bar q_i \pt q_j}\right)_{i,j=1,2},
\end{eqnarray}
where the operators $\frac{\pt }{\pt \bar q_i}$ and $\frac{\pt }{\pt q_j}$ are applied in either order. Then $Hess_\OO(f)$ is an octonionic Hermitian $2\times 2$ matrix.

\section{Octonionic linear algebra.}\label{S:linear-algebra}
\begin{paragraphlist}

\item Let us denote
$$\OO^m:=\{(q_1,\dots,q_m)|\, q_i\in \OO\}.$$
For $\xi=(q_1,\dots,q_m)\in \OO^m$, $x\in \OO$ we will denote by
$\xi\cdot x$ the $m$-tuple $(q_1x,\dots,q_mx)\in\OO^m$. Notice that
usually we will write elements of $\OO^m$ as $m$-columns rather than
rows.

Let us denote by $\ch_n(\RR)$ the space of real symmetric
$(n\times n)$-matrices. The space $\ch_n(\RR)$ is naturally
identified with the space of real valued quadratic forms on
$\RR^n$.

Let us denote by $\ch_m(\OO)$ the space of octonionic Hermitian
$(m\times m)$-matrices. By definition, an $(m\times m)$-matrix
$A=(a_{ij})$ with octonionic entries is called Hermitian if
$a_{ij}=\overline{a_{ji}}$ for any $i,j$. For a matrix $A=(a_{ij})$
denote also $A^*:=(\overline{a_{ji}})$. In what follows we will be only interested
in $2\times 2$ octonionic Hermitian matrices. In this case we have the
following explicit description of such matrices. Namely,
\begin{eqnarray}
\ch_2(\OO)=\left\{ \left[\begin{array}{cc}
                          a     &q\\
                          \bar q&b
                         \end{array}\right]\big| a,b\in \RR, q\in \OO\right\}.
\end{eqnarray}

\item\label{item:lin-alg2} We have the natural $\RR$-linear map
\begin{eqnarray}\label{D:j}
j\colon\ch_2(\OO)\to \ch_{16}(\RR),
\end{eqnarray} which is defined
as follows: for any $A\in \ch_2(\OO)$ the value of the quadratic
form $j(A)$ on any octonionic 2-column $\xi\in \OO^2=
\RR^{16}$ is given by $j(A)(\xi)=Re(\xi^* A\xi)$ (note that the
bracketing inside the right hand side of the formula is not important due to Lemma
\ref{L:oalg1}(i)). It is easy to see that the map $j$ is
injective. Via this map $j$ we will identify $\ch_2(\OO)$ with a
subspace of $\ch_{16}(\RR)$.

\def\hr{\ch_{16}(\RR)}
\def\ho{\ch_2(\OO)}

\item Let us construct now a linear map
$$\theta\colon \ch_{16}(\RR)\to \ch_2(\OO)$$
in the opposite direction and such that
 $\theta\circ j=Id$; it will be useful later. For any
$B\in \ch_{16}(\RR)$ let us denote by $b$ the corresponding
quadratic form on $\RR^{16}= \OO^2$, namely
$$b(x)=\sum_{i,j=1}^{16}x_iB_{ij}x_j.$$

 Define
\begin{eqnarray}\label{Def:theta}
\theta(B):=\frac{1}{16}Hess_\OO(b)=\frac{1}{16}\left(\ddbq\right)_{i,j=1}^2.
\end{eqnarray}
Note that the matrix
in the right hand side of the last formula is independent of a point
in $\OO^2$.

The next lemma is proved in \cite{alesker-octon}, Corollary 1.2.2, but we present a proof for the sake of completeness.
\begin{lemma}\label{L:la2}
The following identity holds.
$$\theta\circ j=Id.$$
\end{lemma}
{\bf Proof.} Let $B=\left[\begin{array}{cc}
                            a&q\\
                            \bar q&b
                           \end{array}\right]$ where $a,b\in \RR, \, q\in \OO$. Below we denote $\Delta_i=\sum_{a=0}^7\frac{\pt^2}{\pt (x_i^a)^2}$. We have
\begin{eqnarray*}
(\theta\circ j)(B)=\frac{1}{16}\left( Hess_\OO(Re(\xi^*B\xi)\right)=\\
\frac{1}{16} Hess_\OO\left(a|\xi_1|^2+b|\xi_2|^2+2 Re(\bar \xi_1 q\xi_2)\right)=\\
\frac{1}{16}\left[\begin{array}{cc}
                 a \Delta_1(|\xi_1|^2)&2\frac{\pt^2}{\pt\bar\xi_1\pt\xi_2}Re(\bar\xi_1 q\xi_2)\\
                 2\frac{\pt^2}{\pt\bar\xi_2\pt\xi_1}Re(\bar\xi_1 q\xi_2)&b \Delta_2(|\xi_2|^2)
                \end{array}\right]=\\
\left[\begin{array}{cc}
        a&\frac{1}{8}\frac{\pt^2}{\pt\bar\xi_1\pt\xi_2}Re(\bar\xi_1 q\xi_2)\\
        \frac{1}{8}\frac{\pt^2}{\pt\bar\xi_2\pt\xi_1}Re(\bar\xi_1 q\xi_2)&b
        \end{array}\right].
\end{eqnarray*}
 We have to show that the term above the diagonal equals $q$ (the term under the diagonal is its conjugate). Indeed, we have
 \begin{eqnarray*}
 \frac{1}{8}\frac{\pt^2}{\pt\bar\xi_1\pt\xi_2}Re(\bar\xi_1 q\xi_2)=\sum_{a,b=0}^7e_a\left(\frac{\pt^2}{\pt x_1^a\pt x_2^b} Re(\bar\xi_1\cdot q\cdot\xi_2)\right)\bar e_b=\\
 \frac{1}{8}\sum_{a,b,c,d=0}^7 e_a\left(\frac{\pt^2}{\pt x_1^a\pt x_2^b} Re(x_1^c\bar e_c\cdot q\cdot x_2^d e_d)\right)\bar e_b=\\
 \frac{1}{8}\sum_{a,b,c,d=0}^7 e_a\left( Re(\delta_{ac}\bar e_c\cdot q\cdot \delta_{bd} e_d)\right)\bar e_b=\\
 \frac{1}{8}\sum_{a,b=0}^7e_a(Re(\bar e_a q e_b))\bar e_b=\frac{1}{8}\sum_b\left(\sum_a e_a(Re(\bar e_a(qe_b)))\right)\bar e_b=\\
 \frac{1}{8}\sum_b (qe_b)\bar e_b=\frac{1}{8}\sum_b q(e_b\bar e_b)=q.
 \end{eqnarray*}
 Lemma follows. \qed

\item There is a useful notion of positive definite octonionic Hermitian matrix.
\begin{definition}\label{D:positive-def}
Let $A\in\ch_2(\OO)$. $A$ is called {\itshape positive definite}
(resp. {\itshape non-negative definite}) if for any 2-column
$\xi\in\OO^2\backslash\{0\}$
$$Re(\xi^* A\xi)>0 \,\, (\mbox{resp. } Re(\xi^* A\xi)\geq 0).$$
\end{definition}

 For a positive definite (resp. non-negative definite) matrix $A$ one
writes as usual $A>0$ (resp. $A\geq 0$).

\item On the class of octonionic Hermitian $2\times 2$-matrices there is
a nice notion of determinant which is defined by
\begin{eqnarray}\label{D:det}
\det\left(\left[\begin{array}{cc}
                a&q\\
                \bar q&b
                \end{array}\right]\right):=ab-|q|^2.
\end{eqnarray}
\begin{remark}
It turns out that a nice notion of determinant does exist also on
octonionic Hermitian matrices of size 3, see e.g. Section 3.4 of
\cite{baez}. Note also that a nice notion of determinant does exist for
{\itshape quaternionic} Hermitian matrices of any size: see the
survey \cite{aslaksen}, the article \cite{gelfand-retakh-wilson},
and for applications to quaternionic plurisubharmonic functions see
\cite{alesker-bsm-03}, \cite{alesker-jga-03},
\cite{alesker-adv-05}, \cite{alesker-alg-anal}.
\end{remark}

\item The following result is a version of the Sylvester criterion for
octonionic matrices of size two, see \cite{alesker-octon}, Proposition 1.2.5.
\begin{proposition}\label{P:sylvester}
Let $A=\left[\begin{array}{cc}
              a&q\\
              \bar q&b
              \end{array}\right]\in\ch_2(\OO)$. Then $A>0$ if and
only if $a>0$ and $\det A>0$.
\end{proposition}

\item Now we remind the notion of {\itshape mixed determinant}
of two octonionic Hermitian matrices in analogy to the classical
real case (see e.g. \cite{schneider-book}). First observe that the
determinant $\det \ch_2(\OO)\to \RR$ is a homogeneous polynomial of
degree 2 on the real vector space $\ch_2(\OO)$. Hence it admits a
unique polarization: a bilinear symmetric map
$$D\colon \ch_2(\OO)\times \ch_2(\OO)\to \RR$$
such that $D(A,A)=\det A$ for any $A\in \ch_2(\OO)$. This map $D$ is
called the mixed determinant. By the abuse of notation it will be
denoted again by $\det$. Explicitly if $A=(a_{ij})_{i,j=1}^2$,
$B=(b_{ij})_{i,j=1}^2$ are octonionic Hermitian matrices then
\begin{eqnarray}\label{E:mixed-det}
\det(A,B)=\frac{1}{2}\left(a_{11}b_{22}+a_{22}b_{11}-2Re(a_{12}b_{21})\right).
\end{eqnarray}
\begin{lemma}\label{L:mix-det-posit}
If $A,B\in \ch_2(\OO)$ are positive (resp. non-negative) definite
then
$$\det(A,B)>0 \, \,(\mbox{resp. } \det(A,B)\geq 0).$$
\end{lemma}
{\bf Proof.} Let us assume that $A,B$ are positive definite. Then by
Proposition \ref{P:sylvester} we have
\begin{eqnarray*}
a_{11}>0,a_{22}>0, |a_{12}|<\sqrt{a_{11}a_{22}},\\
b_{11}>0,b_{22}>0, |b_{12}|<\sqrt{b_{11}b_{22}}.
\end{eqnarray*}
These inequalities imply
\begin{eqnarray}\label{E:iii}
Re(a_{12}b_{21})\leq |a_{12}|\cdot
|b_{21}|<\sqrt{a_{11}a_{22}b_{11}b_{22}}.
\end{eqnarray} Substituting
(\ref{E:iii}) into (\ref{E:mixed-det}) we get
\begin{eqnarray*}
\det(A,B)>\frac{1}{2}\left(a_{11}b_{22}+a_{22}b_{11}-
2\sqrt{a_{11}a_{22}b_{11}b_{22}}\right)\geq 0
\end{eqnarray*}
where the last estimate is the arithmetic-geometric mean inequality.
Thus $\det(A,B)>0$ for positive definite matrices $A,B$. For
non-negative definite matrices the result follows by approximtion. \qed

\begin{remark}
One can also prove the following version of the Aleksandrov
inequality for mixed determinants \cite{aleksandrov-mixed}: The
mixed determinant $\det(\cdot,\cdot)$ is a non-degenerate quadratic
form on the real vector space $\ch_2(\OO)$; its signature type has
one plus and the rest are minuses. Consequently, if $A>0$ then for
any $B\in \ch_2(\OO)$ one has
$$\det(A,B)^2\geq \det A\cdot \det B$$
and the equality holds if and only if $A$ is proportional to
$B$ with a real coefficient.

We do not present a detailed proof of this result since we are not
going to use it. Notice only that this result can be deduced
formally from the corresponding result for quaternionic matrices
proved in \cite{alesker-bsm-03}. Indeed all the entries of $A$ and
$B$ together  contain at most two non-real octonions, hence the
field generated by them is associative by Lemma \ref{L:oalg1}(iv).
\end{remark}

\item \begin{proposition}
Let $A,B\in \ch_2(\OO)$ be positive definite matrices. Then

(1) $A$ is invertible and $A^{-1}$ is positive definite;

(2) $tr(AB)>0$.
\end{proposition}

{\bf Proof.} (1) Entries of $A$ belong to a subfield of $\OO$ isomorphic to $\CC$. Hence by the case of complex matrices, the matrix $A$ is invertible, its entries belong to the same subfield, and $A^{-1}>0$.

(2) All the entries of $A$ and
$B$ together  contain at most two non-real octonions, hence the
field generated by them is associative by Lemma \ref{L:oalg1}(iv). Thus we may assume that $A$ and $B$ have quaternionic entries. Any quaternionic Hermitian matrix is diagonalizable (see e.g. Claim 1.1.7 in \cite{alesker-bsm-03}), namely
$$A=C^*DC,$$
where $D$ is real diagonal and $C^*C=I_2$. Then we have
\begin{eqnarray*}
tr(AB)=tr(C^*DCB)=tr(D(CBC^*)).
\end{eqnarray*}
The matrix $CBC^*\in \ch_2(\HH)$ is also positive definite as follows from the definition (over $\HH)$. Replacing $B$ with the latter matrix and $A$ with $D$ we may assume that $A$ is diagonal.
In that case
\begin{eqnarray*}
tr(AB)=a_{11}b_{11}+a_{22}b_{22}.
\end{eqnarray*}
The last expression is positive since all the terms are positive. \qed

\end{paragraphlist}

\def\hr{\ch_{16}(\RR)}
\def\ho{\ch_2(\OO)}

\section{The groups $SL_2(\OO)$ and $GL_2(\OO)$.}\label{Ss:sl2}
\begin{paragraphlist}
\item We recall in this section the definitions and basic properties of
the groups $SL_2(\OO)$ and $GL_2(\OO)$. We refer to \cite{sudbery} for the proofs and
further details.

\item An octonionic $2\times 2$-matrix is called {\itshape traceless} if
the sum of its diagonal elements is equal to zero. Every $2\times
2$-matrix $A$ with octonionic entries defines an $\RR$-linear
operator on $\OO^2$ by $\xi\mapsto A\cdot\xi$. However, the space of
such operators is not closed under the commutator due to lack of
associativity. One denotes by $\mathfrak{sl}_2(\OO)$ the Lie subalgebra of
$\mathfrak{gl}_{16}(\RR)$ generated by $\RR$-linear operators on
$\OO^2\simeq \RR^{16}$ determined by all traceless octonionic
matrices. This Lie algebra $\mathfrak{sl}_2(\OO)$ turns our to be semi-simple
\cite{sudbery} (see Theorem \ref{T:sl2-descr} below for details).
Note that any semi-simple Lie subalgebra of an algebraic group is a Lie
algebra of a closed algebraic subgroup (see e.g.
\cite{onishik-vinberg}, Ch. 3, \S 3.3). In our case this subgroup is
denoted by $SL_2(\OO)\subset GL(16,\RR)$.

\begin{theorem}[\cite{sudbery}]\label{T:sl2-descr}
(i) The Lie algebra $\mathfrak{sl}_2(\OO)$ is isomorphic to the Lie algebra
$\mathfrak{so}(9,1)$.

(ii) The Lie group $SL_2(\OO)$ is isomorphic to the group
$Spin(9,1)$ (which is the universal covering of the identity
component of the pseudo-orthogonal group $O(9,1)$).
\end{theorem}
\begin{remark}\label{R:spin9}
A maximal compact subgroup of $SL_2(\OO)\simeq Spin(9,1)$ is
isomorphic to the group $Spin(9)$ which is the universal covering of
the special orthogonal group $SO(9)$.
\end{remark}
\begin{remark}
In fact $SL_2(\OO)\subset SL_{16}(\RR)$. Indeed since the Lie algebra $\mathfrak{sl}_2(\OO)$ is semi-simple, it is equal to its own commutant.
Hence $\mathfrak{sl}_2(\OO)\subset \mathfrak{sl}_{16}(\RR)$. Since $SL_2(\OO)$
is connected it follows that $SL_2(\OO)\subset SL(16,\RR)$.
\end{remark}

\item Let us denote by $\mathfrak{gl}_2(\OO):=\RR\oplus \mathfrak{sl}_2(\OO)\subset \mathfrak{gl}_{16}(\RR)$. It is the Lie algebra of the Lie group
$GL_2(\OO):=\RR^{\ast}\cdot SL_2(\OO)$ which is the group generated by the group of homotheties $\RR^{\ast}$ and $SL_2(\OO)$.
$GL_2(\OO)$ is connected since $-I_2\in SL_2(\OO)$ (the latter follows from Lemma 14.77 in \cite{harvey}).

\item\label{item-lie-representation} $\mathfrak{sl}_2(\OO),\mathfrak{gl}_2(\OO) $ and $SL_2(\OO), GL_2(\OO)$ come with their {\itshape
fundamental representations} on $\OO^2=\RR^{16}$.

Furthermore the Lie algebra
$\mathfrak{sl}_2(\OO)$ acts on the space $\ch_2(\OO)$. This action is uniquely
characterized by the following property (see \cite{sudbery} for the
details): if $A$ is a traceless matrix it acts by $A\colon X\mapsto
-A^*X-XA$. Since the group $SL_2(\OO)$ is connected and simply
connected, this representation of $\mathfrak{sl}_2(\OO)$ integrates to a
representation of the group $SL_2(\OO)$ on $\ch_2(\OO)$.

This action naturally extends to $\mathfrak{gl}_2(\OO)$: for $\lam\in \RR\in \mathfrak{gl}_2(\OO)$ one has
$$\lam(X)=-2\lam X.$$
The latter action integrates to the representation of $GL_2(\OO)$ on $\ch_2(\OO)$ such that for any $\lam\in \RR^{\ast}$
$$\lam(X)=\lam^{-2} X.$$
This representation of $GL_2(\OO)$ in $\ch_2(\OO)$ will be denoted by $\rho$.

\begin{proposition}\label{P:det-rho-action}
(1) For any $g\in GL_2(\OO)\subset GL_{16}(\RR)$ one has $\det g>0$ where $\det$ is the usual determinant on $Mat_{16}(\RR)$.
 
(2) For the representation $\rho$ of the group $GL_2(\OO)\subset GL_{16}(\RR)$ in $\ch_2(\OO)$ and for any $g\in GL_2(\OO)$ and any $X\in \ch_2(\OO)$ one has
$$\det(g(X))=\det(g)^{-\frac{1}{4}}\det X,$$
where $\det$ in the left hand side in on $\ch_2(\OO)$, and $\det$ in the right hand side is on $Mat_{16}(\RR)$.
\end{proposition}
{\bf Proof.} (1) The group $GL_2(\OO)$ is connected. Hence $\det$ on it must be positive.

(2) Recall that $GL_2(\OO)=\RR^*\cdot SL_2(\OO)$. Assume first that $g=\lam\cdot I_{16}, \, \lam\in \RR^*$. Then 
$$\det(gX)=\det (\lam^{-2}X)=\lam^{-4}\det X.$$
On the other hand 
$$\det(g)^{-\frac{1}{4}}=\det(\lam I_{16})^{-\frac{1}{4}}=\lam^{-4}.$$

Let us assume now that $g\in SL_2(\OO)$. Then $\det(g(X))=\det X$ by \cite{baez}, Section 3.3.
Since $SL_2(\OO)\subset SL_{16}(\RR)$ then $\det(g)=1$ where $\det$ is on $16\times 16$ real matrix.
\qed

\item The following result is well known, see e.g., in \cite{alesker-octon}, Proposition 1.4.3.
\begin{proposition}\label{C:cone-preserved}
The group $GL_2(\OO)$ preserves the cone of positive definite
octonionic Hermitian matrices.
\end{proposition}
\def\clk{\overline{\ck}}

\item  Let us identify the (real) dual space to $\OO^2=\RR^{16}$ with itself using the standard inner product $\OO^2\times \OO^2\to \RR$ which can be described as follows
\begin{eqnarray}\label{E:euclid-structure}
\langle\xi,\eta\rangle= Re(\xi^*\cdot \eta).
\end{eqnarray}
Let $Mat_{2\times 2}(\OO)$ denote the space of square octonionic matrices of size 2.

\begin{lemma}\label{L:automorphism}
(i) The (Euclidean) dual operator to the operator $\xi\mapsto A\xi$, where $\xi\in\OO^2$ and $A\in Mat_{2\times 2}(\OO)$ is given by
$$\eta\mapsto A^*\eta.$$

(ii) The involutive automorphism of the group $GL_{16}(\RR)$ given by $X\mapsto (X^t)^{-1}$ preserves the subgroup $GL_2(\OO)$. Its restriction to this subgroup will be denoted by
$$\mathcal{S}\colon GL_2(\OO)\tilde\to GL_2(\OO).$$
\end{lemma}
{\bf Proof.} (i) We have
\begin{eqnarray*}
\langle \xi,A\eta\rangle:=Re(\xi^*A\eta)=Re((A^*\xi)^*\eta)=\langle A^*\xi,\eta\rangle.
\end{eqnarray*}

(ii) Recall that $GL_2(\OO)=\RR^*\cdot SL_2(\OO)$. Clearly for $\lam \in \RR^*$ one has $\cs(\lam)=\lam^{-1}\in GL_2(\OO)$. Now it suffices to show that $\cs$ preserves the Lie algebra $\mathfrak{sl}_2(\OO)$ of $SL_2(\OO)$.

The action of $\cs$ on the Lie algebra $\mathfrak{gl}_{16}(\RR)$ is given by $\cs(X)=-X^t$.
The Lie algebra $\mathfrak{sl}_2(\OO)$ is generated by transformations of the form $\hat A(\xi)=A\xi$, where $A\in Mat_{2\times 2}(\OO)$ is traceless, i.e. the sum of its diagonal entries vanishes.
Then, by part (i) one has
$$\cs(\hat A)(\xi)=-A^*\xi.$$
Thus $\cs(\hat A)\in \mathfrak{sl}_2(\OO)$, and the lemma follows. \qed

\item The following  lemma is essentially contained in
\cite{manogue-schray}, and we present a proof for convenience of the reader.
\begin{lemma}[\cite{manogue-schray}]\label{m-sch1}
The linear map $\ct\colon Sym_\RR^2(\OO^2)\to \ch_2(\OO)$ given by
$$\ct(\xi\otimes \eta)\mapsto \frac{1}{2}(\xi\cdot \eta^*+\eta\cdot\xi^*)$$
is $GL_2(\OO)$-equivariant, where the action of $GL_2(\OO)$ on the source space is the second symmetric power of the representation $\cs$,
and its action on the target is the representation $\rho$.\footnote{Here and throughout the paper we denote by $\otimes$ the product both in the symmetric algebra and in the tensor algebra.}
\end{lemma}

{\bf Proof.} The equivariance of the map under $\RR^*\subset GL_2(\OO)$ is evident. It suffices to prove the equivariance under the action of the Lie algebra $\mathfrak{sl}_2(\OO)$.
Moreover it suffices to check that for a set of generators of the Lie algebra $\mathfrak{sl}_2(\OO)$.
 For let $A\in Mat_{2\times 2}(\OO)$ be a traceless matrix.
It suffices to show that for any traceless $A\in Mat_{2\times 2}(\OO)$ and any 2-column $\xi\in \OO^2$ one has the equality of $2\times 2$ matrices
\begin{eqnarray}\label{E:equiv-need}
(A^*\xi)\xi^*+\xi(A^*\xi)^*=A^*(\xi\xi^*)+(\xi\xi^*)A.
\end{eqnarray}
Let us write explicitly
\begin{eqnarray*}
\xi=\left[\begin{array}{c}
          \xi_1\\
          \xi_2
          \end{array}\right],\,\,\, A=\left[\begin{array}{cc}
                                      a&b\\
                                      c&-a
                                       \end{array}\right].
\end{eqnarray*}

The left hand side  of (\ref{E:equiv-need}) is then equal to
\begin{eqnarray*}
lhs=\left(\left[\begin{array}{cc}
                                      \bar a&\bar c\\
                                      \bar b&-\bar a
                                       \end{array}\right]\left[\begin{array}{c}
          \xi_1\\
          \xi_2
          \end{array}\right]\right)[\bar \xi_1,\bar \xi_2]+\left[\begin{array}{c}
          \xi_1\\
          \xi_2
          \end{array}\right]\left([\bar\xi_1,\bar\xi_2]\left[\begin{array}{cc}
                                      a&b\\
                                      c&-a
                                       \end{array}\right]\right)=\\
                                       \left[\begin{array}{cc}
                                       (\bar a\xi_1+ \bar c\xi_2)\bar \xi_1& (\bar a\xi_1 +\bar c\xi_2)\bar \xi_2\\
                                       (\bar b\xi_1-\bar a\xi_2)\bar \xi_1&(\bar b\xi_1-\bar a\xi_2)\bar \xi_2
                                       \end{array}\right]+\left[\begin{array}{cc}
                                                               \xi_1(\bar\xi_1 a+\bar\xi_2 c)&\xi_1(\bar\xi_1b-\bar \xi_2 a)\\
                                                               \xi_2(\bar\xi_1 a+\bar\xi_2 c)&\xi_2(\bar\xi_1b-\bar \xi_2 a)
                                                                \end{array}\right]=\\
                                       \left[\begin{array}{cc}
                                       2|\xi_1|^2Re(a)+2Re(\xi_1\bar\xi_2c)&\bar c|\xi_2|^2+b|\xi_1|^2+(\bar a\xi_1)\bar \xi_2-\xi_1(\bar\xi_2a)\\
                                        c|\xi_2|^2+\bar b|\xi_1|^2+\xi_2(\bar\xi_1 a)-(\bar a\xi_2)\bar \xi_1& -2|\xi_2|^2 Re(a)+2Re(\xi_2\bar\xi_1 b)
                                       \end{array}\right].
\end{eqnarray*}
Let us compute the right hand side of (\ref{E:equiv-need}):
\begin{eqnarray*}
rhs=\left[\begin{array}{cc}
          \bar a&\bar c\\
          \bar b&-\bar a
           \end{array}\right]\left[\begin{array}{cc}
                                   |\xi_1|^2&\xi_1\bar\xi_2\\
                                   \xi_2\bar\xi_1&|\xi_2|^2
                                    \end{array}\right]+\left[\begin{array}{cc}
                                   |\xi_1|^2&\xi_1\bar\xi_2\\
                                   \xi_2\bar\xi_1&|\xi_2|^2
                                    \end{array}\right]\left[\begin{array}{cc}
                                       a&b\\
                                      c&- a
                                       \end{array}\right]=\\
                                       \left[\begin{array}{cc}
                                       \bar a |\xi_1|^2+\bar c(\xi_2\bar\xi_1)&\bar a(\xi_1\bar\xi_2)+\bar c|\xi_2|^2\\
                                       \bar b|\xi_1|^2-\bar a(\xi_2\bar\xi_1)&\bar b(\xi_1\bar\xi_2)-\bar a|\xi_2|^2
                                       \end{array}\right]+
                                       \left[\begin{array}{cc}
                                                           a|\xi_1|^2+(\xi_1\bar \xi_2)c&b|\xi_1|^2-(\xi_1\bar \xi_2)a\\
                                                           (\xi_2\bar\xi_1)a+c|\xi_2|^2&( \xi_2\bar \xi_1)b-a|\xi_2|^2
                                                           \end{array}\right]=\\
                                                           \left[\begin{array}{cc}
                                                           2|\xi_1|^2Re (a)+2Re(\xi_1\bar\xi_2c)&\bar c|\xi_2|^2+b|\xi_1|^2+\bar a(\xi_1\bar\xi_2)-(\xi_1\bar\xi_2)a\\
                                                           c|\xi_2|^2+\bar b|\xi_1|^2+(\xi_2\bar \xi_1)a-\bar a(\xi_2\bar \xi_1)&-2|\xi_2|^2Re(a)+2 Re(\xi_2\bar\xi_1 b)
                                                           \end{array}\right].
\end{eqnarray*}
We have to show that $lhs=rhs$. For the diagonal terms of these matrices the identity is obvious. Since both matrices are Hermitian it suffices to show the equality of the terms with indices $(1,2)$,
 that is
$$(\bar a\xi_1)\bar \xi_2-\xi_1(\bar\xi_2a)=\bar a(\xi_1\bar\xi_2)-(\xi_1\bar\xi_2)a.$$
To simplify the notation we have to show that
\begin{eqnarray}\label{E:1equality}
(\bar a b)c-b(ca)=\bar a(bc)-(bc)a \mbox{ for all } a,b,c\in \OO.
\end{eqnarray}

For $a\in \RR$ both sides of (\ref{E:1equality}) vanish. Thus by linearity we may and will assume that $a$ is purely imaginary, i.e. $\bar a=-a$. We have to prove than that
$$(ab)c+b(ca)=a(bc)+(bc)a.$$
Denote $[a,b,c]:=(ab)c-a(bc)$. The the last equality is equivalent to
\begin{eqnarray}\label{E:associator}
[a,b,c]=[b,c,a].
\end{eqnarray}
By Theorem 15.11(iii) of \cite{adams-lectures-exceptional-lie} it is alternating function of 3 variables. Hence (\ref{E:associator}) holds, and the lemma follows. \qed


\item  The following  lemma is essentially contained in
\cite{manogue-schray}, and we present a proof for convenience of the reader as in \cite{alesker-octon}.
\begin{lemma}[\cite{manogue-schray}]\label{m-sch2}
(i) The group $SL_2(\OO)$ acts naturally on $\OO\PP^1$, namely for
any $\phi\in SL_2(\OO)$ and any $L\in \OO\PP^1$ the subspace
$\phi(L)$ is an octonionic projective line. This action is transitive.

(ii) For any $L\in \OO\PP^1$ and any $\phi\in SL_2(\OO)$ the
restriction
$$\phi|_L\colon L\to \phi(L)$$
is a conformal linear map.
\end{lemma}
{\bf Proof.} Let $\phi\in SL_2(\OO)$. We first show that if for two unit vectors $\xi,\eta\in \OO^2$ one has $\xi\sim\eta$ (as defined in Lemma \ref{L:properties-op1}(iv))
then $|\phi(\xi)|=|\phi(\eta)|$.
To do so observe that for any $v\in \OO^2$ one has
$$|v|^2=v^*v=Tr(vv^*),$$
where $Tr$ stands for the trace and denotes sum of the diagonal elements. Therefore
\begin{eqnarray*}
|\phi(\xi)|^2=Tr(\phi(\xi)\cdot\phi(\xi)^*)\overset{\mbox{Lemma } \ref{m-sch1}}{=}Tr(\phi(\xi\cdot\xi^*))=\\
Tr(\phi(\eta\cdot\eta^*))\overset{\mbox{Lemma } \ref{m-sch1}}{=} Tr(\phi(\eta)\cdot\phi(\eta)^*)=|\phi(\eta)|^2.
\end{eqnarray*}
Thus $|\phi(\xi)|=|\phi(\eta)|$. To prove both parts of the lemma (except transitivity) it remains to show that for $\xi,\eta$ as above one has $\phi(\xi)\sim \phi(\eta)$, namely
$\phi(\xi)\cdot\phi(\xi)^*=\phi(\eta)\cdot\phi(\eta)^*$. Applying Lemma \ref{m-sch1} twice we get
$$\phi(\xi)\cdot\phi(\xi)^*=\phi(\xi\xi^*)=\phi(\eta\eta^*)=\phi(\eta)\cdot \phi(\eta)^*.$$

Finally let us prove the transitivity of the action of $SL_2(\OO)$ on $\OO\PP^1$. It suffices to show that for any non-zero vector $(\alp,b)\in\OO^2$ with $\alp\in \RR$ there exists $\phi\in SL_2(\OO)$ such that
\begin{eqnarray}\label{E:phi-transit}
\phi\left(\left[\begin{array}{c}
                   1\\
                   0
                   \end{array}\right]\right)=\left[\begin{array}{c}
                   \alp\\
                   b
                   \end{array}\right].
                   \end{eqnarray}
The octonion $b$ can be written $b=\beta+s \gamma$, where $\beta,\gamma\in \RR$, and $s$ is a purely imaginary octonion of norm 1. Thus $s^2=-1$. The subalgebra generated by $s$ and $1$ is isomorphic to $\CC$. It suffices to construct a matrix $\phi$
with entries in this subalgebra. Thus for simplicity of notation we may assume that $b\in \CC$. It is well known that $SL_2(\CC)$ acts transitively on $\CC^2\backslash\{0\}$. Hence there is $\phi\in SL_2(\CC)$ satisfying (\ref{E:phi-transit}). Furthermore it is
well known that $\exp\colon \frak{sl}_2(\CC)\to SL_2(\CC)$ is onto (this can be seen using the Jordan form of a complex matrix). Thus $\phi=\exp A$ where $A\in\frak{sl}_2(\CC)\subset \frak{sl}_2(\OO)$.
\qed

\item Recall that the maps $j$ and  $\theta$ were defined in (\ref{D:j}) and (\ref{Def:theta}) respectively. Let us denote by $\ch_{16}^0(\RR)$ the subspace of $\ch_{16}(\RR)$ consisting of real quadratic forms
whose restriction to any projective line is proportional to the square of the standard Euclidean norm on $\OO^2=\RR^{16}$. Lemma \ref{m-sch2} easily implies:
\begin{lemma}\label{L:ch0}
 $\ch^0_{16}(\RR)$ is a $GL_2(\OO)$-invariant subspace of $\ch_{16}(\RR)$.
 \end{lemma}

\item
\begin{lemma}\label{L:image-j}
The following inclusion holds
$$j(\ch_2(\OO))\subset \ch_{16}^0(\RR).$$
\end{lemma}
{\bf Proof.} Observe that any octonionic line (except the one at infinity) is spanned by a vector $\xi=\left[\begin{array}{cc}
                                                            1\\
                                                            a
                                                            \end{array}\right]$.  For any such $\xi$ and any $A=\left[\begin{array}{cc}
                                     \alp&x\\
                                     \bar x &\beta
                                   \end{array}\right]\in \ch_2(\OO)$ we have
\begin{eqnarray*}
j(A)=\left(\left[\begin{array}{cc}
                                                            q\\
                                                            aq
                                                            \end{array}\right]\right)=Re\left([\bar q,\bar q\bar a]\left[\begin{array}{cc}
                                     \alp&x\\
                                     \bar x &\beta
                                   \end{array}\right]\left[\begin{array}{cc}
                                                            q\\
                                                            aq
                                                            \end{array}\right]\right)=\\
Re(\alp|q|^2+\bar q x(aq)+(\bar q\bar a)\bar xq+\beta|aq|^2)=\\
|q|^2(\alp+\beta|a|^2) +2Re(\bar qx (aq))\overset{\mbox{Lemma }\ref{L:oalg1}(v)}{=}\\
|q|^2(\alp+\beta|a|^2+2Re(ax)),
\end{eqnarray*}
where the last equality follows from the fact that any subalgebra generated by any two octonions is associative. The lemma follows. \qed

\item \begin{lemma}\label{L-average}
Let $B\in \ch_{16}(\RR)$. Then, for any $\xi\in\OO^2,\, \xi\ne 0$ one has
\begin{eqnarray}\label{E:lll}
((j\circ\theta)(B))(\xi)=Re(\xi^*\theta(B)\xi)=\left(\int_{\SS_\xi}b(\eta)d\eta\right)\cdot |\xi|^2,
\end{eqnarray}
where $\SS_\xi$ is the unit 7-sphere (with respect to the standard Euclidean metric on $\OO^2=\RR^{16}$) in the octonionic line spanned by $\xi$, and $d\eta$ is the probability Lebesgue measure on this sphere.
\end{lemma}
{\bf Proof.} The first equality follows from the definition, see Section \ref{S:linear-algebra}, paragraph \ref{item:lin-alg2}. In \cite{alesker-octon}, Lemma 1.2.1, it was shown that if $\xi=\left[\begin{array}{c}
                                                                                      1\\
                                                                                      a
                                                                                    \end{array}\right]$ one has \begin{eqnarray*}
Re(\xi^*\theta(B)\xi)=\int_{x\in \SS^7}b(\xi\cdot x)dx,
\end{eqnarray*}
where the integration is with respect to the rotation invariant probability measure on $\SS^7\subset\OO$. The last integral is clearly equal to $\left(\int_{\SS_\xi}b(\eta)d\eta\right)\cdot |\xi|^2$.
Thus, the lemma holds for vectors $\xi$ with the first coordinate equal to 1.

By Lemma \ref{L:image-j} one has     $j(\ch_2(\OO))\subset \ch_{16}^0(\RR)$. The integral in (\ref{E:lll}) defines a function on $\OO^2$
which is proportional to the square of the standard Euclidean norm on each octonionic line which coincides with $(j\circ\theta)(B)$ on each vector $\xi=\left[\begin{array}{c}
                                                                                      1\\
                                                                                      a
                                                                                    \end{array}\right]$. Since any octonionic line contains a vector of such a form,
these two functions (the integral in (\ref{E:lll}) and $(j\circ\theta)(B)$) coincide on $\OO^2$. \qed

\begin{corollary}\label{Cor:isomorphisms}
The linear maps
\begin{eqnarray*}
j\colon \ch_2(\OO)\to \ch_{16}^0(\RR),\\
\theta\big|_{\ch_{16}^0(\RR)}\colon \ch_{16}^0(\RR)\to \ch_2(\OO)
\end{eqnarray*}
are isomorphisms and inverse to each other.
\end{corollary}
{\bf Proof.} First $\theta\circ j=id$ by Lemma \ref{L:la2}. Hence, it remains to show that
$$j\circ \theta\big|_{\ch_{16}^0(\RR)}\colon \ch_{16}^0(\RR)\to \ch_{16}^0(\RR),$$
is the identity.

For a real quadratic form $b$ on $\RR^{16}$ we denote by $B\in \ch_{16}(\RR)$ its matrix in the standard basis.
 By Lemma \ref{L-average} for any $\xi\in \OO^2\backslash \{0\}$ one has
$$((j\circ\theta)(B))(\xi)=\left(\int_{\SS_\xi}b(\eta)d\eta\right)\cdot |\xi|^2.$$
Since $B\in \ch_{16}^0(\RR)$, the function under the integral is constant on the unit sphere $\SS_\xi$. Hence the right hand side is equal to $b(\xi)$. The result follows. \qed

\item
\begin{lemma}\label{L:j-theta-equiv}
The map $$j\circ \theta\colon \ch_{16}(\RR)\to \ch_{16}^0(\RR)$$
is $GL_2(\OO)$-equivariant.
\end{lemma}
{\bf Proof.}  For any $\xi\in\OO^2\backslash\{0\}$ and any
$g\in GL_2(\OO)$ one has
\begin{eqnarray*}
g((j\circ\theta)(B))(\xi)\overset{\mbox{Lemma } \ref{L-average}}{=}\left(\int_{\SS_{g^{-1}\xi}}b(\eta)d\eta\right)\cdot |g^{-1}\xi|^2 =\\
\left(\int_{\SS_\xi}b\left(\frac{g^{-1}\zeta}{|g^{-1}\zeta|}\right)d\zeta\right)|g^{-1}\xi|^2=
\left(\int_{\SS_\xi}b\left(\frac{g^{-1}\zeta}{g^{-1}(\xi/|\xi)|}\right)d\zeta\right)|g^{-1}\xi|^2=\\
\left(\int_{\SS_\xi}b(g^{-1}\zeta)d\zeta\right)|\xi|^2=\left(\int_{\SS_\xi}(gb)(\zeta)d\zeta\right)|\xi|^2\overset{\mbox{Lemma }\ref{L-average}}{=}\\
((j\circ\theta)(gB))(\xi).
\end{eqnarray*}
 \qed

\item
\begin{proposition}\label{L:theta-equivar}
The maps $\theta\colon \hr\to \ho$ and $j\colon \ho\to \hr$ are $GL_2(\OO)$-equivariant with the action of $GL_2(\OO)$ on $\ch_2(\OO)$ as described in paragraph \ref{item-lie-representation} of this subsection,
and its action on $\ch_{16}(\RR)$ is the standard action on quadratic forms: $(gb)(\xi)=b(g^{-1}\xi)$.
\end{proposition}
{\bf Proof.}
The equivariance of $j$ was announced without proof in \cite{alesker-octon}, Remark 1.5.3.
We present a proof below. Recall that $GL_2(\OO)=\RR^*\cdot SL_2(\OO)$. Clearly $j$ is $\RR^*$-equivariant, hence it suffices to prove $SL_2(\OO)$-equivariance.
Namely, we have to show that for any $g\in SL_2(\OO)$ and any $X\in \ch_2(\OO)$ one has
$$j(gX)=g(jX).$$
This is equivalent that for any $\xi\in \RR^{16}=\OO^2$ the following identity holds
$$Re(\xi^*(gX)\xi)=Re\left((g^{-1}\xi)^*X(g^{-1}\xi)\right).$$
It suffices to prove the infinitesimal version of this equality for generators of the Lie algebra $\frak{sl}_2(\OO)$. Namely, one has to show that for any $A\in \frak{sl}_2(\OO)$
$$Re(\xi^*(A^*X+XA)\xi)=Re((A\xi)^*X\xi+\xi^*X(A\xi)).$$
Since in each side the second summand is the Hermitian conjugate of the first one this is equivalent to
\begin{eqnarray}\label{E:jjj}
Re(\xi^*(XA)\xi)=Re(\xi^*X(A\xi)).
\end{eqnarray}

 Let us compute the right hand side of (\ref{E:jjj}) for $A$ being a traceless matrix, i.e. $A=\left[\begin{array}{cc}
                                                                                                     z&x\\
                                                                                                     y&-z
                                                                                                     \end{array}\right]$, and let $X=\left[\begin{array}{cc}
                                                                                                                                   \alp&q\\
                                                                                                                                   \bar q&\beta
                                                                                                                                           \end{array}\right]\in\ch_2(\OO)$:
\begin{eqnarray*}\label{E:mou-2}
\mbox{rhs of } (\ref{E:jjj})=Re\left\{[\bar\xi_1,\bar\xi_2]\left[\begin{array}{cc}
                                                                  \alp&q\\
                                                                  \bar q&\beta
                                                                  \end{array}\right]\left(\left[\begin{array}{cc}
                                                                                                     z&x\\
                                                                                                     y&-z
                                                                                                     \end{array}\right]\left[\begin{array}{c}
                                                                                                                              \xi_1\\
                                                                                                                              \xi_2
                                                                                                                              \end{array}\right]\right)\right\}=\\
Re(\bar\xi_1(\alp z)\xi_1+\bar \xi_1 q(y\xi_1)-\bar\xi_2(\beta z)\xi_2+ \bar \xi_2\bar q (x\xi_2)+\\
\bar\xi_1\alp x\xi_2+\bar \xi_2\beta(y\xi_1)-\bar\xi_1 q(z\xi_2)+\bar \xi_2\bar q(z\xi_1))=\\
Re((\alp z)\xi_1\bar\xi_1+ q(y\xi_1)\bar \xi_1-(\beta z)\xi_2\bar\xi_2 + \bar q (x\xi_2)\bar \xi_2)+\\
Re[\bar\xi_1\alp x\xi_2+(\bar \xi_1\bar y)\beta\xi_2-(\bar\xi_1 q)(z\xi_2)+(\bar \xi_1\bar z)(q\xi_2)]=\\
Re(\alp z+qy)|\xi_1|^2+Re(\bar q x-\beta z)|\xi_2|^2+Re[\bar\xi_1(\alp x+\beta\bar y)\xi_2-(\bar\xi_1 q)(z\xi_2)+(\bar \xi_1\bar z)(q\xi_2)].
\end{eqnarray*}
For the left hand side of (\ref{E:jjj}) we have
\begin{eqnarray*}
\mbox{lhs of } (\ref{E:jjj})= Re\left([\bar\xi_1,\bar\xi_2]\left(\left[\begin{array}{cc}
                                                              \alp&q\\
                                                              \bar q&\beta
                                                              \end{array}\right]\left[\begin{array}{cc}
                                                                                      z&x\\
                                                                                      y&-z
                                                                                      \end{array}\right]\right)\left[\begin{array}{c}
                                                                                                                              \xi_1\\
                                                                                                                              \xi_2
                                                                                                                              \end{array}\right]\right)=\\
Re\left([\bar\xi_1,\bar\xi_2]\left[\begin{array}{cc}
                                                              \alp z+qy&\alp x-qz\\
                                                              \bar q z+\beta y&\bar q x-\beta z
                                                              \end{array}\right]\left[\begin{array}{c}
                                                                                                                              \xi_1\\
                                                                                                                              \xi_2
                                                                                                                              \end{array}\right]\right)=\\
Re\left(\bar\xi_1(\alp z+qy)\xi_1+\bar \xi_2(\bar q x-\beta z)\xi_2 +\bar\xi_1(\alp x-qz)\xi_2+\bar \xi_2(\bar qz+\beta y)\xi_1\right)=\\
Re(\alp z+qy)|\xi_1|^2 +Re(\bar q x-\beta z)|\xi_2|^2+ Re(\bar\xi_1(\alp x+\beta \bar y-qz+\bar z q)\xi_2).
\end{eqnarray*}
Comparing this with the right hand side of (\ref{E:jjj}) it remains to show that
\begin{eqnarray}\label{E:jjj2}
Re[-(\bar\xi_1 q)(z\xi_2)+(\bar \xi_1\bar z)(q\xi_2)]=Re[\bar\xi_1(-qz+\bar z q)\xi_2].
\end{eqnarray}
Note that if $z\in \RR$ then obviously both sides vanish. Thus let us assume that $z$ is purely imaginary, i.e. $\bar z=-z$.
In this case (\ref{E:jjj2}) reads
$$Re[(\bar\xi_1 q)(z\xi_2)+(\bar \xi_1 z)(q\xi_2)]=Re[\bar\xi_1(qz+ z q)\xi_2].$$
This identity holds by Corollary \ref{Cor-moufang}.


The equivariance of  $\theta$  was proven in \cite{alesker-octon}, Lemma 1.5.2. However the proof there contains a gap, so we present a complete proof below.

By Lemma \ref{L:j-theta-equiv} the map $j\circ \theta\colon \ch_{16}(\RR)\to \ch_{16}^0(\RR)$ is $GL_2(\OO)$-equivariant, as well as $j$ as was shown above.
Hence, $\theta=j^{-1}\circ (j\circ \theta)$ is $GL_2(\OO)$-equivariant. \qed

\end{paragraphlist}

\section{Diagonalization.}\label{S:diagonalization}
\begin{paragraphlist}
\item The following proposition seems to be folklore.
\begin{proposition}\label{P:diagonalization}
(1) For any $A\in \ch_2(\OO)$ there exists a transformation from $Spin(9)$ which maps $A$ to a real diagonal matrix.
\newline
(2)Let $A,B\in\ch_2(\OO)$, $A>0$. Then there exists a transformation from $Spin(1,9)$ which maps $A$ to a matrix proportional to $I_2$, and $B$ to a diagonal matrix.
\end{proposition}
{\bf Proof.} (1) The entries of $A$ are contained in a subfield isomorphic to $\CC$. Let us assume, without loss of generality and for simplicity of notation, that the latter subfield is $\CC$.

Recall that any $L\in \frak{sl}_2(\OO)$ acts on $A$ as $A\mapsto -L^*A-AL=:\hat L(A)$. If $L\in \frak{sl}_2(\CC)\subset \frak{sl}_2(\OO)$ it follows that
$$e^{\hat L}(A)=(e^{-L})^*Ae^{-L}.$$

The 1-parametric subgroup $e^{t L}$ is contained in $ SL_2(\OO)$. It acts on $A$ as
$$e^{t\hat L}(A)=\sum_{k=0}^\infty \frac{t^k}{k!} \hat L^k(A)=\sum_{k=0}^\infty \frac{t^k}{k!} \underset{k\mbox{ times}}{\underbrace{\hat L(\dots \hat L(A)\dots)}}.$$
Clearly $\frak{sl}_2(\CC)\subset\frak{sl}_2(\OO)$. If $L\in \frak{sl}_2(\CC)$  then $L$ and $A$ generated an associative subalgebra. But
$$\frac{d}{dt}\big|_0e^{t\hat L}(A)=-L^*A-AL.$$

Hence the subgroup $SL_2(\CC)\subset SL_2(\OO)$ acts on $\ch_2(\OO)$ as
$$X(A)=(X^{-1})^*AX^{-1}.$$
Note that any complex Hermitian matrix can be diagonalized using such transformation with above $X\in SU(2)$ being a complex matrix. The latter subgroup is contained in $Spin(9)$.

\hfill

(2) By part (1) we may apply a transformation from $Spin(9)$ which diagonalizes $A$, i.e. makes it $A'=\left[\begin{array}{cc}
                                                                                                              \lam_1&0\\
                                                                                                              0&\lam_2
                                                                                                              \end{array}\right]$ where $\lam_{1,2}>0$.
Then replace $A'$ with $  \left[\begin{array}{cc}
          r&0\\
          0&r^{-1}
          \end{array}\right]A'\left[\begin{array}{cc}
          r&0\\
          0&r^{-1}
          \end{array}\right]$ with $r\in \RR\backslash\{0\}$.        For an appropriate $r$ we can get $cI_2$.
Since $Spin(9)$ preserves the identity matrix $I_2$, we may apply further an appropriate transformation from $Spin(9)$ and get the first matrix equal to $cI_2$ and the other one diagonal. \qed

\end{paragraphlist}

\section{Further properties of the octonionic
Hessian.}\label{Ss:octHessian} \setcounter{theorem}{0}
\begin{proposition}[\cite{alesker-octon}, Prop. 1.5.1]\label{P:ohs1}
Let $f\colon \OO^2\to \RR$ be a $C^2$-smooth function. Let $A\in
SL_2(\OO)$. Then,
$$\left(\frac{\pt^2}{\pt\bar q_i \pt
q_j}(f(A^{-1}q))\right)=A\left(\dfq (A^{-1}q)\right),$$ where $A$ in the right
hand side denotes the induced action of $A$ on $\ho$.
\end{proposition}
{\bf Proof.} By translation it is enough to check the above
equality at $q=0$. Moreover we may and will assume that $f$ is a
quadratic form. By the definition of $\theta$, the proposition reduces to the $GL_2(\OO)$-equivariance of $\theta$; that holds thanks to Proposition \ref{L:theta-equivar}.
\qed

\section{$G$-affine manifolds.}\label{S:G-mflds}
\begin{paragraphlist}

\item The results of this section are probably folklore. Let $G\subset GL(n,\RR)$ be a subgroup. We introduce a rather special class of manifolds which we call $G$-affine.
\begin{definition}\label{D:G-affine}
(1) We say that an atlas of charts $\{(U_\alp,\phi_\alp)\}$ on a smooth manifold $M^n$, where $\phi_\alp\colon U_\alp \to \phi_\alp(U_\alp)\subset \RR^n$
is a diffeomorphism onto an open subset of $\RR^n$, defines a $G$-affine structure on $M$ if for any $\alp,\beta$ the diffeomorphism
$$\phi_\beta\circ \phi_\alp^{-1}\colon \phi_\alp(U_\alp\cap U_\beta)\to \phi_\beta(U_\alp\cap U_\beta)$$
is a composition of a transformation  from $G$ and a translation.

(2) Two atlases of charts on $M$ are called equivalent if their union defines a $G$-affine structure (in other words, the transition maps between charts
from different atlases are compositions of maps from $G$ and translations.

(3) A $G$-affine structure on $M$ is an equivalence class of atlases of charts.
\end{definition}

The following example is relevant for the current study.
\begin{example}
Let $G=\{Id\}$ be the trivial group. Then a torus $\RR^n/\Lambda$ is a $G$-affine manifold, where $\Lambda\subset \RR^n$ is a lattice.
\end{example}

\begin{remark}
In this work the groups of interests are $G=GL_2(\OO),SL_2(\OO)\simeq  Spin(1,9), Spin(9)\subset GL_{16}(\RR)$.
\end{remark}

\item The following proposition is obvious.
\begin{proposition}\label{P:connection}
On a $G$-affine manifold there exists a unique torsion free flat connection on the tangent bundle such that on each chart of the given atlas $\{(U_\alp,\phi_\alp)\}$ this connection is the trivial one.
The holonomy of this connection is contained in $G$.
\end{proposition}

\item On a $G$-affine manifold $M$ the real Hessian of a $C^2$-smooth function is well defined and takes values in the vector bundle $Sym^2(T^*M)$ whose fiber over
$x\in M$ can be identified with the space of quadratic forms on
the tangent space $T_xM$. This Hessian coincides with the second covariant derivative with respect to the induced flat
connection.

\item Recall that densities are sections of the tensor product of the lines bundles of top degree differential forms and of the orientation bundle. Compactly supported densities can be integrated over the manifold.
The following result is well known, also it can be easily proved by the reader.
 \begin{proposition}\label{P:volume-pr}
Assume that the group $G$ is contained in the group of (linear) volume preserving transformations.
Then on any $G$-affine manifold there exists unique nowhere vanishing density
which is parallel with respect to the connection from Proposition \ref{P:connection}.
\end{proposition}

\item Let us remind the notion of $\lam$-density on a smooth manifold for a real number $\lam$. Let $X^n$ be a smooth manifold (not necessarily $G$-affine). Let $(V_\alp,\psi_\alp)$
be an atlas of charts. We have diffeomorphisms
\begin{eqnarray}\label{E:transition-map}
\psi_\beta\circ \psi_\alp^{-1}\colon \psi_\alp(V_\alp\cap V_\beta)\to \psi_\beta(V_\alp\cap V_\beta).
\end{eqnarray}
Let $j_{\alp\beta}$ denote the absolute value of the determinant of the inverse of the differential of this map. Clearly
\begin{eqnarray}\label{E:cocycle}
j_{\alp\beta}\cdot j_{\beta\gamma}=j_{\alp\gamma}.
\end{eqnarray}
Let us define the real line bundle of $\lam$-densities over $X$ as follows. Over each chart $V_\alp$ we choose the trivial line bundle, and the transition map from $V_\alp$ to $V_\beta$
is equal to $(j_{\alp\beta})^\lam$. Due to (\ref{E:cocycle}) this defines the line bundle. Its sections are called $\lam$-densities. The $1$-densities coincide with densities from the previous paragraph.

The line bundle of $\lam$-densities is canonically oriented: in any chart its orientation coincides with the natural one of $\RR$. A section of this line bundle is called positive (resp. non-negative)
if its value at any point is positive (resp. non-negative) with respect to the given orientation.

Any nowhere vanishing $\lam$-density can be raised to a real power $\mu$, and the result is a $\lam\mu$-density.

\item Let $M^{16}$ be an $GL_2(\OO)$-affine manifold in the sense of Definition \ref{D:G-affine}. The {\itshape tangent bundle} $TM$ can be described as follows.

Let us fix the trivial bundle with fiber $\OO^2$ over each map $V_\alp$.
Let us glue them over each pairwise intersection $V_\alp\cap V_\beta$
according to the transition map (\ref{E:transition-map}) as follows. The differential $d(\psi_\beta\circ \psi_\alp^{-1})$ is a locally constant matrix with values in $GL_2(\OO)$.
 We identify the restriction to $V_\alp\cap V_\beta$ of the trivial bundle $V_\alp\times \OO^2$ with the restriction to the same intersection of the trivial bundle $V_\beta\times \OO^2$ via the map
$$(x,X)\mapsto \left(x,d_x(\psi_\beta\circ \psi_\alp^{-1})\right)X.$$

\item The cotangent bundle $T^*M$ over the $GL_2(\OO)$-affine manifold $M$ is described similarly identifying
the restriction to $V_\alp\cap V_\beta$ of the trivial bundle $V_\alp\times \OO^2$ with the restriction to the same intersection of the trivial bundle $V_\beta\times \OO^2$ via the map
$$(x,X)\mapsto \left(x,\cs(d_x(\psi_\beta\circ \psi_\alp^{-1}))X\right),$$
where $\cs(g)=(g^t)^{-1}$ is the map from Lemma \ref{L:automorphism}(ii).

\item For a $GL_2(\OO)$-affine manifold $M$  let us define a vector bundle $\uch$ over $M$ with the fiber
$\ch_2(\OO)$ which is the space of octonionic Hermitian $2\times 2$ matrices. We will see that $\uch$ can be considered as a sub-bundle of the
bundle $Sym^2_\RR(T^*M)$ of real quadratic forms on the tangent bundle $TM$ and the octonionic Hessian of a real valued function is a section of $\uch$.

Let us fix the trivial bundle with fiber $\ch_2(\OO)$ over each map $U_\alp$. Let us glue them over each pairwise intersection $U_\alp\cap U_\beta$
according to the transition map (\ref{E:transition-map}) as follows.
Let $\rho$ be the representation of $GL_2(\OO)$ in $\ch_2(\OO)$ introduced in Subsection \ref{Ss:sl2}, paragraph \ref{item-lie-representation}.
Let us identify the restriction to $U_\alp\cap U_\beta$ of the trivial bundle $U_\alp\times \ch_2(\OO)$ with the restriction to the same intersection of the trivial bundle $U_\beta\times \ch_2(\OO)$ via the map
\begin{eqnarray}\label{E:abvgd}
(x,X)\mapsto \left(x,\rho\left(d_x(\psi_\beta\circ \psi_\alp^{-1})\right)X\right).
\end{eqnarray}

These identifications define the vector bundle over $M$ which we denote by $\uch$.

\item \label{item1} The $GL_2(\OO)$-equivariant linear map $\ct\colon \OO^2\otimes_\RR \OO^2\to \ch_2(\OO)$ from Lemma \ref{m-sch1} induces the morphism of vector bundles
\begin{eqnarray}\label{T-bundles-morphism}
\ct\colon Sym_\RR^2(T^*M)\to \uch,
\end{eqnarray}
where $T^*M$ is the cotangent bundle of $M$.

\item By Proposition \ref{C:cone-preserved} the group $GL_2(\OO)$ preserves the cone of positive definite octonionic Hermitian matrices. Hence in each fiber of the bundle $\uch$
there is a convex cone which coincides with the latter cone in each coordinate chart. Hence one can talk of positive definite and non-negative definite (say, continuous) sections of $\uch$.

\item It is easy to see that the octonionic Hessian of a real smooth function $f$, which is locally written as $Hess_\OO f$, is a section of $\uch$.

Furthermore given a section $\xi$ of $\uch$, its pointwise determinant $\det \xi$ is a section of the line bundle of $1/4$-densities. 
Indeed, let $\xi$ be obtained using identifications
(\ref{E:abvgd}) where $X\in \ch_2(\OO)$. Then by Proposition \ref{P:det-rho-action} one has the identification
$$(x,\det X)\mapsto (x,\det(d_x(\psi_\beta\circ\psi_\alp^{-1}))^{-\frac{1}{4}}\det X),$$
where the first and the third determinants are defined on $\ch_2(\OO)$, and the second one on $Mat_{16}(\RR)$. 
This implies that $\det \xi$ is a $1/4$-density.

\item \label{item-Theta} If the manifold $M$ is $SL_2(\OO)$-affine manifold, rather than $GL_2(\OO)$, then there is a family of $\lam$-densities $\{\Theta^\lam\}_{\lam\in \RR}$ of $\lam$-densities which is unique up to a proportionality and
\newline
(a) $\Theta^\lam$ never vanishes;
\newline
(b) $\Theta^\lam$ is parallel with respect to the connection from Proposition \ref{P:connection};
\newline
(c) $(\Theta^\lam)^\mu =\Theta^{\lam\mu}$ for any $\lam,\mu\in \RR$.

\end{paragraphlist}

\section{Octonionic analogue of K\"ahler metrics.}\label{S:metrics-octon}
\begin{paragraphlist}

\item Let us introduce on $GL_2(\OO)$-affine manifolds a class of Riemannian metrics which can be considered as an octonionic analogue
of K\"ahler metrics on complex manifolds, or of Hessian metrics on affine manifolds, or and HyperK\"ahler with Torsion (HKT) metrics on hypercomplex manifolds.
The considerations are purely local, so we may and will replace our $GL_2(\OO)$-affine manifold
with a domain $\Ome\subset\OO^2$.

Let us introduce the first order differential operators on $\OO$-valued functions acting on a open subset $\Ome\subset \OO^2$ for $i=1,2$:
\begin{eqnarray*}
\overset{\to}{\pt_{\bar i}}F=\sum_{p=0}^7e_p\frac{\pt F}{\pt x_i^p},\\
F\overset{\leftarrow}{\pt_{\bar i}}=\sum_{p=0}^7\frac{\pt F}{\pt x_i^p}e_p.
\end{eqnarray*}

Note that the first operator coincides with $\frac{\pt}{\pt \bar q_i}$.

\begin{definition}\label{D:closed-current}
Let $T\colon \Ome\to \ch_2(\OO)$ be a $C^\infty$-smooth function.
We call $T$ a {\itshape closed smooth current} if it satisfies the following system of first order linear differential equations:
\begin{eqnarray}\label{E:system-equations}
T_{\bar 12}\overset{\leftarrow}{\pt}_{\bar 2}=\overset{\to}{\pt_{\bar 1}} T_{\bar 22},\,\,\, T_{\bar 2 1}\overset{\leftarrow}{\pt}_{\bar 1}=\overset{\to}{\pt_{\bar 2}} T_{\bar 11}
\end{eqnarray}
\end{definition}

The equations (\ref{D:closed-current}) are analogous to the definition of a closed smooth $(1,1)$-form in the case of 2 complex variables.

\item The main result of this section is
\begin{theorem}\label{T:octon-kahler}
Let $\Ome\subset \OO^2=\RR^{16}$ be an open convex subset. Let $T\colon \Ome\to \ch_2(\OO)$ be a $C^\infty$-smooth function. Then $T$ satisfies the system (\ref{E:system-equations})
if and only if there exists a $C^\infty$-smooth function $u\colon\Ome\to \RR$ such that
$$T=Hess_\OO(u).$$
\end{theorem}

Note that in this theorem $T$ defines a metric if and only if $T>0$, or equivalently $u$ is strictly octonionic plurisubharmonic in the sense of \cite{alesker-octon}.

\item A complex version of Theorem \ref{T:octon-kahler} is classical and is called local $dd^c$-lemma, see e.g. \cite{moroianu}, Proposition 8.8. A quaternionic version of this result was first obtained
by Michelson-Strominger \cite{michelson-strominger} for the flat quaternionic space $\HH^n$, and was proven by Banos-Swann \cite{banos-swann} for arbitrary hypercomplex manifolds.

\item As we will see, the proof of the theorem is a straightforward combination of two main ingredients: an Ehrenpreis' theorem which is a general result
from linear PDE with constant coefficients \cite{hormander}, Theorem 7.6.13, and a computation from commutative algebra which
we performed on a computer using Macauley2.

Let us describe the first mentioned ingredient. This is Theorem 7.6.13 in H\"ormander's book \cite{hormander}. According to this book, this result was first announced by
Ehrenpreis \cite{ehrenpreis} with further contributions by Malgrange \cite{malgrange} and Palamodov \cite{palamodov1}, \cite{palamodov2}.

Let $P_1,\dots, P_J$ be linear differential operators on $\RR^n$ with constant complex coefficients; thus they are complex polynomials in $\pt_1,\dots,\pt_n$.

\begin{theorem}\label{T:ehrenpreis}
Let $\Ome\subset\RR^n$ be an open convex subset. Let $f_1,\dots,f_j\colon \Ome\to \CC$ be $C^\infty$-smooth functions.
Then there exists a $C^\infty$-smooth function $u\colon \Ome\to \CC$ satisfying the system of equations
\begin{eqnarray}\label{E:u-general-solution}
P_ju=f_j, \mbox{ for all }j=1,\dots,J,
\end{eqnarray}
if and only if for any $J$-tuple of linear differential operatorsa with constant complex coefficient $Q_1,\dots,Q_J$ such that
\begin{eqnarray}\label{Q-commut}
\sum_{j=1}^J Q_jP_j=0
\end{eqnarray}
one has
\begin{eqnarray}\label{E-equation-gener-00}
\sum_{j=1}^J Q_jf_j=0.
\end{eqnarray}
\end{theorem}

\begin{remark}\label{R:real-ehrenpreis}
In our application to Theorem \ref{T:octon-kahler} it will be important to have a version of Theorem \ref{T:ehrenpreis} where all the differential operators $P_j$ have real coefficients,
and all the functions $f_j$ and $u$ are real valued. Then exactly the same theorem holds where the polynomials $Q_j$ also have real coefficients. This form is easily deduced from the above one,
in particular a complex solution $u$ is replaced by its real part which is also a solution (\ref{E:u-general-solution}).
\end{remark}

\item We are going to apply Theorem \ref{T:ehrenpreis} (taking into account Remark \ref{R:real-ehrenpreis}) to the system
\begin{eqnarray}\label{E:oct-system}
Hess_\OO(u)=T,
\end{eqnarray}
where $T\colon \Ome\to \ch_2(\OO)$. There are $n=16$ variables and $J=10$ (real) equations.

\item First let us observe that the ring of linear differential operators with constant real coefficients is identified with the ring of real polynomials via the substitution
$\pt_i\mapsto x_i$. In order to apply Theorem \ref{T:ehrenpreis} in a concrete situation one has to describe $J$-tuple of polynomials $(Q_1,\dots, Q_J)$ satisfying (\ref{Q-commut}).
Let us also observe that this set is a submodule of the free module of rank $J$ over the polynomial ring $A:=\RR[x_1,\dots,x_n]$ (now we work with real coefficients).
It is well known that every submodule of a finitely generated $A$-module is finitely generated (see e.g. Proposition 6.5 and Corollary 7.6 in \cite{atiyah-mcdonald}).
Clearly it suffices then to verify (\ref{E-equation-gener-00})
only on generators of that submodule.

The computer software Macaulay2 allows to compute easily generators of such a module. There is a special command 'kernel' which is applied to a matrix whose entries are polynomials with rational coefficients.
This matrix is considered a morphism between free modules over the polynomial ring. The output of the 'kernel' command is generators of the kernel of this morphism.

This computer computation leads to system (\ref{E:system-equations}). We discuss implementation of this procedure and provide a printout of Macaulay2 program in the Appendix.
Let us notice that Macaulay2 works with polynomials with rational (rather than real or complex) coefficients and computes
generators of a submodule over the ring $\QQ[x_1,\dots,x_n]$. In our situation the polynomials $P_1,\dots,P_J$ have indeed rational coefficients. One has
the following easy lemma which implies that it suffices to find generators of a submodule over polynomials with rational coefficients.
\begin{lemma}\label{L:commut-alg}
Let $P_1,\dots,P_J\in \QQ[x_1,\dots,x_n]$. Denote
\begin{eqnarray*}
N_\RR=\{(Q_1,\dots,Q_J)|\,\, Q_j\in \RR[x_1,\dots,x_J]\, \& \,\sum_{j=1}^J Q_jP_j=0\},\\
N_\QQ=\{(Q_1,\dots,Q_J)|\,\, Q_j\in \QQ[x_1,\dots,x_J]\, \&\, \sum_{j=1}^J Q_jP_j=0\}.
\end{eqnarray*}
Then $$N_\RR=\RR\cdot N_\QQ:=\{\sum_\alp \lam_\alp (Q_1^\alp,\dots,Q_J^\alp)|\, \lam_\alp\in \RR \, \&\,  (Q_1^\alp,\dots,Q_J^\alp)\in N_\QQ\}.$$
\end{lemma}
The lemma follows immediately from existence of a basis of $\RR$ over $\QQ$.

\item\label{item:new-metrics} Let us introduce a class of Riemannian metrics on  $GL_2(\OO)$-affine manifolds, called octonionic K\"ahler metrics, analogous to K\"ahler metrics on complex manifold, Hessian metrics on affine manifolds,
and HKT-metrics on hypercomplex manifolds.

\begin{definition}\label{D:o-hermit-metric}
A Riemannian metric on a $GL_2(\OO)$-affine manifold $M$ is called octonionic Hermitian if for any point $x\in M$ and some (equivalently, any) chart containing $x$ the restriction of the
metric to any octonionic line contained in the tangent space $T_xM\simeq \OO^2$ (the latter identification uses the chart)  is proportional to the standard Euclidean metric on $\OO^2$.
\end{definition}

Note that independence of this definition of a chart follows from Lemma \ref{m-sch2}.  Octonionic Hermitian metrics can be described explicitly
in a chart as follows. Let $\Ome\subset \OO^2$ be an open subset.
Let $T=(T_{\bar i j})_{i,j=1}^2$ be a $C^\infty$-smooth function on $\Ome$ with values in the space $\ch_2(\OO)$ of
$2\times 2$ octonionic Hermitian   matrices. Assume that $T>0$ pointwise in $\Ome$. Thus $T$ defines a Riemannian metric $g$ in $\Ome$ as follows
$$|\xi|^2_g=\sum_{i,j=1}^2Re(\bar \xi_i T_{\bar ij}\xi_j)$$
for any $\xi\in \OO^2$. By Corollary \ref{Cor:isomorphisms} this metric $g$ is octonionic Hermitian, and any octonionic Hermitian metric is given by this formula for a unique $T$ as above.

\begin{definition}\label{D:oct-kahler}
An octonionic Hermitian metric $g$ on $M$ is called octonionic K\"ahler if in some (equivalently, any) chart the $\ch_2(\OO)$-valued function $T$ corresponding to $g$ satisfies the system of equations (\ref{E:system-equations}).
\end{definition}
\begin{remark}
The independence of the last definition of a chart can be proved by showing that the system (\ref{E:system-equations}) is invariant under $GL_2(\OO)$ transformations. Alternatively,
by Theorem \ref{T:octon-kahler} $g$ is octonionic K\"ahler if and only if locally $T=Hess_\OO(u)$. But the latter equation is invariant under $GL_2(\OO)$ transformations by Proposition \ref{P:ohs1}.
\end{remark}

\end{paragraphlist}

\section{Integration by parts in $\OO^2$.}\label{S:integration}
\begin{paragraphlist}

\item\label{item:integr-1} In this section we prove a few identities for integrals over $SL_2(\OO)$-affine manifolds to be used later. We start with necessary results from linear algebra.

  Let $A=\left[\begin{array}{cc}
             a_{\bar 11}&a_{\bar 12}\\
             a_{\bar 21}&a_{\bar 22}
             \end{array}\right]$ be an octonionic Hermitian matrix. Like in the commutative case we denote by
             $$adj(A):=\left[\begin{array}{cc}
                             a_{\bar 22}&-a_{\bar 12}\\
                             -a_{\bar 21}&a_{\bar 11}
                              \end{array}\right]$$
                              its adjoint matrix. It is also octonionic Hermitian. One can verify by a direct computation that like in the commutative case if $A$ is invertible then
                              $$adj(A)=(\det A) \cdot A^{-1}.$$
It is easy to see that if $A\geq 0$ then $adj(A)\geq 0$.
The following elementary identity is proved by a straightforward computation:
\begin{eqnarray}\label{mixed-adj}
\det(A,B)=\frac{1}{2}Re(Tr(adj(A)\cdot B))
\end{eqnarray}
for any  $A,B\in\ch_2(\OO)$, where $Tr$ denotes the sum of the diagonal elements.

\item For any $2\times 2$ octonionic matrices $X,Y$ one clearly has
\begin{eqnarray}\label{o2-01}
Re(Tr(X))=Re(Tr(X^*)),\\\label{o2-02}
Re(Tr(XY))=Re(Tr(YX)).
\end{eqnarray}

Let $C\in \ch_2(\OO)$ and let $B$ be any $2\times 2$ octonionic matrix. Then
\begin{eqnarray*}
Re(Tr(CB))\overset{(\ref{o2-01})}{=} Re(Tr(B^*C))\overset{(\ref{o2-02})}{=}Re(Tr(CB^*)).
\end{eqnarray*}
Hence $Re(Tr(CB))=Re(Tr(C\frac{(B+B^*)}{2}))$. This and (\ref{mixed-adj}) imply
\begin{eqnarray}\label{mixed-adj-2}
\det(A,\frac{B+B^*}{2})=\frac{1}{2}Re(Tr(adj(A)\cdot B)).
\end{eqnarray}

\item The following proposition is a version of integration by parts. Here and below we denote
$$\Delta_i:=\sum_{a=0}^7\pt^2_{x_i^a}, \, i=1,2.$$
\begin{proposition}\label{P:integration-parts}
Let $u,v,f\colon \OO^2\to \RR$ be $C^\infty$-smooth functions, and assume that one of them is compactly supported.
Then
$$\int_{\OO^2} \det(Hess_\OO(u),Hess_\OO(v))f dx=- \frac{1}{2} Re\left(\int_{\OO^2}\left([f_1,f_2]adj(Hess_\OO(v))\left[\begin{array}{c}
                                                                                                                   u_{\bar 1}\\
                                                                                                                   u_{\bar 2}
                                                                                                                  \end{array}\right]\right)dx\right).$$
\end{proposition}
{\bf Proof.} Let us rewrite the left hand side in the proposition multiplied by 2 as follows:
\begin{eqnarray*}
\int Re\left(u_{\bar 11}v_{\bar 22}+u_{\bar 22}v_{\bar 11}-(u_{\bar 12}v_{\bar 21}+u_{\bar 21} v_{\bar 12})\right)f dx=\\
-\int\left(Re\{u_{\bar 1}\Delta_2 v_1\}+Re\{u_{\bar 2}\Delta_1v_2\}-Re\{u_{\bar 1}\underset{=\Delta_2 v_1}{\underbrace{\overset{\to}{\pt_2}(v_{\bar 21})}}\}-Re\{u_{\bar 2}\underset{=\Delta_1 v_2}{\underbrace{\overset{\to}{\pt_1}(v_{\bar 12})}}\}\right)fdx-\\
-\int \left(Re\{u_{\bar 1}f_1v_{\bar 22}\}+Re\{u_{\bar 2}f_2 v_{\bar 11}\}-Re\{u_{\bar 1}f_2v_{\bar 21}\}-Re\{u_{\bar 2}f_1v_{\bar 12}\}\right)dx=\\
-\int \left(Re\{f_1v_{\bar 22}u_{\bar 1}\}+Re\{f_2 v_{\bar 11}u_{\bar 2}\}-Re\{f_2v_{\bar 21}u_{\bar 1}\}-Re\{f_1v_{\bar 12}u_{\bar 2}\}\right)dx=\\
-\int Re\left( [f_1,f_2]\left[\begin{array}{cc}
                              v_{\bar 22}&-v_{\bar 12}\\
                              -v_{\bar 21}&v_{\bar 11}
                              \end{array}\right]\left[\begin{array}{c}
                                                       u_{\bar 1}\\
                                                       u_{\bar 2}
                                                       \end{array}\right]\right)dx.
\end{eqnarray*}
The result follows. \qed

\item We are going now to rewrite the statement of Proposition \ref{P:integration-parts} in a different language so that it will make sense on any $SL_2(\OO)$-affine manifold.

First we have
\begin{eqnarray}\label{E:pp1}
Re\left([f_1,f_2]adj(Hess_\OO(v))\left[\begin{array}{c}
                                   u_{\bar 1}\\
                                   u_{\bar 2}
                                   \end{array}\right]\right)=\\
                                   Re Tr\left(adj(Hess_\OO(v))\left[\begin{array}{c}
                                   u_{\bar 1}\\
                                   u_{\bar 2}
                                   \end{array}\right][f_1,f_2]\right)\overset{(\ref{mixed-adj-2})}{=}\\\label{E:pp2}
2\det\left(Hess_\OO(v),\ct\left(\left[\begin{array}{c}
                                                                      u_{\bar 1}\\
                                                                      u_{\bar 2}
                                                                     \end{array}\right]\otimes \left[\begin{array}{c}
                                                                      f_{\bar 1}\\
                                                                      f_{\bar 2}
                                                                     \end{array}\right]\right)\right),
\end{eqnarray}
where $\ct\colon Sym_\RR^2(\OO^2)\to \ch_2(\OO)$ is the map from Lemma \ref{m-sch1}.

Note that the differential $du\in (\RR^{16})^*$ is identified with $\left[\begin{array}{c}
                                                                      u_{\bar 1}\\
                                                                      u_{\bar 2}
                                                                     \end{array}\right]$ when $\RR^{16}$ is identified with $\OO^2$ in the usual sense, i.e.
                                                                     $$(x_0,\dots,x_7,y_0,\dots,y_7)\mapsto (\sum_{p=0}^7x_pe_p,\sum_{p=0}^7y_pe_p),$$
                                                                     and the dual space $(\RR^{16})^*$ is identified with $\RR^{16}$ in the standard coordinate way.

Hence, under these identifications, (\ref{E:pp2}) is equal to $2\det\left(Hess_\OO(v),\ct(du\otimes df)\right)$.
Thus we obtain
\begin{eqnarray}\label{E:integr-by-parts-octon}
Re\left([f_1,f_2]adj(Hess_\OO(v))\left[\begin{array}{c}
                                   u_{\bar 1}\\
                                   u_{\bar 2}
                                   \end{array}\right]\right)=2\det\left(Hess_\OO(v),\ct(du\otimes df)\right).
\end{eqnarray}
Note that the left hand side of (\ref{E:integr-by-parts-octon}) is defined in local coordinates, while the right hand side is defined globally on $M$.

This identification allows us to generalize Proposition \ref{P:integration-parts}
to arbitrary $SL_{2}(\OO)$-manifolds.
\begin{proposition}\label{P:integration-parts2}
Let $M^{16}$ be an $SL_2(\OO)$-affine manifold with a parallel non-vanishing 3/4-density $\Theta$. Let $u,v,f\colon M\to \RR$ be $C^\infty$-smooth functions, and assume that one of them is compactly supported.
Then
 \begin{eqnarray*}
\mbox{(a) }\int_{M} f\det(Hess_\OO(u),Hess_\OO(v))\cdot \Theta =-\int_M \det\left(Hess_\OO(v),\ct(du\otimes df)\right)\cdot \Theta,
\end{eqnarray*}
where now $\ct\colon Sym_\RR^2(T^*M)\to \underline{\ch}_2(M)$ is the map defined in Section \ref{S:G-mflds}, paragraph \ref{item1};
\newline
(b) $\int_M f\det(u,v)\cdot \Theta=\int_M u\det(f,v)\cdot \Theta.$
\end{proposition}
Part (a) follows from \ref{P:integration-parts} using partition of unity. Part (b) immediately follows from part (a).

\begin{corollary}\label{CoR:cor-by-parts}
Let $M^{16}$ be an $SL_2(\OO)$-affine manifold with a parallel non-vanishing 3/4-density $\Theta$. Let $u,f\colon \OO^2\to \RR$ be $C^\infty$-smooth functions, and assume that one of them is compactly supported.
Let $G_0$ be a $C^\infty$-smooth section of $\underline{\ch}_2(M)$ satisfying (\ref{E:system-equations}). Then,
$$\mbox{(a) }\int_{M} f\det(Hess_\OO(u),G_0)\cdot \Theta =-\int_M \det\left(G_0,\ct(du\otimes df)\right)\cdot \Theta,$$
$$\mbox{(b) }\int_M f\det (Hess_\OO(u),G_0)\cdot \Theta=\int_M u\det(Hess_\OO(f),G_0)\cdot \Theta.$$
\end{corollary}
{\bf Proof.} Clearly (b) follows from (a) since the right hand side of (a) is symmetric in $u$ and $f$.
To prove (a) observe that there exists a locally finite open covering $\{U_\alp\}$ by charts with a subordinate to it partition of unity $\{\chi_\alp\}$ such that
$$G_0=Hess_\OO(v_\alp) \mbox{ on } U_\alp,$$
where $v_\alp\colon U_\alp\to \RR$ is a $C^\infty$-smooth function. We have
\begin{eqnarray*}
\int_M f\det(Hess_\OO(v),G_0)\cdot \Theta=\\
\sum_\alp\int_{U_\alp}f\chi_\alp\det(Hess_\OO(u),Hess_\OO(v_\alp))\cdot\Theta\overset{\mbox{Prop. } \ref{P:integration-parts2}(a)}{=}\\
-\sum_\alp\int_{U_\alp}\det(Hess_\OO(v_\alp), \ct(du\otimes d(f\chi_\alp)))\cdot \Theta=\\
-\sum_\alp\int_{U_\alp}\det(G_0, \ct(du\otimes d(f\chi_\alp)))\cdot \Theta=\\
-\int_M\det(G_0, \ct(du\otimes df))\cdot \Theta.
\end{eqnarray*}
\qed

\end{paragraphlist}


\section{The Monge-Amp\`ere equation and the main result.}\label{S:MA-main}
\begin{paragraphlist}

\item
 Let $M$ be a $GL_2(\OO)$-affine manifold. Using the language described in Section \ref{S:G-mflds} we can interpret the Monge-Amp\`ere equation globally.
Let $G_0$ be a $C^\infty$-smooth section of the bundle $\uch$. Let $\phi,f$ be $C^\infty$-smooth functions. The  Monge-Amp\`ere equation which is written in local coordinates as
$$\det(G_0+Hess_\OO(\phi))=e^f\det(G_0),$$
is well defined globally when both sides are interpreted as $1/4$-densities.

 To prove solvability of the last equation we will assume that $G_0$ satisfies (\ref{D:closed-current}).
 By Theorem \ref{T:octon-kahler} this is equivalent to the following condition:
for every point $p\in M$ there exists a neighborhood $U$
and a $C^\infty$-function $u$ such that
\begin{eqnarray}\label{E:condition}
G_0(x)=Hess_\OO u(x) \mbox{ for any } x\in U.
\end{eqnarray}

\item Let us state a general conjecture a special case of which we will prove in this paper. This conjecture is obviously motivated by the Calabi-Yau theorem for K\"ahler manifolds \cite{yau}.
\begin{conjecture}\label{CY-theorem}
Let $M^{16}$ be a compact $GL_2(\OO)$-affine manifold. Let $G_0$ be a $C^\infty$-smooth positive section of the bundle $\uch$ satisfying condition (\ref{E:condition}). Let $f\colon M\to \RR$ be a $C^\infty$-smooth function.
Then there exist a $C^\infty$-smooth function $\phi\colon M\to \RR$
and a constant $A>0$ such that
 \begin{eqnarray}\label{E:conj}
 \det(G_0+Hess_\OO(\phi))=A e^f\det(G_0) \mbox{ on } M.
 \end{eqnarray}
\end{conjecture}

\item Let us assume that $M$ is an $SL_2(\OO)$-affine manifold rather than $GL_2(\OO)$. Then by the results of Section \ref{S:G-mflds}, paragraph \ref{item-Theta},
there exists a parallel and non-vanishing $3/4$-density $\Theta^{3/4}$, which will be denoted just $\Theta$ in the rest of the paper.
\begin{lemma}\label{L:const-A}
Let $M^{16}$ be a compact $SL_2(\OO)$-affine manifold with a non-vanishing parallel $3/4$-density $\Theta$. Let a $C^4$-smooth function $\phi$ and a constant $A$ satisfy (\ref{E:conj}). Then,
$$A=\frac{\int_M \det G_0\cdot \Theta}{\int_M e^f\det G_0\cdot \Theta}.$$
\end{lemma}
{\bf Proof.} Let us multiply both sides of (\ref{E:conj}) by $\Theta$ and integrate over $M$. Using the assumption (\ref{E:condition}) and Corollary \ref{CoR:cor-by-parts} we get in the left hand side
\begin{eqnarray*}
\int_M\big(\det G_0+2\det(G_0,Hess_\OO\phi)+\det(Hess_\OO\phi,Hess_\OO\phi)\big)\cdot\Theta=\\
\int_M\big(\det G_0+\phi \det(G_0,Hess_\OO 1)+\phi  \det(Hess_\OO\phi,Hess_\OO 1)\big)\cdot\Theta=\\
\int_M\det G_0\cdot \Theta.
\end{eqnarray*}
In the right hand side we have $A\int_M e^f \det G_0\cdot \Theta.$ \qed

\item From now on dealing with $SL_2(\OO)$-affine manifolds we will replace $f$ with $f+\ln A$ thus getting the equation
\begin{eqnarray}\label{E:ma-2}
\det(G_0+Hess_\OO(\phi))= e^f\det(G_0),
\end{eqnarray}
where
\begin{eqnarray}\label{E:normalization}
\int_M (e^f-1) \det G_0\cdot\Theta=0.
\end{eqnarray}

\item Let us formulate the main result of this paper. Note first that a $Spin(9)$-affine manifold is $SL_2(\OO)$-affine manifold, and hence it carries a non-vanishing parallel $3/4$-density $\Theta$.
\begin{theorem}\label{T:main-result}
Let $M^{16}$ be a compact $Spin(9)$-affine manifold. Let $f\colon M\to \RR$ be a $C^\infty$-smooth function. Let $G_0$ be a $C^\infty$-smooth positive section of $\uch$ satisfying (\ref{E:condition}).
Assume that $f,G_0$ satisfy (\ref{E:normalization}). Then there exists a unique (up to an additive constant)
$C^\infty$-smooth function $\phi\colon M\to \RR$ satisfying the Monge-Amp\`ere equation (\ref{E:ma-2}).
\end{theorem}

\end{paragraphlist}

\section{Uniqueness of solution.}\label{S:uniqueness}
\begin{theorem}\label{T:uniqueness-thm} Let $M^{16}$ be a compact connected $GL_2(\OO)$-affine manifold. Let $G_0$ be a continuous positive section of $\underline{\ch}_2(M)$. Let $f\colon M\to \RR$ be a continuous function.
Let $\phi_1,\phi_2\colon M\to \RR$ be $C^2$-smooth functions satisfying
$$\det(G_0+Hess_\OO(\phi_{p}))=e^f\det(G_0) \mbox{ where } p=1,2.$$
Then $\phi_1-\phi_2$ is a constant.
\end{theorem}
We start with the following lemma.
\begin{lemma}\label{L:ellipticity}
Let $M$ be a compact connected $GL_2(\OO)$-affine manifold. Let $G$ be a continuous section of $\uch$ such that $\det G>0$ everywhere on $M$, and $G>0$ at least at one point.
Then $G>0$ everywhere.
\end{lemma}
{\bf Proof.} Set $\cu:=\{z\in M| G(u)>0\}$. By assumption $\cu$ is non-empty. Clearly $\cu$ is open. It suffices to show that $\cu$ is closed. Let $z\in \bar\cu$. Then $G(z)\geq 0$.
If $G(z)$ is not positive definite then $\det G(z)=0$ contradicting the assumption. To show the latter statement note that we may assume the $G(z)$ is a complex Hermitian matric since its entries
are contained in a subfield of $\OO$ isomorphic to $\CC$. Then $G(z)\geq 0$ but not positive definite also as a complex Hermitian matrix by the Sylvester criterion (Proposition \ref{P:sylvester})).
For complex matrices that implies that $\det G(z)=0$. \qed

\hfill

{\bf Proof of Theorem \ref{T:uniqueness-thm}.} Subtracting the equations we get
\begin{eqnarray*}
\det(Hess_\OO(\phi_1-\phi_2), (G_0+Hess_\OO(\phi_1))+(G_0+Hess_\OO(\phi_2)))=0.
\end{eqnarray*}
 Since $G_0+Hess_\OO(\phi_{1,2})>0$ by Lemma \ref{L:ellipticity}, the function $\phi_1-\phi_2$ satisfies a second order linear elliptic equation with no free term and with continuous coefficients.
By the strong maximum principle
(see \cite{gilbarg-trudinger}, Theorem 3.5) this function has to be constant. \qed

\section{Zero order estimate.}\label{S:zero-order}
The main result of this section is Theorem \ref{T:l-infty} which is a $C^0$-estimate of a solution of the Monge-Amp\`ere
equation (\ref{E:ma-2}). The method of the proof is close to the original one due to Yau \cite{yau} with further modifications explained in  \cite{joyce-book}.

Throughout this section we assume the following:  $M^{16}$ is a compact connected $SL_2(\OO)$-affine manifold with a positive parallel $3/4$-density $\Theta$;
 $f\in C^\infty(M,\RR)$;
 $G_0$ is a $C^\infty$-smooth section of the bundle $\uch$ satisfying condition (\ref{E:condition}). In addition assume that the condition (\ref{E:normalization}) is satisfied
and $\phi$ is a $C^2$-smooth solution of (\ref{E:ma-2}).

\begin{lemma}\label{ze-1}
The solution $\phi$ satisfies the following estimates
\newline (a)
$||\nabla|\phi|^{p/2}||_{L^2}^2\leq
c_1\cdot\frac{p^2}{(p-1)} \int_M(1-e^f)\phi
|\phi|^{p-2}\det G_0\cdot \Theta\mbox{ for }p>2,$,
\newline (b)$
||\nabla\phi||_{L^2}^2\leq
c_1 \int_M(1-e^f)\phi
\det G_0\cdot \Theta,$
\newline where the constant $c_1$ depends on $M, G_0$, and $\Theta$ only.
\end{lemma}
{\bf Proof.} Let us fix a finite triangulation $\{\cs_\alp\}_\alp$ of $M$ into simplices such that each simplex is contained in a coordinate chart.
Note that $$(\phi|\phi|^{p-2})_{x_i^a}=(p-1)\phi_{x_i^a}|\phi|^{p-2}.$$
We have
\begin{eqnarray*}
\int_M (1-e^f)\phi|\phi|^{p-2}\det G_0\cdot \Theta=\\
\int_M \phi|\phi|^{p-2}(\det G_0-\det(G_0+Hess_\OO\phi))\cdot \Theta=\\
-\int_M\phi|\phi|^{p-2}\det(2G_0+Hess_\OO\phi,Hess_\OO\phi)\cdot \Theta\overset{\mbox{Corol. \ref{CoR:cor-by-parts}(a)}}{=}\\
(p-1)\int_M |\phi|^{p-2}\det\left(2G_0+Hess_\OO(\phi),\ct(d\phi\otimes d\phi)\right)\cdot \Theta\overset{(\ref{E:integr-by-parts-octon})}{=}\\
\frac{(p-1)}{2}\sum_\alp\int_{\cs_\alp} |\phi|^{p-2}\left(\left[\phi_1,\phi_2\right](adj(2G_0+Hess_\OO(\phi)))\left[\begin{array}{c}
                                                                 \phi_{\bar 1}\\
                                                                 \phi_{\bar 2}
                                                                   \end{array}\right]\right)(dq)_\alp,
\end{eqnarray*}
where $(dq)_\alp$ denotes a Lebesgue measure on a coordinate chart containing $\cs_\alp$; its normalization depends on the restriction of $\Theta$ to $\cs_\alp$.
Since for $2\times 2$ matrices $adj(A+B)=adj(A)+adj(B)$, the last expression can be rewritten as
\begin{eqnarray*}
\frac{(p-1)}{2}\sum_\alp\int_{\cs_\alp} |\phi|^{p-2}\left(\left[\phi_1,\phi_2\right](adj( G_0)+\underset{\geq 0}{\underbrace{adj(G_0+Hess_\OO(\phi))}})\left[\begin{array}{c}
                                                                 \phi_{\bar 1}\\\label{E:ddd3}
                                                                 \phi_{\bar 2}
                                                                   \end{array}\right]\right)(dq)_\alp.
\end{eqnarray*}
Than last expression is greater than or equal to
\begin{eqnarray}\label{E:ddd5}
 \frac{(p-1)}{2}\sum_\alp\int_{\cs_\alp} |\phi|^{p-2}\left(\left[\phi_1,\phi_2\right]adj( G_0)\left[\begin{array}{c}
  \phi_{\bar 1}\\
  \phi_{\bar 2}
  \end{array}\right]\right)(dq)_\alp.
   \end{eqnarray}
Taking  $p=2$ in (\ref{E:ddd5}) we have
\begin{eqnarray*}
\frac{1}{2}\sum_\alp\int_{\cs_\alp} \left(\left[\phi_1,\phi_2\right]adj( G_0)\left[\begin{array}{c}
  \phi_{\bar 1}\\
  \phi_{\bar 2}
  \end{array}\right]\right)(dq)_\alp\geq const ||\nabla\phi||^2_{L^2},
\end{eqnarray*}
where $const$ depends only on $M,G_0,\Theta$, and the triangulation only.

Now let us take $p>2$ in (\ref{E:ddd5}). Since
$$(|\phi|^{p/2})_{x_i^a}=\pm\frac{p}{2} \phi_{x_i^a}|\phi|^{p/2-1}$$
pointwise, we get
\begin{eqnarray*}
(\ref{E:ddd5})=
\frac{2(p-1)}{p^2}   \sum_\alp\int_{\cs_\alp} \left(\left[(|\phi|^{p/2})_1,(|\phi|^{p/2})_2\right]adj(G_0)\left[\begin{array}{c}
                                                                 (|\phi|^{p/2})_{\bar 1}\\
                                                                 (|\phi|^{p/2})_{\bar 2}
                                                                   \end{array}\right]\right) (dq)_\alp \geq\\
const       \frac{(p-1)}{p^2} ||\nabla(|\phi|^{p/2})||_{L^2}^2,
\end{eqnarray*}
where $const$ again depends only on $M,G_0,\Theta$, and the triangulation only. \qed

\begin{lemma}\label{ze-3}
There exists a constant $c_2$ depending on $M, G_0$, and $\Theta$ only,
such that for any function $\psi$ from the Sobolev space $ W^{1,2}(M)$,
$$||\psi||^2_{L^{\frac{16}{7}}}\leq c_2(||\nabla
\psi||_{L^2}^2+||\psi||^2_{L^2}).$$ Moreover, if  $\psi$ satisfies
$\int_M\psi \det G_0\cdot \Theta=0$, one has
$$||\psi||^2_{L^{\frac{16}{7}}}\leq c_2||\nabla \psi||_{L^2}^2.$$
\end{lemma}

{\bf Proof:} By the Sobolev imbedding theorem there exists a
constant $C'$ such that, for any function $\psi\in W^{1,2}(M)$, one
has
$$||\psi||^2_{L^{\frac{16}{7}}}\leq C'(||\nabla
\psi||^2_{L^2}+||\psi||^2_{L^2}).$$
That gives the first inequality in the lemma.

If the function $\psi$
satisfies $\int_M\psi\det G_0\cdot \Theta=0$ then one has
$||\psi||_{L^2}^2\leq \tilde C ||\nabla \psi||_{L^2}^2$ since the
second eigenvalue of the Laplacian on $M$ is strictly positive. This gives
Lemma \ref{ze-3} proven. \qed

\begin{lemma}\label{ze-4}
Let us use the notation of the beginning of this section and let us assume again that
$$\int_M \phi\det G_0\cdot \Theta=0.$$
Then  there exists a constant $c_3$ depending on $M,G_0,\Theta$ and
$||f||_{C^0}$ only such that for any $p\in[2,\frac{16}{7}]$ one has
$||\phi||_{L^p}\leq c_3$.
\end{lemma}

{\bf Proof:} By Lemma \ref{ze-1}(b) we get
\begin{eqnarray*}
 ||\nabla \phi||^2_{L^2}\leq &  const
 (1+\exp({||f||_{C^0}}))||\phi||_{L^1}\\ \leq & const'\cdot
(1+ \exp({||f||_{C^0}}))||\phi||_{L^2}
\end{eqnarray*}
where the second
 inequality follows from the H\"older inequality. Since
\[ \int_M\phi\det  G_0\cdot \Theta=0,\] we have
 $$||\phi||_{L^2}\leq c_2||\nabla \phi||_{L^2}.$$
Hence $||\nabla \phi||_{L^2}\leq c_2 const'\cdot
(1+ \exp({||f||_{C^0}}))$. Therefore by Lemma \ref{ze-3}
there exists a constant $c_3'$ depending on $M,G_0,\Theta$, and
$||f||_{C^0}$ only such that $$||\phi||_{L^{\frac{16}{7}}} \leq
c_3'.$$ Hence by the H\"older inequality $||\phi||_{L^{p}}\leq
c_3$ for $p\in[2,\frac{16}{7}]$. \qed


\begin{lemma}\label{ze-5}
There exist positive constants $Q,c_4$ depending on $M,$ $G_0,\Theta$, and
$||f||_{C^0}$ only such that for any $p\geq 2$
\[ ||\phi||_{L^p}\leq Q(c_4p)^{-\frac{8}{p}}.\]
\end{lemma}

{\bf Proof.}  Obviously there exists a constant $c_5>0$ depending only on $M,G_0,\Theta$ such that
\begin{eqnarray}\label{E:c-5-def}
(\int_M \det G_0 \cdot \Theta)^{\frac{1}{p}}<c_5 \mbox{ for all } p\geq 2.
\end{eqnarray}
Define $$c_4:=(8/7)^7c_2(4  c_1c_5   e^{||f||_{C^0}} +1)$$ where $c_1,c_2$ are from Lemmas \ref{ze-1}, \ref{ze-3}. Choose
$Q$ so that $Q>c_3(c_4p)^{\frac{8}{p}}$ for $2\leq p\leq
\frac{16}{7}$ and $Q>(c_4p)^{\frac{8}{p}}$ for $2\leq p<\infty$.

We will prove the result by induction on $p$. By Lemma \ref{ze-4}, if $2\leq
p\leq \frac{16}{7}$ then $||\phi||_{L^p}\leq c_3\leq
Q(c_4p)^{-\frac{8}{p}}$. For the inductive step, suppose that
$$||\phi||_{L^p}\leq Q(c_4p)^{-\frac{8}{p}} \mbox{ for } 2\leq p\leq
k,\, \mbox{ where } k\geq \frac{16}{7} \mbox{ is a real number}.$$

We will show that
$$||\phi||_{L^q}\leq Q(c_4q)^{-\frac{8}{q}} \mbox{ for } 2<q\leq\frac{16}{7}k,$$ and, therefore, by induction Lemma \ref{ze-5} will be
proven. Let $p\in(2,k]$. By Lemma \ref{ze-1}(a) we get
\begin{eqnarray}\label{ine1}
||\nabla|\phi|^{p/2}||_{L^2}^2\leq
c_1\frac{p^2}{(p-1)}(1+e^{||f||_{C^0}})||\phi||^{p-1}_{L^{p-1}}.
\end{eqnarray}

Applying Lemma \ref{ze-3} to $\psi=|\phi|^{p/2}$ we get
\begin{eqnarray}\label{ine2}
||\phi||^p_{L^{\frac{8}{7} p}}\leq
c_2(||\nabla|\phi|^{p/2}||^2_{L^2}+||\phi||^p_{L^p}).
\end{eqnarray}

Combining (\ref{ine2}) and (\ref{ine1}) we obtain
$$||\phi||^p_{L^{\frac{8}{7} p}}\leq
c_2(4p  c_1
e^{||f||_{C^0}}||\phi||^{p-1}_{L^{p-1}}+||\phi||^p_{L^p}).$$ Let
$q:=\frac{8}{7} p$. Since $2< p\leq k$ we have by the induction assumption
\begin{eqnarray}\label{E:0-0-1}
||\phi||_{L^p}\leq Q(c_4p)^{-\frac{8}{p}}.
\end{eqnarray}

By the H\"older inequality $||\phi||_{L^{p-1}}^{p-1}\leq
||\phi||_{L^p}^{p-1}(\int_M \det G_0 \cdot \Theta)^{\frac{1}{p}}\overset{(\ref{E:c-5-def})}{\leq} c_5 ||\phi||_{L^p}^{p-1}$.

Hence we get
\begin{eqnarray*}
||\phi||^p_{L^q}\leq  c_2\left(4p  c_1c_5
e^{||f||_{C^0}}||\phi||^{p-1}_{L^{p}}+
||\phi||_{L^p}^p\right)\leq\\
c_2\left(4p c_1c_5
e^{||f||_{C^0}}(Q(c_4p)^{-\frac{8}{p}})^{p-1}+\left(Q(c_4p)^{-\frac{8}{p}}\right)^p\right).
\end{eqnarray*}

We have  by the choice of $Q$ $$Q(c_4p)^{-\frac{8}{p}}\geq 1.$$  Hence,
$$||\phi||^p_{L^q}\leq c_2Q^p(c_4p)^{-8}(4p  c_1c_5
e^{||f||_{C^0}} +1).$$ It remains to show that the last expression
is at most $Q^p(c_4q)^{-\frac{8}{q} p}.$ Equivalently it suffices to check
that $$c_2(c_4p)^{-8}(4p  c_1c_5  e^{||f||_{C^0}} +1)\leq
(c_4q)^{-7}=(\frac{8}{7}c_4 p)^{-7}.$$
In view that the left hand side is at most $c_2(c_4p)^{-8}p(4  c_1c_5
  e^{||f||_{C^0}} +1)$, it is enough to check that
$$c_2c_4^{-8}(4  c_1c_5   e^{||f||_{C^0}} +1)\leq
(\frac{8}{7}c_4 )^{-7}.$$ Namely,
$$c_4\geq (8/7)^7c_2(4  c_1c_5   e^{||f||_{C^0}} +1)$$ which holds by the definition of $c_4$.
\qed

\hfill

When $p\to +\infty$ we get the main result of this section.
\begin{theorem}\label{T:l-infty}
Let $M^{16}$ be a compact $SL_2(\OO)$-affine manifold with a positive parallel $3/4$-density $\Theta$. Let $G_0$ be a $C^\infty$-smooth positive definite section of the bundle $\uch$ satisfying condition (\ref{E:condition}).
Let $\phi$ be a $C^2$-smooth solution of the Monge-Amp\`ere equation (\ref{E:ma-2}) satisfying $\int_M\phi \det G_0\cdot \Theta=0$. Then there exists a constant $C$ depending on $M,G_0,\Theta$, and $||f||_{C^0}$ such that
$$||\phi||_{C^0}\leq C.$$
\end{theorem}

\section{Elementary inequality.}
The goal of this section is to prove Proposition \ref{P:elem-inequality} which will be used later on in the proof of the second order estimates.

 First recall that a $C^2$-smooth function $u\colon\Ome\to\RR$ in a domain $\Ome\subset\OO^2$ is called {\itshape strictly
octonionic plurisubharmonic} (strictly oPSH) provided $Hess_\OO(u)>0$ in $\Ome$ pointwise.
\begin{proposition}\label{P:elem-inequality}
Let $u\in C^4$ be a strictly oPSH function such that at
a given point $z$ its octonionic Hessian $(u_{\bar i j})$ is
diagonal. Then at this point $z$ one has
\begin{eqnarray}\label{C49}
\sum_{p=0}^7\sum_{i,k=1}^2\frac{|u_{\bar k kx^p_i}|^2}{u_{\bar i
i}u_{\bar k k}}\leq 4\sum_{p=0}^7\sum_{i,k,l=1}^2\frac{|u_{\bar k i
x^p_l}|^2}{u_{\bar i i}u_{\bar k k}}.
\end{eqnarray}
\end{proposition}
{\bf Proof.} Let us fix now the indices $i,k$ and compare the
summands containing this pair of indices in both sides of
(\ref{C49}).

First consider the case $i=k$. In the left hand side we have
\begin{eqnarray}\label{C50}
\sum_{p=0}^7\frac{|u_{\bar k kx^p_k}|^2}{u_{\bar k k}u_{\bar k k}}.
\end{eqnarray}
In the right hand side of (\ref{C49}) we have
\begin{eqnarray}\label{C51}
4\sum_{p=0}^7\sum_{l=1}^2\frac{|u_{\bar k kx^p_l}|^2}{u_{\bar k
k}u_{\bar k k}}
\end{eqnarray}
It is clear that $(\ref{C50})\leq (\ref{C51})$.

Let us consider the case now $i\ne k$. The left hand side of
(\ref{C49}) contains two summands with the pair $i,k$:
\begin{eqnarray}\label{C52}
\frac{1}{u_{\bar i i}u_{\bar k k}}\sum_p(|u_{\bar k
kx^p_i}|^2+|u_{\bar i ix^p_k}|^2).
\end{eqnarray}
The right hand side of (\ref{C49}) contains two summands with the
pair $i,k$:
\begin{eqnarray}\label{C53}
\frac{4}{u_{\bar i i}u_{\bar k k}}\sum_{p,l}(|u_{\bar k
ix^p_l}|^2+|u_{\bar i kx^p_l}|^2)=\frac{8}{u_{\bar i i}u_{\bar k
k}}\sum_{p,l}|u_{\bar k ix^p_l}|^2.
\end{eqnarray}
To finish the proof of proposition, it suffices now to show that
$(\ref{C52})\leq (\ref{C53})$. Explicitly, after cancelling out
the term $u_{\bar i i}u_{\bar k k}$ on both sides, it reduces to
\begin{eqnarray}\label{C54}
\sum_p(|u_{\bar k kx^p_i}|^2+|u_{\bar i ix^p_k}|^2)\leq
8\sum_{p,l}|u_{\bar k ix^p_l}|^2.
\end{eqnarray}
In order to show such a general inequality it suffices to sum up in
the right hand side over $l=i,k$ only. Thus (\ref{C54}) follows from
\begin{eqnarray}\label{C55}
\sum_p(|u_{\bar k kx^p_i}|^2+|u_{\bar i ix^p_k}|^2)\leq
8\sum_p(|u_{\bar k ix^p_i}|^2+|u_{\bar k ix^p_k}|^2).
\end{eqnarray}
In the last inequality we may separate summands containing
derivatives $kki$ and $kii$. These two inequalities are completely
symmetric and obtained one from the other by exchange $i$ by $k$.
Thus it suffices to show that
\begin{eqnarray}\label{C56}
\sum_p|u_{\bar k kx^p_i}|^2\leq 8\sum_p|u_{\bar k ix^p_k}|^2=8\sum_p|u_{\bar i kx_k^p}|^2.
\end{eqnarray}
Let us define two operators
$\overset{\leftarrow}{\pt_k},\overset{\leftarrow}{\pt_{\bar k}}$
acting on the space of octonionic valued functions:
\begin{eqnarray*}
\Phi\overset{\leftarrow}{\pt_k}:=\sum_{p=0}^7\frac{\pt\Phi}{\pt
x^p_k}\bar e_p,\\
\Phi\overset{\leftarrow}{\pt_{\bar k}}:=\sum_{p=0}^7\frac{\pt
\Phi}{\pt x^p_k} e_p,
\end{eqnarray*}
where $\bar e_p$ denotes the octonionic conjugate of the
octonionic unit $e_p$ (here $p=0,...,7$).

Recall that $\Delta_k=\sum_{p=0}^7\frac{\pt^2}{(\pt x^k_p)^2}$ is the
Laplacian with respect to the $k$-th octonionic variable. Note that for real valued function $u$ one has $\Delta_ku=u_{\bar k k}$.

Let us show that for an $\OO$-valued function $\Psi$ one has
\begin{eqnarray}\label{E:derivatives}
\Delta_k\Psi=(\Psi\overset{\leftarrow}{\pt_{\bar
k}})\overset{\leftarrow}{\pt_k}=(\Psi\overset{\leftarrow}{\pt_k})\overset{\leftarrow}{\pt_{\bar k}}.
\end{eqnarray}
By Lemma \ref{L:oalg1}(iii) for any octonions $a,b,c$ one has
$$(ca)b+(c\bar b)\bar a=c(ab+\bar b\bar a).$$

 Using this identity we have
\begin{eqnarray*}
(\Psi\overset{\leftarrow}{\pt_k})\overset{\leftarrow}{\pt_{\bar k}}=\sum_{p,q=0}^7 (\Psi_{x_k^px_k^q}\bar e_p) e_q=\frac{1}{2}\sum_{p,q=0}^7((\Psi_{x_k^px_k^q}\bar e_p) e_q+(\Psi_{x_k^px_k^q}\bar e_q) e_p)=\\
\frac{1}{2}\sum_{p,q=0}^7\Psi_{x_k^px_k^q}(\bar e_p e_q+\bar e_q e_p)=\\
\sum_p\Psi_{x_k^px_k^p}+\frac{1}{2}\sum_{p\ne q}\Psi_{x_k^p x_k^q}(\bar e_p e_q+\bar e_q e_p)=\Delta_k\Psi+\sum_{p\ne q}\Psi_{x_k^px_k^q} Re(\bar e_p e_q)=\\
\Delta_k \Psi.
\end{eqnarray*}
Thus (\ref{E:derivatives}) is proven.

Let us take $\Phi=u_{\bar i}(=\sum_qe_qu_{x^q_i})$. Then (\ref{C56})
can be rewritten
\begin{eqnarray}\label{C57}
|\Delta_k\Phi|^2\leq 8
\sum_p|\Phi_{x^p_k}\overset{\leftarrow}{\pt_k}|^2.
\end{eqnarray}
We have
$$\Delta_k\Phi=(\Phi\overset{\leftarrow}{\pt_k})\overset{\leftarrow}{\pt_{\bar k}}=
\sum_q(\Phi\overset{\leftarrow}{\pt_k})_{x^q_k} e_q.$$ Denote
$\Psi:=\Phi\overset{\leftarrow}{\pt_k}$. Then (\ref{C57}) takes the form:
\begin{eqnarray}\label{C58}
|\sum_q\Psi_{x^q_k} e_q|\leq \sqrt{8\sum_q|\Psi_{x^q_k}|^2}.
\end{eqnarray}

Indeed by the Cauchy-Schwartz inequality we have $$|\sum_{q=0}^7\Psi_{x^q_k}
e_q|\leq \sum_{q=0}^7 |\Psi_{x^q_k}|\leq \sqrt{8}\cdot
\sqrt{\sum_{q=0}^7|\Psi_{x^q_k}|^2}.$$ That completes the proof. \qed

\section{A few elementary identities.}

\begin{claim}\label{inverse}
Let $X\in \ch_2(\OO)$. Then there exists a unique $Y\in \ch_2(\OO)$ such that $XY=YX=I_2$ if and only if $\det X\ne 0$. (Such $Y$ is denoted by $X^{-1}$ as usual.)
\end{claim}
{\bf Proof.} The subalgebra generated by elements of $X$ is commutative, thus is isomorphic either to $\RR$ or $\CC$. Over these field $X$ is invertible, and the inverse matrix is necessarily Hermitian.

Conversely, let $XY=YX=I_2$. The subalgebra generated by elements of $X$ and $Y$ is associative, thus one may assume that it is isomorphic to $\HH$. Over $\HH$ the uniqueness of $X^{-1}$ is well known and proved
as usual. Thus elements of $Y$ belong to the commutative field generated by elements of $X$, in fact $Y=\frac{1}{\det X}adj(X)$. This proves the uniqueness of $Y$. \qed

\hfill

In local coordinates we denote by $f_a$ the derivative of a function $f$ with respect to a real direction $a$. Furthermore we denote
$$tr(X):=Re(Tr(X)),$$
where $Tr$ is the sum of diagonal elements as previously.
\begin{claim}\label{Cl:2}
Let $U$ be a smooth function with values in $\ch_2(\OO)$. Then,

(1) $(\log(\det U))_a=tr(U^{-1}U_a)$;

(2) $(U^{-1})_a=-U^{-1}U_a U^{-1}$.\footnote{The brackets in the right hand side are not important since the algebra generated by elements of two Hermitian matrices is associative.}
\end{claim}

{\bf Proof.} (1) We have
\begin{eqnarray*}
(\log(\det U))_a=(\det U)^{-1}\cdot 2\det(U,U_a)\overset{(\ref{mixed-adj})}{=}\\
(\det U)^{-1}\cdot tr(adj(U)U_a)=tr(U^{-1}U_a).
\end{eqnarray*}

(2) is elementary. \qed

\hfill

The proof of the following claim is also standard.

\begin{claim}\label{Cl:3}
Assume $\det U=F$. Then
\newline
(1) $tr(U^{-1}U_a)=(\log F)_a$.
\newline
(2) $tr(U^{-1}U_{ab})-tr((U^{-1}U_bU^{-1})U_a)=(\log F)_{ab}.$
\end{claim}
\begin{remark}
Given 3 matrices $X,Y,Z\in \ch_2(\OO)$, $tr(XYXZ)$ does not depend on bracketing inside. This follows from two facts: (1) $Re(abc)$ does not depend on bracketing for $a,b,c\in\OO$;
(2) $XYX$ does not depend on bracketing for $X,Y\in\ch_2(\OO)$.
\end{remark}

\section{Second order estimates.}\label{S:second-order}
\begin{paragraphlist}
\item Let $M$ be a compact connected $SL_2(\OO)$-affine manifold with a parallel positive $3/4$-density $\Theta$. Let $G_0$ be a $C^\infty$-smooth positive definite section of $\uch$ satisfying (\ref{E:condition}).
Also assume that $f\in C^\infty(M,\RR)$ and $\phi \in C^4(M,\RR)$ satisfy
\begin{eqnarray}\label{E:ma-second-order}
\det(G_0+Hess_\OO(\phi))=e^f\det(G_0).
\end{eqnarray}

In this section we make the following strong assumption on $M$: $M$ is a $Spin(9)$-affine manifold.
It follows that there exists a positive $C^\infty$-smooth section $G_{00}$ of $\underline{\ch}_2(M)$ which is flat with respect to the natural (corresponding to the $GL_2(\OO)$-structure) connection on $M$.

\item\label{item:DEF-Delta} Let us define the second order linear differential operator on real valued functions
$$\Delta f=tr(G_{00}^{-1}Hess_\OO(f)),$$
(recall that $tr=Re\,Tr$).
$\Delta$ is globally defined. Indeed, this operator is equal to the derivative of $U\mapsto \log\frac{\det U}{\Theta^{1/3}}$ at $U=G_{00}$ in the direction $Hess_\OO(f)$.

By Proposition \ref{P:diagonalization} locally there exist coordinates (consistent with the $GL_2(M)$- rather than the $Spin(9)$-structure on $M$) in which $G_{00}=Id$.
In such coordinates  the operator $\Delta$ is the usual flat Laplacian
$$\Delta f=f_{\bar 11}+f_{\bar 22}.$$

Let us denote $U:=G_0+Hess_\OO \phi$. Lemma \ref{L:ellipticity} implies that $U>0$. Define another linear second order differential operator
$$\Delta' f=Tr(U^{-1}Hess_\OO(f)).$$
Both $\Delta$ and $\Delta'$ are elliptic with vanishing free term.

\item Below we will make some local computations in coordinates system where $G_{00}=Id$ and $G_0+Hess_\OO \phi=U=Hess_\OO u$. Thus $u$ is strictly oPSH.
Clearly, $\Delta u=tr(U)$.

In such coordinate system we have
$$\Delta'(f)=\sum_{i,j=1}^2 Re((U^{-1})_{\bar ij}f_{\bar ji}).$$

\begin{proposition}\label{P:flat-double-Laplacian}
Let $u\in C^4$ be a strictly oPSH function. Let $\gamma$ be a smooth function of one real variable. Choose
(locally) coordinates such that $G_{00}\equiv Id$.  Then, in these
coordinates,
\begin{eqnarray*}
\,\,[\gamma(\Delta u)]_{ab}=\gamma''(\Delta u)\cdot tr(U_a)\cdot
tr(U_b)+\gamma'(\Delta u)tr(U_{ab}).
\end{eqnarray*}
If moreover at a point $z$ the matrix $U(z)=(u_{\bar i j}(z))$ is
diagonal then at this point $z$ one has
\begin{eqnarray*}
\Delta'(\gamma(\Delta u))=\gamma''(\Delta
u)\sum_{i,p}\frac{1}{u_{\bar i i}}(trU_{x^p_i})^2+\gamma'(\Delta
u)\sum_{i,k}\frac{u_{\bar i i \bar k k}}{u_{\bar i i}}.
\end{eqnarray*}
\end{proposition}
{\bf Proof.} The first part is straightforward. For the second part, using the first part, we have
\begin{eqnarray*}
\Delta'(\gamma(\Delta u))=\sum_i \frac{(\gamma(\Delta u))_{\bar ii}}{u_{\bar ii}}=\sum_{i,p}\frac{(\gamma(\Delta u))_{x_i^px_i^p}}{u_{\bar ii}}=\\
\gamma''(\Delta u)\sum_{i,p} \frac{(Tr(U_{x_i^p}))^2}{u_{\bar ii}}+\gamma'(\Delta u)\sum_{i,p} \sum_k \frac{(u_{\bar kk})_{x_i^px_i^p}}{u_{\bar ii}}=\\
\gamma''(\Delta u)\sum_{i,p} \frac{(Tr(U_{x_i^p}))^2}{u_{\bar ii}}+\gamma'(\Delta u)\sum_{i,k}  \frac{u_{ \bar ii\bar kk}}{u_{\bar ii}}
\end{eqnarray*}
\qed

\item
\begin{lemma}\label{L:4-order}
(1) Let $\det U=F$. Then one has\footnote{The bracketing in the expression below is not important.}
$$Tr(U^{-1}\Delta U)-\sum_{i,p}Tr(U^{-1}U_{x_i^p}U^{-1}U_{x_i^p})=\Delta (\log F).$$
\newline
(2) Assume that at a point $z$ the matrix $U(z)=(u_{\bar i j}(z))$ is
diagonal. Then, the equality above is equivalent to
$$\sum_{i,k}\frac{u_{\bar ii\bar kk}}{u_{\bar ii}}-\sum_{i,l,n,p}\frac{|u_{\bar lnx_i^p}|^2}{u_{\bar ll}u_{\bar nn}}=\Delta (\log F).$$
\end{lemma}
{\bf Proof.} Part (2) follows directly from part (1). Part (1) follows from Claim \ref{Cl:3}(2) by taking there $a=b$ corresponding to $x_i^p$ and summing over $i,p$. \qed

\hfill
\begin{proposition}\label{P:Laplace-twice}
Let $u\in C^4$ be strictly oPSH.
At any point $z\in M$ one has
$$\Delta'((\Delta u)^{3/4})\geq \frac{3}{4}(\Delta u)^{-1/4}\Delta (\log F),$$
where $F$ is as in Lemma \ref{L:4-order}.
\end{proposition}

{\bf Proof.} Fix $z\in M$. By Proposition \ref{P:diagonalization} we may and will choose coordinates near $z$ such that $G_{00}(z)=Id$ and $U(z)$ is diagonal.

Combining results of Lemma \ref{L:4-order} into Proposition \ref{P:flat-double-Laplacian} we get
\begin{eqnarray*}
\Delta'(\gamma(\Delta u))=\gamma''(\Delta
u)\sum_{i,p}\frac{1}{u_{\bar i i}}(TrU_{x_p^i})^2+
\gamma'(\Delta u)(\sum_{i,l,n,p}\frac{|u_{\bar lnx_i^p}|^2}{u_{\bar ll}u_{\bar nn}}+\Delta (\log F)).
\end{eqnarray*}
Equivalently
\begin{eqnarray}\label{E:01}
\Delta'(\gamma(\Delta u))=\gamma''(\Delta
u)\sum_{i,p}\frac{1}{u_{\bar i i}}(\sum_{k} u_{\bar kk x_i^p})^2+\gamma'(\Delta u)(\sum_{i,l,n,p}\frac{|u_{\bar lnx_i^p}|^2}{u_{\bar ll}u_{\bar nn}}+\Delta (\log F)).
\end{eqnarray}

Let us assume that $\gamma$ is concave, i.e. $\gamma''\leq 0$. Then,
\begin{eqnarray*}
\Delta'(\gamma(\Delta u))= \gamma''(\Delta
u)\sum_{i,p}\frac{1}{u_{\bar i i}}\left(\sum_k \sqrt{u_{\bar kk}}\frac{u_{\bar kk x_i^p}}{\sqrt{u_{\bar kk}}}\right)^2+\gamma'(\Delta u)(\sum_{i,l,n,p}\frac{|u_{\bar lnx_i^p}|^2}{u_{\bar ll}u_{\bar nn}}+\Delta (\log F))\geq\\
\gamma''(\Delta u)\sum_{i,p}\frac{1}{u_{\bar i i}}\cdot \Delta u\cdot\left( \sum_k\frac{ |u_{\bar kk x_i^p}|^2}{u_{\bar kk}}\right)+\gamma'(\Delta u)(\sum_{i,l,n,p}\frac{|u_{\bar lnx_i^p}|^2}{u_{\bar ll}u_{\bar nn}}+\Delta (\log F))=\\
\Delta u \cdot \gamma''(\Delta u)\sum_{i,k,p}\frac{ |u_{\bar kk x_i^p}|^2}{u_{\bar ii}u_{\bar kk}}+\gamma'(\Delta u)\left(\sum_{i,l,n,p}\frac{|u_{\bar lnx_i^p}|^2}{u_{\bar ll}u_{\bar nn}}+\Delta (\log F)\right).
\end{eqnarray*}

Let us finally take $\gamma(t)=t^{3/4}$. It is concave for $t>0$. Then we get
\begin{eqnarray*}
\Delta'((\Delta u)^{3/4})\geq \frac{3}{4}(\Delta u)^{-1/4}\left(\sum_{i,l,n,p}\frac{|u_{\bar lnx_i^p}|^2}{u_{\bar ll}u_{\bar nn}}-
\frac{1}{4}\sum_{i,k,p}\frac{|u_{\bar kkx_i^p}|^2}{u_{\bar ii}u_{\bar kk}}\right)+\frac{3}{4}(\Delta u)^{-1/4}\Delta (\log F).
\end{eqnarray*}
By Proposition \ref{P:elem-inequality} the first summand is non-negative. The result follows. \qed

\item Let us state the main result of this section.
\begin{theorem}\label{T:second-order}
Let $M^{16}$ be a compact connected $Spin(9)$-affine manifold. Let $G_{00}$ be a $C^\infty$-smooth positive definite section $\uch$ parallel with respect to the flat connection corresponding
to the affine structure of $M$.\footnote{Such a section $G_{00}$ does exist on any $Spin(9)$-affine manifold. It necessarily satisfies the condition (\ref{E:condition}).}
 Let $G_0$ be another $C^\infty$-smooth positive definite section of $\uch$ satisfying the condition (\ref{E:condition}). Assume also that $f\in C^2(M,\RR)$ and
$\phi\colon M\to \RR$ is a $C^4$-smooth
solution of the octonionic Monge-Amp\`ere equation
\begin{eqnarray}\label{E:monge-ampere}
\det(G_0+Hess_\HH\phi)=e^f\det G_0.
\end{eqnarray}
Then there exists a constant
$C$ depending on $M,G_{00},G_0$, and $||f||_{C^2(M)}$ such that
$$||\Delta\phi||_{C^0}\leq C.$$
\end{theorem}
{\bf Proof.} The statement of the theorem does not change if we add any constant to $\phi$.
Hence we may assume that $\int_M\phi\det G_0\cdot\Theta =0$. Then by Theorem \ref{T:l-infty} $||\phi||_{C^0}$ is bounded by a constant depending on $M,G_0,\Theta$, and $||f||_{C^0}$ only.

  Let us fix a finite covering $\{U_\alp\}$ by open charts inducing the $Spin(9)$-affine structure on $M$, i.e. there are diffeomorphisms
$$f_\alp\colon U_\alp\tilde\to V_\alp$$
with open subsets $V_\alp\subset \RR^{16}$ such that $f_\beta\circ f_\alp^{-1}$ is a composition of a translation and of a transformation from the group $Spin(9)$. Furthermore all $f_\alp$
can be chosen so that the corresponding to $G_{00}$ matrix on $V_\alp$ is equal to $I_2\in \ch_2(\OO)$. We can choose  compact subsets
$ K_\alp \subset U_\alp$
such that $\cup_\alp K_\alp=M$. Clearly there exists $C_1>0$ such that
\begin{eqnarray}\label{E:C-1-p}
G_0(x)>C_1 \cdot I_2 \mbox{ for any } \alp \mbox{ and any } x\in K_\alp.
\end{eqnarray}
Here and below we tacitly identify $U_\alp$ with its image $V_\alp$ under $f_\alp$  and sections of $\uch$ with corresponding matrix valued functions.


Let $z_0$ be a point of minimum of $\phi$. Then one has $Hess_\OO(\phi(z_0))\geq 0$. Hence $G_0(z_0)+Hess_\OO\phi(z_0)\geq 0$. But $\det(G_0+Hess_\OO(\phi))>0$ on $M$. Hence by Lemma \ref{L:ellipticity}
\begin{eqnarray*}\label{E:lower-estimate-Laplace}
Hess_\OO(\phi)> -G_0.
\end{eqnarray*}
Hence it remains to obtain an upper bound on $\Delta \phi$.

Let us consider the function
$$T:=\frac{4}{3}[tr(G^{-1}_{00}\cdot (G_0+Hess_\OO(\phi)))]^{3/4}-\phi.$$
This function is globally defined (this can be proven using a variation argument in the same way as in paragraph \ref{item:DEF-Delta} of this section for $\Delta$). In order to
prove the theorem, it suffices to show that $T$ is bounded
from above (due to the boundedness of $||\phi||_{C^0}$). Let $z\in M$ be a point of maximum of the function $T$. We clearly have
$$\Delta'(T)(z)\leq 0.$$
Near $z$ let us represent $G_0=Hess_\OO (u_0)$. Set $u:=u_0+\phi$ near the point $z$; this function is strictly oPSH. Denote $U:=Hess_\OO(u)=G_0+Hess_\OO(\phi)$.


Let us choose $\alp$ such that  $z\in K_\alp$.  Then we have
\begin{eqnarray}\label{E:ineq-C1}
tr(U^{-1}G_0)|_z\overset{(\ref{E:C-1-p})}{\geq} C_1 \, tr(U^{-1}(z)).
\end{eqnarray}

By Proposition  \ref{P:diagonalization} there exists a transformation $g\in Spin(9)$ which diagonalizes $U(z)$ (or rather the matrix corresponding to
$U(z)$ under $f_\alp$).


Then, near $z$ one has
$T=\frac{4}{3}(\Delta u)^{3/4}-\phi.$
Consequently,
\begin{eqnarray*}
0\geq \Delta'(\frac{4}{3}(\Delta u(z))^{3/4})-\Delta'\phi=\Delta'(\frac{4}{3}(\Delta u(z))^{3/4})-tr(U^{-1}(U-G_0))|_z=\\
\Delta'(\frac{4}{3}(\Delta u(z))^{3/4})-2+tr(U^{-1}G_0)\big|_z\overset{(\ref{E:ineq-C1})}{\geq} \\
\Delta'(\frac{4}{3}(\Delta u(z))^{3/4})-2+C_1\sum_{i=1}^2\frac{1}{u_{\bar ii}(z)}.
\end{eqnarray*}
 Applying now Proposition \ref{P:Laplace-twice} we get for $C_2>0$
\begin{eqnarray*}
0\geq (\Delta u(z))^{-1/4}(\Delta \log F)-2+C_1\sum_i\frac{1}{u_{\bar ii}(z)}\geq -C_2(\Delta u(z))^{-1/4}-2+C_1\sum_i\frac{1}{u_{\bar ii}(z)},
\end{eqnarray*}
where $F=e^f\det G_0$ as above. This is equivalent to
\begin{eqnarray}\label{E:ooo1}
(u_{\bar 11}(z)+u_{\bar 22}(z))^{-1/4}\geq C_3(\frac{1}{u_{\bar 11}(z)}+\frac{1}{u_{\bar 22}(z)})-C_4,
\end{eqnarray}
where $C_3,C_4>0$. Recall also that
\begin{eqnarray}\label{E:ooo2}
u_{\bar 11}(z)u_{\bar 22}(z)=F(z).
\end{eqnarray}
Without loss of generality we may assume that
\begin{eqnarray}\label{E:bolshe}
u_{\bar 11}(z)\geq u_{\bar 22}(z).
\end{eqnarray}
Then (\ref{E:ooo1}) implies that
\begin{eqnarray}\label{E:ooo4}
[u_{\bar 11}(z)]^{-1/4}\geq \frac{C_3}{F}u_{\bar 11}(z)-C_4.
\end{eqnarray}
If $u_{\bar 11}(z)\leq \frac{C_4+1}{C_3}F$ then the theorem is proven. Let us assume otherwise. Then (\ref{E:ooo4}) implies
$$[u_{\bar 11}(z)]^{-1/4}\geq 1.$$
Hence $u_{\bar 11}(z)\leq 1$ and the theorem is proven again. \qed

\end{paragraphlist}

\def\lam{\lambda}
\def\Lam{\Lambda}
\section{A weak Harnack inequality.}\label{S:Harnack}
Let us recall now the weak Harnack inequality (see e.g. \cite{gilbarg-trudinger}, Theorem 8.18, or
\cite{han-lin}, Theorem 4.15). Below we normalize everywhere the Lebesgue measure
on $\RR^N$ by $vol(B_1) = 1$ where $B_1$ is Euclidean unit ball.
\begin{theorem}[weak Harnack inequality]\label{harnack}
Let $B_R\subset \RR^N$ be an open Euclidean ball of radius $R$. Let
$(a_{ij})_{i,j=1}^N\in L^\infty(B_R)\cap C^1(B_R)$ be a symmetric matrix satisfying
 uniform elliptic estimate
$$\lam||\xi||^2\leq \sum_{i,j}a_{ij}(x)\xi_i\xi_j\leq
\Lam||\xi||^2,\mbox{ for all }\xi\in \RR^N,$$ with $\lam,\Lam>0$. Let
$v\in C^2(B_R)$ be a function satisfying
\begin{eqnarray*}
v\geq 0,\\
\sum_{i,j}D_j(a_{ij}D_iv)\leq \psi,
\end{eqnarray*}
where $\psi\in L^\infty(B_R)$. Then, for any $0<\theta<\tau<1$ we
have
$$\inf_{B_{\theta R}}v+R||\psi||_{L^\infty(B_R)}\geq
C\left(R^{-N}\int_{B_{\tau R}}v\right),$$ where the constant $C$
depends only on $\lam,\Lam,\theta,\tau,N$.
\end{theorem}
\begin{remark}\label{R:harnack}
 Theorem \ref{harnack} also holds in the following weaker form. We
will take $R=4r, \theta=1/4,\tau=1/2$. Then,
\begin{eqnarray}
r^{-N}\int_{B_r}v\leq C'(\inf_{B_r} v+r),
\end{eqnarray}
where the constant $C'$ depends on
$\lam,\Lam,||\psi||_{L^\infty(B_{4r})},N$ only.
\end{remark}

In the next Lemma \ref{L:real-operator} and Proposition \ref{P:divergence-form} we show that the operator $\cd h:=\det U\cdot tr(U^{-1}\cdot Hess_\OO (h))$, where $U=Hess_\OO(u)$,
can be written in the divergence form, and hence later we will be able to apply to it Theorem \ref{harnack}.

\begin{lemma}\label{L:real-operator}
Let $u,h\colon \OO^2\to \RR$ be smooth functions.
Denote $U=Hess_\OO(u), H=Hess_\OO(h)\in \ch_2(\OO)$. Then,
$$\det U\cdot tr(U^{-1}H)=u_{\bar 22}h_{\bar 11}+u_{\bar 11}h_{\bar 22}-2\sum_{r,s}Re\left(\bar e_r u_{\bar 12} e_s\right)h_{x_1^rx_2^s}.$$
\end{lemma}
{\bf Proof.} By Section \ref{S:integration}, paragraph \ref{item:integr-1}, we have
$$(\det U)\cdot U^{-1}=adj(U)=\left[\begin{array}{cc}
                             u_{\bar 22}&-u_{\bar 12}\\
                             -u_{\bar 21}&u_{\bar 11}
                              \end{array}\right].$$

Hence the left hand side of the above equality is equal to
\begin{eqnarray*}
tr\left(\left[\begin{array}{cc}
                             u_{\bar 22}&-u_{\bar 12}\\
                             -u_{\bar 21}&u_{\bar 11}
                              \end{array}\right]\cdot \left[\begin{array}{cc}
                             h_{\bar 11}&h_{\bar 12}\\
                             h_{\bar 21}&h_{\bar 22}
                              \end{array}\right]\right)=
                             u_{\bar 22}h_{\bar 11}+u_{\bar 11}h_{\bar 22}-2Re(u_{\bar 12}h_{\bar 21})=\\
                             u_{\bar 22}h_{\bar 11}+u_{\bar 11}h_{\bar 22}-2Re(u_{\bar 12}(\sum_{rs} e_s h_{x_1^rx_2^s}\bar e_r))=\\
                              u_{\bar 22}h_{\bar 11}+u_{\bar 11}h_{\bar 22}-2\sum_{rs} Re(u_{\bar 12} e_s \bar e_r)h_{x_1^rx_2^s}.
\end{eqnarray*}
The lemma follows since $Re(ab)=Re(ba)$. \qed

\hfill

Let $\cd (h)=\sum_{m,n}a_{mn}\pt_m\pt_n h$ be a second order differential operator in the real space $\RR^N$ with $a_{mn}=a_{nm}$. We say that $\cd$ can be written in the divergence form if
$$\cd (h)=\sum_{m,n}\pt_m(a_{mn}\pt_n h),$$
or equivalently
\begin{eqnarray}\label{E:divergence}
\sum_{m}\pt_m a_{mn}=0 \mbox{ for any } n.
\end{eqnarray}

\begin{proposition}\label{P:divergence-form}
Let $U=Hess_\OO(u)$ where $u$ is strictly oPSH. Let $$\cd (h)=\det U\cdot tr(U^{-1}\cdot Hess_\OO (h)).$$
Then $\cd$ can be written in the divergence form.
\end{proposition}
{\bf Proof.} By Lemma \ref{L:real-operator}
\begin{eqnarray*}
\cd=u_{\bar 22}\sum_{r}\pt_{x_1^r}\pt_{x_1^r}+u_{\bar 11}\sum_s \pt_{x_2^s}\pt_{x_2^s}-2\sum_{r,s} Re(\bar e_ru_{\bar 12}e_s)\pt_{x_1^r}\pt_{x_2^s}=\\
u_{\bar 22}\sum_{r}\pt_{x_1^r}\pt_{x_1^r}+u_{\bar 11}\sum_s \pt_{x_2^s}\pt_{x_2^s}-\sum_{r,s} Re(\bar e_ru_{\bar 12}e_s)\pt_{x_1^r}\pt_{x_2^s}-\sum_{r,s} Re(\bar e_ru_{\bar 12}e_s)\pt_{x_2^s}\pt_{x_1^r}.
\end{eqnarray*}

First assume in (\ref{E:divergence}) that $x_n=x_1^r$. Then the left hand side of (\ref{E:divergence}) is equal to
\begin{eqnarray*}
\pt_{x_1^r}u_{\bar 22}-\sum_s \pt_{x_2^s}\left(Re(\bar e_r u_{\bar 12}e_s)\right)=\\
\Delta_2u_{x_1^r}-Re(\bar e_r(\sum_s u_{\bar 12,x_2^s}e_s))=\\
\Delta_2u_{x_1^r}-Re\left[\bar e_r\left((u_{\bar 1}\overset{\leftarrow}{\pt_{ 2}})\overset{\leftarrow}{\pt_{\bar 2}}\right)\right]\overset{(\ref{E:derivatives})}{=}\\
\Delta_2u_{x_1^r}-Re(\bar e_r\Delta_2 u_{\bar 1})=\\
\Delta_2u_{x_1^r}-\Delta_2(Re(\bar e_ru_{\bar 1}))=0.
\end{eqnarray*}

The case $x_n=x_2^s$ is considered in the same manner. \qed

\section{An auxiliary estimate for oPSH functions.}\label{S:estimate}
Let $L\subset \OO^2$ be an affine octonionic line. Let $v_0,v_1,\dots,v_7$ be an orthonormal basis of the linear subspace parallel to $L$.
Consider the linear differential operator
$$\Delta_L:=\sum_{p=0}^7 \frac{\pt^2}{\pt v_p^2}.$$
Clearly $\Delta_L$ is independent of the choice of a basis.
\begin{lemma}\label{L:p-kvadrat}
Let $u\in C^4$ be strictly oPSH function on an open subset $\co\subset \OO^2$ such that
$$MA_\OO(u)=F>0,$$
where $F\in C^\infty(M,\RR)$.
Then for any affine octonionic line $L$  we have
\begin{eqnarray}
tr(U^{-1}\cdot\Delta_L U)\geq \Delta_L (\log F).
\end{eqnarray}
\end{lemma}
\def\xp1{x_p^1x_p^1}
{\bf Proof.} Let $v_0,\dots,v_7$ be an orthonormal basis of the linear subspaces parallel to $L$. Then
$\Delta_L=\sum_{p=0}^7\frac{\pt^2}{(\pt v_p)^2}$. It is enough to show
that
\begin{eqnarray}\label{4}
tr\left(U^{-1}\cdot\frac{\pt^2U}{\pt v_p^2} \right)\geq \frac{\pt^2(\log F)}{\pt v_p^2}
\end{eqnarray}
for any $p=0,\dots, 7$. By Claim \ref{Cl:3}(2) we have
\begin{eqnarray*}\label{5}
tr(U^{-1}\cdot \frac{\pt^2U}{\pt v_p^2})=
\frac{\pt^2(\log F)}{\pt v_p^2}+tr(U^{-1}\frac{\pt U}{\pt v_p}U^{-1}\frac{\pt U}{\pt v_p}).
\end{eqnarray*}
In order to prove the lemma, it suffices to show that
$tr(U^{-1}\frac{\pt U}{\pt v_p}U^{-1}\frac{\pt U}{\pt v_p})\geq 0$. More generally let us
show that if $A,B\in \ch_2(\OO)$ and $A>0$ then
$$tr(A^{-1}BA^{-1}B)\geq 0.$$
All elements of $A,B$ belong to the associative subfield generated by $a_{\bar 12}, b_{\bar 12}$. Any such subfield is contain in a subfield isomorphic to quaternions $\HH$. Note that
the matrices $A$ and $B$ can be diagonalized simultaneously over $\HH$. More precisely there exists an invertible quaternionic matrix $X$ and a real diagonal matrix $D$ such that $A=X^*X$ and $B=X^*DX$.

Then
\begin{eqnarray*}
tr(A^{-1}BA^{-1}B)=tr(X^{-1}D^2 X)=tr(D^2)\geq 0.
\end{eqnarray*}
Lemma is proved. \qed

\section{A proposition from linear algebra.}
\begin{paragraphlist}
\item Let $A\in \ch_2(\OO)$. Its characteristic polynomial $p_A(t)=\det(A-tI_2)$ has two real roots which are called spectrum of $A$. Indeed, the subfield generated by $a_{\bar 12}$ is associative and commutative, and the result
follows from the case of complex Hermitian matrices.

The action of $Spin(9)$ on $\ch_2(\OO)$ preserves $\det$ and hence the characteristic polynomial. Hence spectrum of a matrix does not change under the $Spin(9)$-action.

\item Let $\zeta\in \OO^2$ be a vector. Define the octonionic Hermitian matrix
$$\zeta\otimes \zeta^*:=(\zeta_i\bar\zeta_j).$$

\begin{claim}\label{Cl:pos-def-rank-one}
The matrix $\zeta\otimes \zeta^*$ is non-negative definite.
\end{claim}
{\bf Proof.}  It suffices to show that for any $\eps>0$ the matrix
$$\zeta\otimes\zeta^*+\eps I_2=\left[\begin{array}{cc}
                                     |\zeta_1|^2+\eps&\zeta_1\bar\zeta_2\\
                                     \zeta_2\bar\zeta_1&|\zeta_2|^2+\eps
                                     \end{array}\right]$$
is positive definite. This follows from Proposiiton \ref{P:sylvester} (the Sylvester criterion). Indeed the diagonal entries are positive, and the determinant is equal to
$$(|\zeta_1|^2+\eps)(|\zeta_2|^2+\eps)-|\zeta_1\bar\zeta_2|^2=\eps^2+\eps(|\zeta_1|^2+|\zeta_2|^2)>0.$$
\qed

\item We will need the following linear algebraic proposition which is
completely analogous to the real, complex, and quaternionic cases (see
\cite{gilbarg-trudinger}, Lemma 17.13, for the real case; for
the complex case \cite{siu}, p. 103, or \cite{blocki-lecture-notes},
Lemma 5.17; and for the quaternionic case see  \cite{alesker-calabi-yau}, Lemma 4.9.

\begin{proposition}\label{L:hermit}
Fix $0<\lam<\Lam<\infty$. Then there exist a natural number $N$,
unit vectors $\zeta_1,\dots,\zeta_N\in\OO^2$ such that at least one of the two coordinates of each $\zeta_i$ is real, and
$0<\lam_*<\Lam_*<\infty$ depending on $\lam,\Lam$ only such that
every  $A\in \ch_2(\OO)$ whose spectrum lies in
$[\lam,\Lam]$ can be written
$$\sum_{k=1}^N\beta_k\zeta_k\otimes \zeta_k^* \mbox{ with }
\beta_k\in [\lam_*,\Lam_*].$$ Moreover, the vectors $\zeta_1,\dots,\zeta_N$
can be chosen to contain the vectors $(1,0),(0,1)\in\OO^2$.
\end{proposition}
{\bf Proof.} The proof is standard in most part. For convenience of the reader we outline the argument. Let us denote by $K$ the set of octonionic
Hermitian matrices whose spectrum lies in $[\lam,\Lam]$. This is a
compact subset of the interior of the cone of positive definite
octonionic Hermitian matrices; we denote the latter open cone by $\cc$. Then
there exists a convex compact polytope $P\subset \cc$ which contains
a neighborhood of $K$. Using a diagonalization (over a subfield of $\OO$ isomorphic to $\CC$), every matrix $X\in \cc$, in particular any vertex of $P$, can be written
in the form
\begin{eqnarray}\label{E:2-1}
X=\sum_{i=1}^2\alp_i(X) \left(\xi_i(X)\otimes \xi_i^*(X)\right),
\end{eqnarray}
where $\alp_i(X)>0$, and $\xi_i(X)\in\OO^2$ are unit vectors with one of the two coordinates being real. Let us
fix such a decomposition for any vertex of $P$. Let us define a new
finite subset $S_1\subset \cc$ as follows
\begin{eqnarray*}
S_1:=\left\{\left(\sum_j\alp_j(X)\right)\cdot
\left(\xi_i(X)\otimes\xi_i^*(X)\right)|\, X\in Vert(P),
\alp_j(X),\xi_i(X) \mbox{ satisfy } (\ref{E:2-1})\right\}.
\end{eqnarray*}
Let us add to $S_1$ matrices of the form $e_i\otimes e_i^*$,
$i=1,2$, where $e_1=(1,0),e_2=(0,1)$. This set we denote by $S$. It is clear that
$P\subset conv(S)$. Hence $conv(S)$ contains a neighborhood of $K$.
Now the lemma follows from the following general result (where one takes $Q=conv(S)$).
\begin{lemma}\label{L:convex-lemma}
 Let $K$ be a
compact subset of a finite dimensional real vector space $V$ which is contained in the interior of a
compact convex polytope $Q$. Then there exists $\eps>0$ such that
any point $x\in K$ can be written as a convex combination of
vertices of $Q$:
$$x=\sum_v \beta(v)v \mbox{ with } \beta(v)>\eps,$$
where the sum runs over all vertices of $Q$, $\sum_v\beta(v)=1$.
\end{lemma}
It remains to prove Lemma \ref{L:convex-lemma}.

{\bf Proof.} First let us show that any point from the interior of $Q$ can be presented as a convex combination of vertices with strictly positive coefficients. 

Let us denote by $Vert$ the set of vertices of $Q$. Applying a translation we may and will assume that $0\in Vert$.
Consider the vector subspace $V\subset \RR^{Vert}$ defined by 
$$\cx:=\{\{x_v\}_{v\in Vert}|\,\, \sum_v x_v =0.\}.$$
Consider the linear map $T\colon \cx\to V$ given by
\begin{eqnarray}\label{E:convex-map}
T(\{x_v\})= \sum_v x_v\cdot v.
\end{eqnarray}
This map is onto since $Q$ has non-empty interior and $0\in Vert$. Consider the simplex 
$$\Delta:=\{\{x_v\}|\,\, \sum_v x_v=1,\, x_v>0\}.$$
Then $T(\Delta)$ is an open convex subset of the interior of $Q$. Clearly its closure is
$$\overline{T(\Delta)}=Q.$$
But for any convex open set $A$ one has $A=int(\bar A)$. Hence 
$$T(\Delta)=int(Q).$$
In other words every point from the interior of $Q$ can be presented as a convex combination of vertices of $Q$ with strictly positive coefficients. 

To finish proving the lemma we will now assume that $0\in int (Q)$. By what we just proven there exists a presentation
$$0=\sum_v\alp_v\cdot v, \mbox{ with } \sum_{v\in Vert}\alp_v=1,\, \alp_v>0.$$
There exists $0<\lam <1$ such that $K\subset \lam Q$. Let $x\in K$ be any point. Then there exists a presentation
$$x=\lam \sum_v\gamma_v v\mbox{ with } \sum_v\gamma_v=1,\,\gamma_v\geq 0.$$
Then we have
\begin{eqnarray*}
x=\lam\sum_v\gamma_v \cdot v+(1-\lam)\sum_v \alp_v\cdot v=\\
\sum_v (\lam \gamma_v+(1-\lam)\alp_v)\cdot v.
\end{eqnarray*}
This is a convex combination of vertices and 
$$\inf_v (\lam \gamma_v+(1-\lam)\alp_v)\geq (1-\lam)\inf_v\alp_v.$$
Taking $\eps:=(1-\lam)\inf_v \alp_v>0$ we get the lemma.
\qed

\hfill

Let us state one more estimate on matrices to be used later.
\begin{lemma}\label{L:linear-algebra}
Let $A,B\in\ch_2(\OO)$ and $A,B>0$. Then
$$tr(A^{-1}(A-B))\leq \frac{2}{\sqrt{\det A}}(\sqrt{\det A}-\sqrt{\det B}).$$
\end{lemma}
{\bf Proof.} The elements of $A,B$ generate an associative subfield contained in a subfield isomorphic to $\HH$. Over $\HH$ the statement was proved in Lemma 4.8 in \cite{alesker-calabi-yau} (using simultaneous diagonalization of $A$ and $B$
and directly applying the arithmetic-geometric mean inequality). \qed

\end{paragraphlist}

\def\dz{\Delta_\zeta}
\section{$C^{2,\alp}$-estimate.}\label{S:local-estimate}
\begin{paragraphlist}
\item In this section $M$ is a compact connected $Spin(9)$-affine manifold with a flat $C^\infty$-smooth positive definite section $G_{00}$ of $\uch$. We denote as in Section \ref{S:second-order}, paragraph \ref{item:DEF-Delta},
the operator
$$\Delta f=tr(G_{00}^{-1}Hess_\OO(f)).$$
$\Delta$ is elliptic, and one can choose local coordinates consistent with the $Spin(9)$-structure such that $\Delta$ is proportional to the usual flat Laplacian on $\OO^2=\RR^{16}$.

Furthermore we assume as previously that $f\in C^\infty(M,\RR)$, and $G_0$ is a $C^\infty$-smooth positive definite section of $\uch$ satisfying condition (\ref{E:condition}). Let $\phi\in C^4$ satisfy the Monge-Amp\`ere equation
$$\det (G_0+Hess_\OO(\phi))=e^f\det G_0.$$

\item The main result of this section is
\begin{theorem}\label{thm-main}
 There exist $\alp\in (0,1)$ and $C>0$, both depending on
$M$,$G_0$,
$||f||_{C^2(M)}$,$||\phi||_{C^0(M)}$,$||\Delta\phi||_{C^0(M)}$ only,
such that
$$||\phi||_{C^{2,\alp}(M)}\leq C.$$
\end{theorem}
This theorem immediately follows from the below local estimate.

\begin{theorem}\label{thm-local}
Let $u\in C^4$ be a strictly oPSH function in an open subset
$\co\subset \OO^2$ satisfying
\begin{eqnarray}\label{2}
\det(Hess_\OO(u))=F>0
\end{eqnarray}
with $F\in C^2(\co)$. Let $\co'\subset \co$ be a relatively compact open
subset. Then there exist $\alp\in(0,1)$ depending on
$||u||_{C^0(\co)},||\Delta u||_{C^0(\co)},$ $ ||\log F||_{C^2(\co)}$
only (here $\Delta$ is the flat Laplacian on $\OO^2=\RR^{16}$) and a constant $C>0$ depending in addition on
$dist(\co',\pt\co)$ such that
$$||u||_{C^{2,\alp}(\co')}\leq C.$$
\end{theorem}

\item First let us introduce a notation. Let $\zeta\in\OO^2$ be a unit vector such that at least one of its coordinates is real.
Let $\dz$ be the Laplacian on the octonionic line (or its translates) $\{\zeta\cdot q|\, q\in \OO\}$ spanned by $\zeta$ (see Section \ref{S:estimate}).
The following lemma was proved in \cite{alesker-octon}, see equality (28) there.
\begin{lemma}\label{L:laplace-line}
For any smooth real values function $f$ one has
$$\dz f(x)=Re(\zeta^*Hess_\OO(f)\zeta).$$
\end{lemma}

\item {\bf Proof of Theorem \ref{thm-local}.}
\newline
\underline{Step 1.} Lemma \ref{L:linear-algebra} implies that for any point $x,y\in \co$
\begin{eqnarray}\label{E:II}
F(y)tr\left(U^{-1}(y)(U(y)-U(x))\right)\leq \\\label{E:II-1}
2F^{1/2}(y)(F^{1/2}(y)-F^{1/2}(x))\leq
C_1||x-y||,
\end{eqnarray}
where $C_1$ depends on $||\log F||_{C^1(\co)}$ only.

\hfill

\underline{Step 2.} Let $\cd$ be the operator from Proposition \ref{P:divergence-form}, i.e.
$$\cd (h)=\det U\cdot tr(U^{-1}\cdot Hess_\OO (h))=F\cdot tr(U^{-1}\cdot Hess_\OO (h)).$$
 We claim that $\cd$ is uniformly elliptic in $\co$ with constants
$\lam,\Lam$ depending on $||\log F||_{C^0(\co)},
||\Delta u||_{C^0(\co)}$ only.

Indeed, fix a point $z\in \co$. By Propositions \ref{P:ohs1} and \ref{P:diagonalization} we can make change of variables $\xi'=T\xi$ where $T\in Spin(9)$ and assume that $U(z)$ is diagonal. Then
$$\cd h(z)=u_{\bar 11}(z)h_{\bar 22}(z)+u_{\bar 22}(z)h_{\bar 11}(z).$$
The uniform ellipticity follows from this, the estimate $0<u_{\bar ii}(z)\leq ||\Delta u||_{C^0(\co)}$, and the equality $u_{\bar 11}(z)u_{\bar 22}(z)=F(z)$.

\hfill

\underline{Step 3.} Fix a point $z_0\in\co'$ and a ball $B_{4r}=B(z_0,4r)\subset\co$. We are going to apply the weak
Harnack inequality to the operator $\cd$ acting on the function
$$v:=\sup_{B_{4r}}\dz u-\dz u.$$

The function $v$ satisfies by Lemma \ref{L:p-kvadrat}
$$\cd v\leq -\dz (\log F)$$
where we have used that $\Delta_\zeta(u_{\bar i j})=(\Delta_\zeta
u)_{\bar i j}$.

Applying the weak Harnack inequality in the form given in
Remark \ref{R:harnack} we get:
\begin{eqnarray}\label{harnack-main}
r^{-16}\int_{B_r}(\sup_{B_{4r}}\dz u-\dz u)\leq C(\sup_{B_{4r}}\dz
u-\sup_{B_r}\dz u +r)
\end{eqnarray}
where $C$ depends on $||\log F||_{C^2(\co)}$ and $||\Delta u||_{C^0(\co)}$ only. We apply this inequality shortly.

\hfill

\underline{Step 4.} The eigenvalues of
$U=Hess_\OO(u)$ belong to $[\lam_1,\Lam_1]$ with
$0<\lam_1<\Lam_1<\infty$ depending on $||\Delta u||_{C^0(\co)},
||\log F||_{C^0(\co)}$ only. Hence the eigenvalues of $FU^{-1}$ are in
$[\lam_2,\Lam_2]$ with $0<\lam_2<\Lam_2<\infty$ are under control.\footnote{In fact one can take $\lam_2=\lam_1$ and $\Lam_2=\Lam_1$.}
By Proposition
\ref{L:hermit} there exist a positive integer $N$, unit vectors
$\zeta_1,\dots,\zeta_N\in \OO^2$ such that each $\zeta_k$ has at least one real coordinate, and $0<\lam_*<\Lam_*<\infty$ (under control) such
that for any $y\in\co$
$$F(y)U^{-1}(y)=\sum_{i=1}^N\beta_k(y)\zeta_k\otimes \zeta_k^*
\mbox{ with }\beta_k(y)\in [\lam_*,\Lam_*].$$

Observe also that for a unit vector $\zeta\in \OO^2$ such that at least one of its coordinates is real one has
$$tr((\zeta\otimes\zeta^*)(Hess_\OO(u)))=tr(\zeta^*(Hess_\OO(u))\zeta)\overset{\mbox{Lemma }\ref{L:laplace-line}}{=}\Delta_\zeta u.$$
This and (\ref{E:II})-(\ref{E:II-1})
imply
\begin{eqnarray}\label{E:III}
\sum_{k=1}^N\beta_k(y)(\Delta_{\zeta_k}u(y)-\Delta_{\zeta_k}u(x))\leq
C_1||x-y|| \mbox{ for } x,y\in \co.
\end{eqnarray}

\underline{Step 5.}
Let $\co''$ be a compact neighborhood of $\co'$ in $\co$. Let us fix now $z_0\in \co''$ and a ball $B_{4r}=B(z_0,4r)\subset \co$ centered at $z_0$.
Let us denote
$$M_{k,r}:=\sup_{B_r}\Delta_{\zeta_k}u,\,\,\,
m_{k,r}:=\inf_{B_r}\Delta_{\zeta_k}u,$$
$$\eta(r):=\sum_{k=1}^N(M_{k,r}-m_{k,r}).$$ We will show that for
some $\alp\in (0,1), r_0>0,C>0$ under control
\begin{eqnarray}\label{E:eta-estim}
\eta(r)\leq Cr^\alp \mbox{ for } 0<r<r_0.
\end{eqnarray}

Since $\zeta_1,\dots,\zeta_N$ can be chosen to contain the vectors $(1,0),(0,1)\in \OO^2$,
this will imply an estimate on
$||\Delta u||_{C^\alp(\co'')}$. Then the Schauder estimates
(\cite{gilbarg-trudinger}, Theorem 4.6) will imply an estimate on
$||u||_{C^{2,\alp}(\co')}$.

The condition $\eta(r)\leq Cr^\alp$ is equivalent to
$$\eta(r)\leq \delta \eta(4r)+r,\, 0<r<r_1$$
where $\delta\in (0,1), r_1$ are under control
(\cite{gilbarg-trudinger}, Lemma 8.23). Summing up
(\ref{harnack-main}) over $\zeta=\zeta_l$ with $l\ne k$ for any
fixed $k$ we get
\begin{eqnarray}\label{E:IV}
r^{-16}\int_{B_r}\sum_{l\ne k}(M_{l,4r}-\Delta_{\zeta_l}u)\leq
C_3(\eta(4r)-\eta(r)+r).
\end{eqnarray}

By (\ref{E:III}) we have for any $x\in B_{4r},y\in B_r$
\begin{eqnarray*}
\beta_k(y)(\Delta_{\zeta_k}u(y)-\Delta_{\zeta_k}u(x))\leq
C_1||x-y||+\sum_{l\ne
k}\beta_l(y)(\Delta_{\zeta_l}u(x)-\Delta_{\zeta_l}u(y))\leq\\
C_4r+\Lam_*\sum_{l\ne k}(M_{l,4r}-\Delta_{\zeta_l}u(y)).
\end{eqnarray*}

Optimizing in $x$ we get $$\Delta_{\zeta_k}u(y)-m_{k,4r}\leq
\frac{1}{\lam_*}\left(C_4r+\Lam_*\sum_{l\ne
k}(M_{l,4r}-\Delta_{\zeta_l}u(y))\right).$$ Integrating the last
inequality over $y\in B_r$ and using (\ref{E:IV}) we obtain
\begin{eqnarray}\label{E:V}
r^{-16}\int_{B_r}(\Delta_{\zeta_k}u(y)-m_{k,4r})\leq
C_5(\eta(4r)-\eta(r)+r).
\end{eqnarray}
Let us estimate the left hand side of (\ref{E:V}) from below. Since
we have normalized the Lebesgue measure on $\OO^2$ so that
$vol(B_1)=1$, we have
\begin{eqnarray*}
r^{-16}\int_{B_r}(\Delta_{\zeta_k}u(y)-m_{k,4r})\overset{(\ref{harnack-main})}{\geq}
-m_{k,4r}+M_{k,4r}-C(M_{k,4r}-M_{k,r}+r)=\\
C(M_{k,r}-m_{k,r})-(C-1)(M_{k,4r}-m_{k,4r})+C(m_{k,r}-m_{k,4r})-Cr\geq\\
C(M_{k,r}-m_{k,r})-(C-1)(M_{k,4r}-m_{k,4r})-Cr.
\end{eqnarray*}
Substituting this back into (\ref{E:V}) we obtain
\begin{eqnarray*}
C(M_{k,r}-m_{k,r})-(C-1)(M_{k,4r}-m_{k,4r})-Cr\leq
C_5(\eta(4r)-\eta(r)+r).
\end{eqnarray*}
Summing this up over $k$ we get
\begin{eqnarray*}
C\eta(r)-(C-1)\eta(4r)\leq C_6(\eta(4r)-\eta(r)+r).
\end{eqnarray*}
Hence, $$\eta(r)\leq \frac{C+C_6-1}{C+C_6}\eta(4r)+r.$$ Theorem
\ref{thm-local} is proved. \qed

\end{paragraphlist}

\section{Proof of main Theorem \ref{T:main-result}.}\label{S:main-thm}
\begin{paragraphlist}
\item For any integer $k\geq 1$ and $\beta\in (0,1)$ let us define
\begin{eqnarray*}
\cu^{k,\beta}:=\{\phi\in C^{k,\beta}(M)|\, \det G_0+Hess_\OO\phi>0
\mbox{ and } \int_M \phi\cdot
\det G_0\cdot\Theta=0 \},\\
\cv^{k,\beta}:=\{\chi\in C^{k,\beta}(M)|\, \chi>0 \mbox{ and }
\int_M(\chi-1)\cdot\det G_0\cdot\Theta=0\},
\end{eqnarray*}
and set
$$\cm(\phi):=\frac{\det(G_0+Hess_\OO(\phi))}{\det G_0}.$$
The continuity of $$\cm\colon \cu^{k+2,\beta}\to \cv^{k,\beta}$$ is obvious. The only thing to check
is that for any $\phi\in C^{k+2,\beta}(M)$ one has
\begin{eqnarray*}\label{E:pp1}
\int_M (\cm(\phi)-1)\det G_0\cdot \Theta=0.
\end{eqnarray*}
This easily follows from Lemma \ref{L:const-A}.

\item Note that $\cu^{k+2,\beta}$ and $\cv^{k,\beta}$ are open
subsets in Banach spaces with the induced H\"older norms. The following proposition holds.
\begin{proposition}\label{P:open-image}
Let $M^{16}$ be a compact connected $SL_2(\OO)$-affine manifold with a parallel $3/4$-density $\Theta$.
Let $k\geq 2$ be an integer, and let $\beta\in (0,1)$. Then the map
$\cm\colon \cu^{k+2,\beta}\to \cv^{k,\beta}$ is locally a
diffeomorphism of Banach spaces, and in particular its image is an
open subset.
\end{proposition}
{\bf Proof.} By the inverse function theorem for Banach spaces, it
suffices to show that the differential of $\cm$ at any $\phi\in
\cu^{k+2,\beta}$ is an isomorphism of tangent spaces. The tangent
space to $\cu^{k+2,\beta}$ at any point is
$$\tilde C^{k+2,\beta}(M):=\{\phi\in C^{k+2,\beta}(M)|\, \int_M
\phi\cdot \det G_0\cdot \Theta=0\}.$$ The tangent space to
$\cv^{k,\beta}$ at $\cm(\phi)$ is
$$\tilde C^{k,\beta}(M):=\{\chi\in C^{k,\beta}(M)|\, \int_M \chi
\cdot \det G_0\cdot\Theta =0\}.$$

The differential of $\cm$ at $\phi$ is equal to
\begin{eqnarray*}
D\cm_\phi(\psi)=2\frac{\det(G_0+Hess_\OO\phi,Hess_\OO\psi)}{\det G_0}.
\end{eqnarray*}
Being defined by this formula we consider $D\cm_\phi$ as a map
$C^{k+2,\beta}(M)\to C^{k,\beta}(M)$ (without restricting to $\tilde
C$). Then obviously $D\cm_\phi$ is a linear differential second order
elliptic operator. It has no free term (i.e.
$D\cm_\phi(1)=0$), and its coefficients belong to $C^{k,\beta}$.

By the strong maximum principle (see e.g. \cite{gilbarg-trudinger},
Theorem 8.19) and since $M$ is connected, the kernel of $D\cm_\phi$
consists only of constant functions; in particular it is one dimensional. The
image of $D\cm_\phi$ is a closed subspace of $C^{k,\beta}(M)$ by
\cite{joyce-book}, Theorem 1.5.4, which is a version of the Schauder
estimates.

We claim that the index of $D\cm_\phi$ is zero. That would imply that $$codim\, Im(D\cm_\phi)=1.$$
Since $Im(D\cm_\phi)\subset
\tilde C^{k,\beta}(M)$, and since $\dim C^{k,\beta}(M)/\tilde
C^{k,\beta}(M)=1$, it follows that $Im(D\cm_\phi)=\tilde
C^{k,\beta}(M)$.

It remains to prove that the index of $D\cm_\phi$ is zero. For linear elliptic operators with infinitely smooth coefficients on real valued functions this is always true
since the adjoint operator has the same symbol as the original operator.
The operators with H\"older coefficients require more care. By \cite{joyce-book}, Theorem 1.5.4, the image of a linear elliptic operator of second order with $C^{k,\alp}$-coefficients, $k\geq 2$,
$$P\colon C^{k+2,\alp}(M)\to C^{k,\alp}(M)$$
is equal to the orthogonal complement of the kernel of its adjoint $P^*$. Let us show that in our case the adjoint to $D\cm_\phi$ is equal to itself using Corollary \ref{CoR:cor-by-parts}(b):
\begin{eqnarray*}
(D\cm_\phi(\psi),\xi)=2\int_M\xi\det(G_0+Hess_\OO\phi,Hess_\OO\psi) \Theta=\\
2\int_M\psi\det(G_0+Hess_\OO\phi,Hess_\OO\xi) \Theta=(D\cm_\phi(\xi),\psi).
\end{eqnarray*}
Thus $D\cm_\phi$ is self-adjoint and the kernel of its adjoint is also 1-dimensional. \qed

\item \begin{proposition}\label{P:closed-image}
Let $M^{16}$ be a compact $Spin(9)$-affine manifold. Let $k\geq 2$ be an integer, $\beta\in (0,1)$. Then the map
$$\cm\colon \cu^{k+2,\beta}\to \cv^{k,\beta}$$
is a diffeomorphism of Banach manifolds, in particular it is onto.
\end{proposition}
{\bf Proof.} $\cm$ is one-to-one by the uniqueness of the solution (Theorem \ref{T:uniqueness-thm}).


By Proposition \ref{P:open-image} it suffices to show that
$\cm$ is onto. Since $\cv^{k,\beta}$ is obviously connected (and
even convex), and since the image of $\cm$ is open by Proposition
\ref{P:open-image}, it suffices to show that the image of $\cm$ is a
closed subset of $\cv^{k,\beta}$.

Let us have a sequence of points in the image
$\cm\phi_i\overset{C^{k,\beta}}{\to} e^f\in \cv^{k,\beta}$ where
$\phi_i\in \cu^{k+2,\beta}$. By Theorems \ref{T:l-infty}, \ref{T:second-order},
\ref{thm-main}, there exist $\alp\in (0,1)$
and a constant $C$ both depending on $||f||_{C^2}$, $M$,
$G_0$, $\Theta$  such that
$||\phi_i||_{C^{2,\alp}}<C$ for $i\gg 1$. By the Arzel\`a-Ascoli
theorem choosing a subsequence we may assume that $\phi_i\to \phi$
in $C^2(M)$. Clearly $\phi\in C^{2,\alp}(M)$ and one has
$\cm(\phi)=e^f$. Explicitly one has
\begin{eqnarray}\label{E:limit-solut}
\det(G_0+Hess_\OO\phi)=e^f\det G_0.
\end{eqnarray}
Moreover, $G_0+Hess_\OO\phi\geq 0$. This, the Sylvester criterion (Proposition \ref{P:sylvester}), and
(\ref{E:limit-solut}) imply that the inequality is strict, i.e.
$G_0+Hess_\OO\phi>0$, and hence the equation (\ref{E:limit-solut}) is
elliptic with $C^\infty$-smooth coefficients on the left hand side and with
$C^{k,\beta}$ on the right hand side. Hence Lemma 17.16 from
\cite{gilbarg-trudinger} implies that $\phi\in C^{k+2,\beta}(M)$. Thus $\phi\in
\cv^{k+2,\beta}$ and $e^f=\cm(\phi)\in Im(\cm)$. \qed

\item Finally let us state the main result of the paper which is an
immediate consequence of Proposition \ref{P:closed-image}.
\begin{theorem}\label{T:MA-solution}
Let $M^{16}$ be a compact connected $Spin(9)$-affine manifold. Let $G_0$ be $C^\infty$-smooth positive definite section of the bundle $\uch$ satisfying condition (\ref{E:condition}).
Let $k\geq 2$ be an integer, $\beta\in (0,1)$, and $f\in C^{k,\beta}(M)$.

Then there is a unique constant $A>0$ such that the octonionic
Monge-Amp\`ere equation
\begin{eqnarray}\label{E:MA-in-main-thm}
\det(G_0+Hess_\OO\phi)= A e^f\det G_0
\end{eqnarray}
has a unique, up to a constant, $C^2$-smooth solution $\phi$ which
necessarily belongs to $C^{k+2,\beta}(M)$. If $f$ is $C^\infty$-smooth, then the solution $\phi$ is also $C^\infty$-smooth.
\end{theorem}

\begin{remark}
Recall that the constant $A$ in the theorem is defined by the following
condition
$$\int_M(Ae^f-1)\cdot \det G_0\cdot \Theta=0,$$
where $\Theta$ is a parallel positive $3/4$-density on $M$.
\end{remark}
\end{paragraphlist}

\section{Appendix}

In this appendix we describe the commutative algebra computer computation
necessary to prove the existence of a solution for  system \eqref{E:oct-system} for arbitrary closed smooth current $T$ satisfying conditions \eqref{E:system-equations}.
 Ehrenpreis' Theorem \ref{T:ehrenpreis} combined with this computation implies existence of a solution.

Setting $q_1:=x$ and $q_2:=y$ we have
\begin{eqnarray}
x=\sum_{k=0}^7 x_k e_k, \quad y=\sum_{k=0}^7 y_k e_k.
\end{eqnarray}
In terms of these variables system  \eqref{E:oct-system} reduces to the following system of second order differential equations
\begin{eqnarray}
\left\{
\begin{array}{l}
 \frac{\partial^2 u}{\partial \bar x \partial x}=\sum_{i=0}^7 \frac{\partial^2 u}{\partial x_i^2}=T_{\bar 1 1},\\
 \\
 \frac{\partial^2 u}{\partial \bar y \partial y}=\sum_{i=0}^7 \frac{\partial^2 u}{\partial y_i^2}=T_{\bar 2 2},\\
 \\
 \frac{\partial^2 u}{\partial \bar x \partial y}=\sum_{i,j=0}^7 e_i \bar e_j \frac{\partial^2 u}{\partial x_i\partial y_j}=T_{\bar 1 2}.
 \end{array}
 \right.
\end{eqnarray}
This system is equivalent to a following system of $10$ scalar equations
\begin{eqnarray}
&& P_1 u= T_{\bar 1 1},\quad P_2 u= T_{\bar 2 2}, \nonumber \\
&& P_{k+3} u= T_{\bar 1 2}^{k}, \quad k=0,\cdots,7,
\end{eqnarray}
where an upper index $k$ in  $T_{\bar 1 2}$ stands for a $k$-th component of octonion $T_{\bar 1 2}$ and  $P_i$ are
the second order differential operators. These operators, after substitution $\partial_{x_i} \longmapsto x_i$,    $\partial_{y_i} \longmapsto y_i,$
become  quadratic polynomials that read
\begin{eqnarray}
P_1=x\bar x=\sum_{i=0}^7 x_i^2, \quad P_2= y \bar y=\sum_{i=0}^7 y_i^2,
\end{eqnarray}
\begin{eqnarray}
\left \{\begin{array}{r}
P_3 = x_0y_0 + x_1y_1 + x_2y_2 + x_3y_3 + x_4y_4 + x_5y_5 + x_6y_6 + x_7y_7,\\
P_4 = x_1y_0 - x_0y_1 + x_3y_2 - x_2y_3 + x_5y_4 - x_4y_5 - x_7y_6 + x_6y_7,\\
P_5 = x_2y_0 - x_3y_1 - x_0y_2 + x_1y_3 + x_6y_4 + x_7y_5 - x_4y_6 - x_5y_7,\\
P_6 = x_3y_0 + x_2y_1 - x_1y_2 - x_0y_3 + x_7y_4 - x_6y_5 + x_5y_6 - x_4y_7,\\
P_7 = x_4y_0 - x_5y_1 - x_6y_2 - x_7y_3 - x_0y_4 + x_1y_5 + x_2y_6 + x_3y_7,\\
P_8 = x_5y_0 + x_4y_1 - x_7y_2 + x_6y_3 - x_1y_4 - x_0y_5 - x_3y_6 + x_2y_7,\\
P_9 = x_6y_0 + x_7y_1 + x_4y_2 - x_5y_3 - x_2y_4 + x_3y_5 - x_0y_6 - x_1y_7,\\
P_{10}= x_7y_0 - x_6y_1 + x_5y_2 + x_4y_3 - x_3y_4 - x_2y_5 + x_1y_6 - x_0y_7.
\end{array}
\right.
\end{eqnarray}
We note that polynomials $P_{k+3}$ for $k=0,\cdots,7$ are coefficient in front of $e_k$  in the expression $x\bar y.$

Setting $f_1:=  T_{\bar 1 1},$ $f_2:= T_{\bar 2 2}$ and $f_{k+3}:=T_{\bar 1 2}^{k}, \quad k=0,\cdots,7,$  we now will apply Theorem \ref{T:ehrenpreis}.
To verify that the assumptions of this theorem are met we implement the strategy described in Section \ref{S:metrics-octon} part 6.

 Setting  polynomial ring   $A:=\mathbb{Q}[x_0,\cdots,x_7,y_0,\cdots, y_7]$ (by Lemma \ref{L:commut-alg} it suffices to consider rational rather than real coefficients) we now look for
 generators of the submodule $N_\QQ$ of the free $A$-module $A^{10}$
 $$N_\QQ=\{(Q_1,\dots,Q_{10})\in A^{10}|\, \sum_{i=1}^{10} Q_iP_i=0\}.$$
   To obtain such  generators we use symbolic computations which we implement in a well known program Macaulay2.
 We give a printout of an executed program for Macaulay2 for computing generators of the submodule $N$ at the end of this appendix. The program was executed at Macaulay2Web interface located at \\https:/ /www.unimelb-macaulay2.cloud.edu.au/.
Interested reader  can download our program from  \\http://www.math.kent.edu/$\sim$gordon/MA/generator\_of\_submoduleXYb.txt\\
The output of the program produces $16$ first order polynomials as generators. These generators result in $16$ first order polynomials
which after back-substitution  $x_i\longmapsto \partial_{x_i}$, $y_i\longmapsto \partial_{y_i}$ give a system of $16$ first order differential equations that read:
\begin{eqnarray}
\left(
\begin{array}{l}\nonumber
 \partial_{x_{0}}\\
 \partial_{x_{1}}\\
 \partial_{x_{2}}\\
 \partial_{x_{3}}\\
 \partial_{x_{4}}\\
 \partial_{x_{5}}\\
 \partial_{x_{6}}\\
 \partial_{x_{7}}
 \end{array}
\right)
T_{\bar 2 2}
=\left(
\begin{array}{rrrrrrrr}
\partial_{y_{0}}& -\partial_{y_{1}}&-\partial_{y_{2}}&-\partial_{y_{3}}&-\partial_{y_{4}}&-\partial_{y_{5}}&-\partial_{y_{6}}&-\partial_{y_{7}}\\
\partial_{y_{1}}& \partial_{y_{0}}&\partial_{y_{3}}&-\partial_{y_{2}}&\partial_{y_{5}}&-\partial_{y_{4}}&-\partial_{y_{7}}&\partial_{y_{6}}\\
\partial_{y_{2}}& -\partial_{y_{3}}&\partial_{y_{0}}&\partial_{y_{1}}&\partial_{y_{6}}&\partial_{y_{7}}&-\partial_{y_{4}}&-\partial_{y_{5}}\\
\partial_{y_{3}}& \partial_{y_{2}}&-\partial_{y_{1}}&\partial_{y_{0}}&\partial_{y_{7}}&-\partial_{y_{6}}&\partial_{y_{5}}&-\partial_{y_{4}}\\
\partial_{y_{4}}&- \partial_{y_{5}}&-\partial_{y_{6}}&-\partial_{y_{7}}&\partial_{y_{0}}&\partial_{y_{1}}&\partial_{y_{2}}&\partial_{y_{3}}\\
\partial_{y_{5}}&\partial_{y_{4}}&-\partial_{y_{7}}&\partial_{y_{6}}&-\partial_{y_{1}}&\partial_{y_{0}}&-\partial_{y_{3}}&\partial_{y_{2}}\\
\partial_{y_{6}}& \partial_{y_{7}}&\partial_{y_{4}}&-\partial_{y_{5}}&-\partial_{y_{2}}&\partial_{y_{3}}&\partial_{y_{0}}&-\partial_{y_{1}}\\
\partial_{y_{7}}&- \partial_{y_{6}}&\partial_{y_{5}}&\partial_{y_{4}}&-\partial_{y_{3}}&-\partial_{y_{2}}&\partial_{y_{1}}&\partial_{y_{0}}
\end{array}
\right)
\left(
\begin{array}{c}
T_{\bar 1 2}^0\\
T_{\bar 1 2}^1\\
T_{\bar 1 2}^2\\
T_{\bar 1 2}^3\\
T_{\bar 1 2}^4\\
T_{\bar 1 2}^5\\
T_{\bar 1 2}^6\\
T_{\bar 1 2}^7
\end{array}
\right)
\end{eqnarray}

\begin{eqnarray}
\left(
\begin{array}{l}\nonumber
 \partial_{y_{0}}\\
 \partial_{y_{1}}\\
 \partial_{y_{2}}\\
 \partial_{y_{3}}\\
 \partial_{y_{4}}\\
 \partial_{y_{5}}\\
 \partial_{y_{6}}\\
 \partial_{y_{7}}
 \end{array}
\right)
T_{\bar 11}
=\left(
\begin{array}{rrrrrrrr}
\partial_{x_{0}}& \partial_{x_{1}}&\partial_{x_{2}}&\partial_{x_{3}}&\partial_{x_{4}}&\partial_{x_{5}}&\partial_{x_{6}}&\partial_{x_{7}}\\
\partial_{x_{1}}& -\partial_{x_{0}}&-\partial_{x_{3}}&\partial_{x_{2}}&-\partial_{x_{5}}&\partial_{x_{4}}&\partial_{x_{7}}&-\partial_{x_{6}}\\
\partial_{x_{2}}& \partial_{x_{3}}&-\partial_{x_{0}}&-\partial_{x_{1}}&-\partial_{x_{6}}&-\partial_{x_{7}}&\partial_{x_{4}}&\partial_{x_{5}}\\
\partial_{x_{3}}&- \partial_{x_{2}}&\partial_{x_{1}}&-\partial_{x_{0}}&-\partial_{x_{7}}&\partial_{x_{6}}&-\partial_{x_{5}}&\partial_{x_{4}}\\
\partial_{x_{4}}& \partial_{x_{5}}&\partial_{x_{6}}&\partial_{x_{7}}&-\partial_{x_{0}}&-\partial_{x_{1}}&-\partial_{x_{2}}&-\partial_{x_{3}}\\
\partial_{x_{5}}&- \partial_{x_{4}}&\partial_{x_{7}}&-\partial_{x_{6}}&\partial_{x_{1}}&-\partial_{x_{0}}&\partial_{x_{3}}&-\partial_{x_{2}}\\
\partial_{x_{6}}&- \partial_{x_{7}}&-\partial_{x_{4}}&\partial_{x_{5}}&\partial_{x_{2}}&-\partial_{x_{3}}&-\partial_{x_{0}}&\partial_{x_{1}}\\
\partial_{x_{7}}&\partial_{x_{6}}&-\partial_{x_{5}}&-\partial_{x_{4}}&\partial_{x_{3}}&\partial_{x_{2}}&-\partial_{x_{1}}&-\partial_{x_{0}}
\end{array}
\right)
\left(
\begin{array}{c}
T_{\bar 1 2}^0\\
T_{\bar 1 2}^1\\
T_{\bar 1 2}^2\\
T_{\bar 1 2}^3\\
T_{\bar 1 2}^4\\
T_{\bar 1 2}^5\\
T_{\bar 1 2}^6\\
T_{\bar 1 2}^7
\end{array}
\right)
\end{eqnarray}

Noting that $T_{\bar 2 1}=\bar T_{\bar 1 2}$ one can directly  verify  that the systems of equations above is a scalar form of the closed smooth current condition \eqref{E:system-equations}.

\begin{landscape}
\begin{verbatim}
Macaulay2, version 1.20.0.1
with packages: ConwayPolynomials, Elimination, IntegralClosure, InverseSystems, Isomorphism,
LLLBases, MinimalPrimes, OnlineLookup, PrimaryDecomposition, ReesAlgebra, Saturation, TangentCone

i1 : A=QQ[x_0..y_7];

i2 : P_1 = x_0*x_0 + x_1*x_1 + x_2*x_2 + x_3*x_3 + x_4*x_4 + x_5*x_5 + x_6*x_6 + x_7*x_7;

i3 : P_2 = y_0*y_0 + y_1*y_1 + y_2*y_2 + y_3*y_3 + y_4*y_4 + y_5*y_5 + y_6*y_6 + y_7*y_7;

i4 : P_3 = x_0*y_0 + x_1*y_1 + x_2*y_2 + x_3*y_3 + x_4*y_4 + x_5*y_5 + x_6*y_6 + x_7*y_7;

i5 : P_4 = x_1*y_0 - x_0*y_1 + x_3*y_2 - x_2*y_3 + x_5*y_4 - x_4*y_5 - x_7*y_6 + x_6*y_7;

i6 : P_5 = x_2*y_0 - x_3*y_1 - x_0*y_2 + x_1*y_3 + x_6*y_4 + x_7*y_5 - x_4*y_6 - x_5*y_7;

i7 : P_6 = x_3*y_0 + x_2*y_1 - x_1*y_2 - x_0*y_3 + x_7*y_4 - x_6*y_5 + x_5*y_6 - x_4*y_7;

i8 : P_7 = x_4*y_0 - x_5*y_1 - x_6*y_2 - x_7*y_3 - x_0*y_4 + x_1*y_5 + x_2*y_6 + x_3*y_7;

i9 : P_8 = x_5*y_0 + x_4*y_1 - x_7*y_2 + x_6*y_3 - x_1*y_4 - x_0*y_5 - x_3*y_6 + x_2*y_7;

i10 : P_9 = x_6*y_0 + x_7*y_1 + x_4*y_2 - x_5*y_3 - x_2*y_4 + x_3*y_5 - x_0*y_6 - x_1*y_7;

i11 : P_10= x_7*y_0 - x_6*y_1 + x_5*y_2 + x_4*y_3 - x_3*y_4 - x_2*y_5 + x_1*y_6 - x_0*y_7;

i12 : N=matrix{{P_1,P_2,P_3,P_4,P_5,P_6,P_7,P_8,P_9,P_10}};
\end{verbatim}
\smallskip
\noindent  {\tt o12 :}  \hspace{.1cm}{\tt Matrix}  $A^1 \longleftarrow A^{10}$
\begin{verbatim}
i13 : kernel N
\end{verbatim}
$\texttt{o13~=~image}\ \begin{array}{l}\left\{2\right\}\vphantom{00000000-y_{0}y_{1}y_{2}y_{3}y_{4}y_{5}y_{6}y_{7}}\\\left\{2\right\}\vphantom{x_{0}x_{1}x_{2}x_{3}x_{4}x_{5}x_{6}x_{7}00000000}\\\left\{2\right\}\vphantom{-y_{0}-y_{1}-y_{2}-y_{3}-y_{4}-y_{5}-y_{6}-y_{7}x_{0}-x_{1}-x_{2}-x_{3}-x_{4}-x_{5}-x_{6}-x_{7}}\\\left\{2\right\}\vphantom{y_{1}-y_{0}y_{3}-y_{2}y_{5}-y_{4}-y_{7}y_{6}x_{1}x_{0}-x_{3}x_{2}-x_{5}x_{4}x_{7}-x_{6}}\\\left\{2\right\}\vphantom{y_{2}-y_{3}-y_{0}y_{1}y_{6}y_{7}-y_{4}-y_{5}x_{2}x_{3}x_{0}-x_{1}-x_{6}-x_{7}x_{4}x_{5}}\\\left\{2\right\}\vphantom{y_{3}y_{2}-y_{1}-y_{0}y_{7}-y_{6}y_{5}-y_{4}x_{3}-x_{2}x_{1}x_{0}-x_{7}x_{6}-x_{5}x_{4}}\\\left\{2\right\}\vphantom{y_{4}-y_{5}-y_{6}-y_{7}-y_{0}y_{1}y_{2}y_{3}x_{4}x_{5}x_{6}x_{7}x_{0}-x_{1}-x_{2}-x_{3}}\\\left\{2\right\}\vphantom{y_{5}y_{4}-y_{7}y_{6}-y_{1}-y_{0}-y_{3}y_{2}x_{5}-x_{4}x_{7}-x_{6}x_{1}x_{0}x_{3}-x_{2}}\\\left\{2\right\}\vphantom{y_{6}y_{7}y_{4}-y_{5}-y_{2}y_{3}-y_{0}-y_{1}x_{6}-x_{7}-x_{4}x_{5}x_{2}-x_{3}x_{0}x_{1}}\\\left\{2\right\}\vphantom{y_{7}-y_{6}y_{5}y_{4}-y_{3}-y_{2}y_{1}-y_{0}x_{7}x_{6}-x_{5}-x_{4}x_{3}x_{2}-x_{1}x_{0}}\end{array}\left(\!\begin{array}{cccccccccccccccc}
\vphantom{\left\{2\right\}}0&0&0&0&0&0&0&0&-y_{0}&y_{1}&y_{2}&y_{3}&y_{4}&y_{5}&y_{6}&y_{7}\\
\vphantom{\left\{2\right\}}x_{0}&x_{1}&x_{2}&x_{3}&x_{4}&x_{5}&x_{6}&x_{7}&0&0&0&0&0&0&0&0\\
\vphantom{\left\{2\right\}}-y_{0}&-y_{1}&-y_{2}&-y_{3}&-y_{4}&-y_{5}&-y_{6}&-y_{7}&x_{0}&-x_{1}&-x_{2}&-x_{3}&-x_{4}&-x_{5}&-x_{6}&-x_{7}\\
\vphantom{\left\{2\right\}}y_{1}&-y_{0}&y_{3}&-y_{2}&y_{5}&-y_{4}&-y_{7}&y_{6}&x_{1}&x_{0}&-x_{3}&x_{2}&-x_{5}&x_{4}&x_{7}&-x_{6}\\
\vphantom{\left\{2\right\}}y_{2}&-y_{3}&-y_{0}&y_{1}&y_{6}&y_{7}&-y_{4}&-y_{5}&x_{2}&x_{3}&x_{0}&-x_{1}&-x_{6}&-x_{7}&x_{4}&x_{5}\\
\vphantom{\left\{2\right\}}y_{3}&y_{2}&-y_{1}&-y_{0}&y_{7}&-y_{6}&y_{5}&-y_{4}&x_{3}&-x_{2}&x_{1}&x_{0}&-x_{7}&x_{6}&-x_{5}&x_{4}\\
\vphantom{\left\{2\right\}}y_{4}&-y_{5}&-y_{6}&-y_{7}&-y_{0}&y_{1}&y_{2}&y_{3}&x_{4}&x_{5}&x_{6}&x_{7}&x_{0}&-x_{1}&-x_{2}&-x_{3}\\
\vphantom{\left\{2\right\}}y_{5}&y_{4}&-y_{7}&y_{6}&-y_{1}&-y_{0}&-y_{3}&y_{2}&x_{5}&-x_{4}&x_{7}&-x_{6}&x_{1}&x_{0}&x_{3}&-x_{2}\\
\vphantom{\left\{2\right\}}y_{6}&y_{7}&y_{4}&-y_{5}&-y_{2}&y_{3}&-y_{0}&-y_{1}&x_{6}&-x_{7}&-x_{4}&x_{5}&x_{2}&-x_{3}&x_{0}&x_{1}\\
\vphantom{\left\{2\right\}}y_{7}&-y_{6}&y_{5}&y_{4}&-y_{3}&-y_{2}&y_{1}&-y_{0}&x_{7}&x_{6}&-x_{5}&-x_{4}&x_{3}&x_{2}&-x_{1}&x_{0}
\end{array}\!\right)$

\bigskip\bigskip

\noindent {\tt o13 :}  \hspace{.1cm}{\tt  $A$-module, submodule of $A^{10}$ }

  \end{landscape}


\end{document}